\pdfoutput=1

\documentclass[12pt,a4paper]{article}
\usepackage{graphicx} 
\usepackage{eqsection,indent}
\usepackage{array, url}
\usepackage{amsmath, amssymb, amsfonts,bm}
\usepackage{color}
\usepackage{subeqnarray}

\usepackage{tikz,subfigure,float}

\usetikzlibrary{positioning}

\usepackage{changepage}

\title{\vspace*{-1cm}Thron-type continued fractions (T-fractions) \\
       for some classes of increasing trees
         \\[3mm]}

\date{12 December 2024}

\author{
      \hspace*{-1cm}
      {\large Veronica Bitonti${}^{1,2}$, Bishal Deb${}^{3}$,
            and Alan D.~Sokal${}^{1,4}$}
   \\[5mm]
     \hspace*{-1.25cm}
      \normalsize
           ${}^1$Department of Mathematics, University College London,
                    London WC1E 6BT, UK   \\[2mm]
     \hspace*{-3.8cm}
      \normalsize
           ${}^2$Trinity College, University of Cambridge,
                    Cambridge CB2 1TQ, UK   \\[2mm]
     \hspace*{-6.3cm}
      \normalsize
           ${}^3$Sorbonne Universit\'e and Universit\'e Paris Cit\'e, CNRS,
       \\
      \normalsize 
     \hspace*{-1.6cm}
          Laboratoire de Probabilit\'es,
           Statistique et Mod\'elisation, Paris 75005, FRANCE\\[2mm]
     \normalsize
     \hspace*{-2.9cm}
           ${}^4$Department of Physics, New York University,
                    New York, NY 10003, USA
       \\[4mm]
     \hspace*{-1cm}
     {\tt vb435@cam.ac.uk},
     {\tt bishal@gonitsora.com},
     {\tt sokal@nyu.edu}  \\[4mm]
}

\begin{document}

\maketitle
\thispagestyle{empty}

\begin{abstract}
We introduce some classes of increasing labeled and multilabeled trees,
and we show that these trees provide combinatorial interpretations
for certain Thron-type continued fractions with coefficients that are 
quasi-affine of period~$2$.
Our proofs are based on bijections from trees to
labeled Motzkin or Schr\"oder paths;
these bijections extend the well-known bijection of
Fran{\c{c}}on--Viennot (1979)
interpreted in terms of increasing binary trees.
This work can also be viewed as a sequel to the recent work 
of Elvey Price and Sokal~(2020), where they provide combinatorial
interpretations for Thron-type continued fractions 
with coefficients that are affine.
Towards the end of the paper, we conjecture an
equidistribution of vincular patterns on permutations.
\end{abstract}

\medskip
\noindent
{\bf Key Words:}
Continued fraction, S-fraction, J-fraction, T-fraction,
increasing tree, multilabeled tree,
Motzkin path, Schr\"oder path,
exponential Riordan array, permutation, vincular pattern.

\medskip
\noindent
{\bf Mathematics Subject Classification (MSC 2020) codes:}
05A19 (Primary);
05A05, 05A15, 05A30, 05C05, 05C30, 30B70 (Secondary).

\vspace*{1cm}

\newtheorem{theorem}{Theorem}[section]
\newtheorem{proposition}[theorem]{Proposition}
\newtheorem{lemma}[theorem]{Lemma}
\newtheorem{corollary}[theorem]{Corollary}
\newtheorem{definition}[theorem]{Definition}
\newtheorem{conjecture}[theorem]{Conjecture}
\newtheorem{question}[theorem]{Question}
\newtheorem{problem}[theorem]{Problem}
\newtheorem{openproblem}[theorem]{Open Problem}
\newtheorem{example}[theorem]{Example}
\newtheorem{remark}[theorem]{Remark}

\renewcommand{\theenumi}{\alph{enumi}}
\renewcommand{\labelenumi}{(\theenumi)}
\def\eop{\hbox{\kern1pt\vrule height6pt width4pt
depth1pt\kern1pt}\medskip}
\def\prf{\par\noindent{\bf Proof.\enspace}\rm}
\def\rmk{\par\medskip\noindent{\bf Remark\enspace}\rm}

\newcommand{\textbfit}[1]{\textbf{\textit{#1}}}

\newcommand{\be}{\begin{equation}}
\newcommand{\ee}{\end{equation}}
\newcommand{\<}{\langle}
\renewcommand{\>}{\rangle}
\newcommand{\widebar}{\overline}
\def\reff#1{(\protect\ref{#1})}
\def\spose#1{\hbox to 0pt{#1\hss}}
\def\ltapprox{\mathrel{\spose{\lower 3pt\hbox{$\mathchar"218$}}
    \raise 2.0pt\hbox{$\mathchar"13C$}}}
\def\gtapprox{\mathrel{\spose{\lower 3pt\hbox{$\mathchar"218$}}
    \raise 2.0pt\hbox{$\mathchar"13E$}}}
\def\textprime{${}^\prime$}
\def\proof{\par\medskip\noindent{\sc Proof.\ }}
\def\firstproof{\par\medskip\noindent{\sc First Proof.\ }}
\def\secondproof{\par\medskip\noindent{\sc Second Proof.\ }}
\def\alternateproof{\par\medskip\noindent{\sc Alternate Proof.\ }}
\def\algebraicproof{\par\medskip\noindent{\sc Algebraic Proof.\ }}
\def\graphicalproof{\par\medskip\noindent{\sc Graphical Proof.\ }}
\def\combinatorialproof{\par\medskip\noindent{\sc Combinatorial Proof.\ }}
\def\proofof#1{\bigskip\noindent{\sc Proof of #1.\ }}
\def\firstproofof#1{\bigskip\noindent{\sc First Proof of #1.\ }}
\def\secondproofof#1{\bigskip\noindent{\sc Second Proof of #1.\ }}
\def\thirdproofof#1{\bigskip\noindent{\sc Third Proof of #1.\ }}
\def\algebraicproofof#1{\bigskip\noindent{\sc Algebraic Proof of #1.\ }}
\def\combinatorialproofof#1{\bigskip\noindent{\sc Combinatorial Proof of #1.\ }}
\def\completionofproofof#1{\bigskip\noindent{\sc Completion of the Proof of #1.\ }}
\def\sketchofproof{\par\medskip\noindent{\sc Sketch of proof.\ }}
 \def\qed{\hbox{\hskip 6pt\vrule width6pt height7pt depth1pt \hskip1pt}\bigskip}
\renewcommand{\qed}{ $\square$ \bigskip}
\newcommand{\myendremark}{ $\blacksquare$ \bigskip}
\def\half{ {1 \over 2} }
\def\third{ {1 \over 3} }
\def\twothird{ {2 \over 3} }
\def\smfrac#1#2{{\textstyle{#1\over #2}}}
\def\smhalf{ {\smfrac{1}{2}} }
\newcommand{\real}{\mathop{\rm Re}\nolimits}
\renewcommand{\Re}{\mathop{\rm Re}\nolimits}
\newcommand{\imag}{\mathop{\rm Im}\nolimits}
\renewcommand{\Im}{\mathop{\rm Im}\nolimits}
\newcommand{\sgn}{\mathop{\rm sgn}\nolimits}
\newcommand{\tr}{\mathop{\rm tr}\nolimits}
\newcommand{\supp}{\mathop{\rm supp}\nolimits}
\newcommand{\disc}{\mathop{\rm disc}\nolimits}
\newcommand{\diag}{\mathop{\rm diag}\nolimits}
\newcommand{\tridiag}{\mathop{\rm tridiag}\nolimits}
\newcommand{\AZ}{\mathop{\rm AZ}\nolimits}
\newcommand{\EAZ}{\mathop{\rm EAZ}\nolimits}
\newcommand{\NC}{\mathop{\rm NC}\nolimits}
\newcommand{\PF}{{\rm PF}}
\newcommand{\rk}{\mathop{\rm rk}\nolimits}
\newcommand{\perm}{\mathop{\rm perm}\nolimits}
\def\hboxscript#1{ {\hbox{\scriptsize\em #1}} }
\renewcommand{\emptyset}{\varnothing}
\newcommand{\eqdef}{\stackrel{\rm def}{=}}

\newcommand{\restrict}{\upharpoonright}

\newcommand{\compinv}{{\langle -1 \rangle}}   

\newcommand{\scra}{{\mathcal{A}}}
\newcommand{\scrb}{{\mathcal{B}}}
\newcommand{\scrc}{{\mathcal{C}}}
\newcommand{\scrd}{{\mathcal{D}}}
\newcommand{\scrdtilde}{{\widetilde{\mathcal{D}}}}
\newcommand{\scre}{{\mathcal{E}}}
\newcommand{\scrf}{{\mathcal{F}}}
\newcommand{\scrg}{{\mathcal{G}}}
\newcommand{\scrh}{{\mathcal{H}}}
\newcommand{\scri}{{\mathcal{I}}}
\newcommand{\scrj}{{\mathcal{J}}}
\newcommand{\scrk}{{\mathcal{K}}}
\newcommand{\scrl}{{\mathcal{L}}}
\newcommand{\scrlbar}{{\overline{\mathcal{L}}}}
\newcommand{\scrm}{{\mathcal{M}}}
\newcommand{\scrn}{{\mathcal{N}}}
\newcommand{\scro}{{\mathcal{O}}}
\newcommand\scroo{
  \mathchoice
    {{\scriptstyle\mathcal{O}}}
    {{\scriptstyle\mathcal{O}}}
    {{\scriptscriptstyle\mathcal{O}}}
    {\scalebox{0.6}{$\scriptscriptstyle\mathcal{O}$}}
  }
\newcommand{\scrp}{{\mathcal{P}}}
\newcommand{\scrq}{{\mathcal{Q}}}
\newcommand{\scrr}{{\mathcal{R}}}
\newcommand{\scrs}{{\mathcal{S}}}
\newcommand{\scrt}{{\mathcal{T}}}
\newcommand{\scrv}{{\mathcal{V}}}
\newcommand{\scrw}{{\mathcal{W}}}
\newcommand{\scrz}{{\mathcal{Z}}}
\newcommand{\SP}{{\mathcal{SP}}}
\newcommand{\ST}{{\mathcal{ST}}}

\newcommand{\bfa}{{\mathbf{a}}}
\newcommand{\bfA}{{\mathbf{A}}}
\newcommand{\bfb}{{\mathbf{b}}}
\newcommand{\bfc}{{\mathbf{c}}}
\newcommand{\bfd}{{\mathbf{d}}}
\newcommand{\bfe}{{\mathbf{e}}}
\newcommand{\bfff}{{\mathbf{f}}}
\newcommand{\bfh}{{\mathbf{h}}}
\newcommand{\bfj}{{\mathbf{j}}}
\newcommand{\bfi}{{\mathbf{i}}}
\newcommand{\bfk}{{\mathbf{k}}}
\newcommand{\bfl}{{\mathbf{l}}}
\newcommand{\bfL}{{\mathbf{L}}}
\newcommand{\bfm}{{\mathbf{m}}}
\newcommand{\bfn}{{\mathbf{n}}}
\newcommand{\bfp}{{\mathbf{p}}}
\newcommand{\bfr}{{\mathbf{r}}}
\newcommand{\bft}{{\mathbf{t}}}
\newcommand{\bfu}{{\mathbf{u}}}
\newcommand{\bfv}{{\mathbf{v}}}
\newcommand{\bfw}{{\mathbf{w}}}
\newcommand{\bfx}{{\mathbf{x}}}
\newcommand{\bfX}{{\mathbf{X}}}
\newcommand{\bfy}{{\mathbf{y}}}
\newcommand{\bfz}{{\mathbf{z}}}
\renewcommand{\k}{{\mathbf{k}}}
\newcommand{\n}{{\mathbf{n}}}
\newcommand{\vv}{{\mathbf{v}}}
\newcommand{\w}{{\mathbf{w}}}
\newcommand{\x}{{\mathbf{x}}}
\newcommand{\y}{{\mathbf{y}}}
\newcommand{\cc}{{\mathbf{c}}}
\newcommand{\zero}{{\mathbf{0}}}
\newcommand{\one}{{\mathbf{1}}}
\newcommand{\bmm}{{\mathbf{m}}}

\newcommand{\ba}{{\bm{a}}}
\newcommand{\bz}{{\bm{z}}}

\newcommand{\balpha}{{\bm{\alpha}}}
\newcommand{\balphapre}{{\bm{\alpha}^{\rm pre}}}
\newcommand{\bbeta}{{\bm{\beta}}}
\newcommand{\bgamma}{{\bm{\gamma}}}
\newcommand{\bdelta}{{\bm{\delta}}}
\newcommand{\bkappa}{{\bm{\kappa}}}
\newcommand{\blambda}{{\bm{\lambda}}}
\newcommand{\bmu}{{\bm{\mu}}}
\newcommand{\bnu}{{\bm{\nu}}}
\newcommand{\bphi}{{\bm{\phi}}}

\newcommand{\sfa}{{{\sf a}}}
\newcommand{\sfahat}{{\widehat{\sf a}}}
\newcommand{\sfb}{{{\sf b}}}
\newcommand{\sfbhat}{{\widehat{\sf b}}}
\newcommand{\sfc}{{{\sf c}}}
\newcommand{\sfd}{{{\sf d}}}
\newcommand{\sfe}{{{\sf e}}}
\newcommand{\sff}{{{\sf f}}}
\newcommand{\sfg}{{{\sf g}}}
\newcommand{\sfh}{{{\sf h}}}
\newcommand{\sfi}{{{\sf i}}}
\newcommand{\sfL}{{\sf L}}
\newcommand{\bsfa}{{\mbox{\textsf{\textbf{a}}}}}
\newcommand{\bsfahat}{{\widehat{\mbox{\textsf{\textbf{a}}}}}}
\newcommand{\bsfb}{{\mbox{\textsf{\textbf{b}}}}}
\newcommand{\bsfbhat}{{\widehat{\mbox{\textsf{\textbf{b}}}}}}
\newcommand{\bsfc}{{\mbox{\textsf{\textbf{c}}}}}
\newcommand{\bsfd}{{\mbox{\textsf{\textbf{d}}}}}
\newcommand{\bsfe}{{\mbox{\textsf{\textbf{e}}}}}
\newcommand{\bsff}{{\mbox{\textsf{\textbf{f}}}}}
\newcommand{\bsfg}{{\mbox{\textsf{\textbf{g}}}}}
\newcommand{\bsfh}{{\mbox{\textsf{\textbf{h}}}}}
\newcommand{\bsfi}{{\mbox{\textsf{\textbf{i}}}}}

\newcommand{\bfscra}{{\bm{\mathcal{A}}}}
\newcommand{\bfscrb}{{\bm{\mathcal{B}}}}
\newcommand{\bfscrc}{{\bm{\mathcal{C}}}}
\newcommand{\bfscrap}{{\bm{\mathcal{A}'}}}
\newcommand{\bfscrbp}{{\bm{\mathcal{B}'}}}
\newcommand{\bfscrcp}{{\bm{\mathcal{C}'}}}
\newcommand{\bfscrapp}{{\bm{\mathcal{A}''}}}
\newcommand{\bfscrbpp}{{\bm{\mathcal{B}''}}}
\newcommand{\bfscrcpp}{{\bm{\mathcal{C}''}}}

\newcommand{\C}{{\mathbb C}}
\newcommand{\D}{{\mathbb D}}
\newcommand{\Z}{{\mathbb Z}}
\newcommand{\N}{{\mathbb N}}
\newcommand{\Q}{{\mathbb Q}}
\newcommand{\PP}{{\mathbb P}}
\newcommand{\R}{{\mathbb R}}
\newcommand{\RR}{{\mathbb R}}
\newcommand{\E}{{\mathbb E}}

\newcommand{\scrgg}{{\mathscr{g}}}  

\newcommand{\bzero}{{\bm{0}}}
\newcommand{\bone}{{\bm{1}}}

\newcommand{\scrss}{{\mathscr{s}}}  

\newcommand{\val}{{\rm val}}
\newcommand{\pk}{{\rm pk}}
\newcommand{\dasc}{{\rm dasc}}
\newcommand{\ddes}{{\rm ddes}}

\newcommand{\Val}{{\rm Val}}
\newcommand{\Pk}{{\rm Pk}}
\newcommand{\Dasc}{{\rm Dasc}}
\newcommand{\Ddes}{{\rm Ddes}}

\newcommand{\lev}{{\rm lev}}
\newcommand{\wt}{{\rm wt}}

\newcommand{\nid}{{\rm nid}}
\newcommand{\croix}{{\rm croix}}

\newcommand{\lodd}{{\rm lodd}}
\newcommand{\leven}{{\rm leven}}

\newcommand{\omegahat}{\widehat{\omega}}
\newcommand{\omegatilde}{\widetilde{\omega}}

\newcommand{\RT}{\mathcal{RT}}
\newcommand{\IRT}{\mathcal{IRT}}

\newcommand*{\Scale}[2][4]{\scalebox{#1}{$#2$}}

\newcommand*{\Scaletext}[2][4]{\scalebox{#1}{#2}} 

\clearpage 

\tableofcontents

\clearpage

\section{Introduction}

Let $(a_n)_{n\geq 0}$ be a sequence of combinatorial numbers or polynomials
with \hbox{$a_0=1$}.
In~this paper we are interested in expressing the ordinary generating
function $\sum_{n=0}^\infty a_n t^n$ as a continued fraction of Thron-type
(T-fraction for short):
\be
\sum_{n=0}^\infty a_n t^n 
\;=\;
\cfrac{1}{1 - \delta_1 t - \cfrac{\alpha_1 t}{1 - \delta_2 t - \cfrac{\alpha_2 t}{ 1 - \delta_3 t - \cfrac{\alpha_3 t}{1 - \cdots}}}}
   \;\;.
 \label{eq.Tfrac.def}
\ee
(Both sides of this expression are to be interpreted as formal power series
in the indeterminate $t$.)

The study of T-fractions, especially those in which all the coefficients
$\delta_i$ and $\alpha_i$ are nonzero, is comparatively rare
in the combinatorial literature.
The most commonly studied continued fractions are 
those of Stieltjes-type (S-fraction),
\be
\sum_{n=0}^\infty a_n t^n 
\;=\;
\cfrac{1}{1 - \cfrac{\alpha_1 t}{1 - \cfrac{\alpha_2 t}{ 1 - \cfrac{\alpha_3 t}{1 - \cdots}}}}
   \;\;,
 \label{def.Stype}
\ee
and Jacobi-type (J-fraction),
\be
\sum_{n=0}^\infty a_n t^n 
\;=\;
\cfrac{1}{1 - \gamma_0 t - \cfrac{\beta_1 t^2}{1 - \gamma_2 t - \cfrac{1 - \beta_2 t^2}{ 1 - \gamma_3 t - \cfrac{\beta_3 t^2}{1 - \cdots}}}}
   \;\;.
 \label{def.Jtype}
\ee
Clearly, T-fractions are a generalization of the S-fractions,
and reduce to them when $\delta_i=0$ for all $i$.
A nontrivial example of a T-fraction along with combinatorial interpretations
was obtained recently by Elvey Price and one of us
\cite{Elvey-Price-Sokal_wardpoly}:
we considered the T-fraction with coefficients that are affine in $n$,
\begin{subeqnarray}
\alpha_n  &=& x+ (n-1)u \\
\delta_n &=& z + (n-1)(w'+w'')
   \label{eq.alphadelta.elvey}
\end{subeqnarray}
and we showed \cite[Theorem~1.2]{Elvey-Price-Sokal_wardpoly}
that the Taylor coefficients $a_n$
are multivariate generalizations of the Ward polynomials 
that enumerate {\em super-augmented perfect matchings}\/
or {\em phylogenetic trees}\/
with respect to suitable statistics.
Note that specializing to $z=w'=w''=0$ gives the 
S-fraction for perfect matchings with a weight~$x$ for each 
record and a weight~$u$ for each cycle-peak non-record
\cite[Theorem~4.1]{Sokal-Zeng_masterpoly}.
The most general result in \cite{Elvey-Price-Sokal_wardpoly}
--- the so-called ``master T-fraction''
\cite[Theorem~2.1]{Elvey-Price-Sokal_wardpoly} ---
enumerates super-augmented perfect matchings
with respect to an infinite number of statistics.

Our initial goal in this project
was to generalize \reff{eq.alphadelta.elvey}
to the case in which the coefficients $\alpha_n$ and $\delta_n$,
rather than being affine in $n$,
are instead ``quasi-affine of period~2'':
\begin{subeqnarray}
\alpha_{2k-1}  &=& x + (k-1)u \\
\alpha_{2k}    &=& y + (k-1)v \\
\delta_{2k-1}  &=& a + (k-1)c \\
\delta_{2k}    &=& b + (k-1)d
   \label{eq.alphadelta.quasiaffine}
\end{subeqnarray}
Note that specializing to $a=b=c=d=0$ in this T-fraction
gives the S-fraction enumerating permutations
with respect to exclusive records, antirecords and excedances
\cite[Theorem~2.1]{Sokal-Zeng_masterpoly}.
And this is, in turn, a special case of an S-fraction or J-fraction
enumerating permutations with respect to an infinite number of statistics
\cite[Theorems 2.9 and 2.11]{Sokal-Zeng_masterpoly}.

When all eight parameters in \reff{eq.alphadelta.quasiaffine} equal~1,
the sequence $(a_n)_{n \ge 0}$ is
\be
   (a_n)_{n \ge 0}  \;=\;
   1, 2, 6, 24, 124, 800, 6208, 56240, 582272, 6781888, 87769632, \ldots\,,
\ee
which is not in the OEIS \cite{OEIS} and for which we do not have
any natural combinatorial interpretation.\footnote{
   We have an ``unnatural'' combinatorial interpretation
   by taking $\bsfa=\bsfb=\bsfc=\bsfd=\bsff=\one$
   and $\sfe_k = k+1$ in Theorem~\ref{thm.sbt.master}.
   This corresponds to increasing interval-labeled restricted ternary trees 
   with the non-minimal labels in vertices being ``multicolored''.
}
However, for certain special cases we are able to find
a natural combinatorial interpretation.
For instance, when $c=0$ and the other seven variables equal~1, we have
\be
   (a_n)_{n \ge 0}  \;=\;
   1, 2, 6, 23, 109, 632, 4390, 35621, 330545, 3451774, 40059838, \ldots\,;
  \label{eq.c=0}
\ee
and although this sequence is not currently in the OEIS,
we will show that it enumerates
increasing interval-labeled restricted ternary trees
--- a class that will be defined later in this introduction.
More generally, with the constraints $c=0$, $x=u$, $y=v$, $b=d$
we are able to find a combinatorial interpretation of
the polynomials in four variables generated by
\be
   \alpha_{2k-1} \:=\: kx \,,\quad
   \alpha_{2k} \:=\: ky \,,\quad
   \delta_{2k-1} \:=\: a \,,\quad
   \delta_{2k} \:=\: kb
\ee
as enumerating increasing interval-labeled restricted ternary trees
with respect to some statistics counting {\em node types} and {\em label surplus.}


Our combinatorial interpretations will thus involve
several classes of increasing trees,
i.e.~labeled rooted trees in which the labels increase from parent to child.
We will prove each J-fraction or T-fraction by constructing a bijection
from the given class of increasing trees
to a suitable class of labeled Motzkin or Schr\"{o}der paths.
These bijections will generalize
the classical Fran\c{c}on--Viennot \cite{Francon_79} bijection
from permutations to labeled Motzkin paths,
reformulated \cite{Flajolet_80,Goulden_83}
as a bijection from increasing binary trees to labeled Motzkin paths.
Some of our tree constructions have been studied previously
by Kuba and Panholzer \cite{Kuba_16},
albeit not in the context of continued fractions.
%

In the remainder of this introduction
we introduce our families of increasing trees
in increasing levels of generality
(Section~\ref{subsec.intro.comb.interpret})
and then provide an outline of the rest of the paper
(Section~\ref{subsec.intro.outline}).



\subsection{Combinatorial models for T-fractions}
 \label{subsec.intro.comb.interpret}

Our main combinatorial objects, as the title suggests,
are several families of increasing trees.
We list them here:

\medskip

{\bf Increasing binary trees.}
Our first family is the well-known set of all
\textbfit{increasing binary trees}
on the vertex set $[n] \eqdef \{1,\ldots,n\}$.
That is, there is a binary rooted tree with $n$ vertices,
in which the vertices carry distinct labels from the label set $[n]$;
furthermore, the label of a child is always greater than
the label of its parent.
(In particular, the root must have label~1,
and the vertex with label~$n$ is necessarily a leaf.)
It is well known \cite[p.~45]{Stanley_12}
that the cardinality of such trees is $n!$.
The sequence $(n!)_{n\geq 0}$ has the well-known S-fraction
with $\alpha$'s given by $1,1,2,2,3,3,4,4,5,5,\ldots$,
due to Euler \cite[section~21]{Euler_1760}.\footnote{
   The paper \cite{Euler_1760},
   which is E247 in Enestr\"om's \cite{Enestrom_13} catalogue,
   was probably written circa 1746;
   it~was presented to the St.~Petersburg Academy in 1753,
   and published in 1760.
}
This is~\reff{eq.alphadelta.quasiaffine}
with $x=y=u=v=1$ and $a=b=c=d=0$.
We denote by $\mathcal{B}_n$
the set of all increasing binary trees on the vertex set $[n]$.

\begin{figure}[t]
\centering
\newcommand{\nodea}{\node[draw,circle] (a) {$1$}
;}\newcommand{\nodeb}{\node[draw,circle] (b) {$2$}
;}\newcommand{\nodec}{\node[draw,circle] (c) {$3$}
;}\newcommand{\noded}{\node[draw,circle] (d) {$4$}
;}\newcommand{\nodee}{\node[draw,circle] (e) {$5$}
;}\newcommand{\nodef}{\node[draw,circle] (f) {$6$}
;}\newcommand{\nodeg}{\node[draw,circle] (g) {$7$}
;}\newcommand{\nodeh}{\node[draw,circle] (h) {$8$}
;}
\begin{tikzpicture}
\matrix[column sep=.1cm, row sep=0.5cm,ampersand replacement=\&]{
\& \& \nodea \& \& \\
\& \nodec \& \& \nodeb \& \\
\nodee \& \& \nodef \& \& \noded \\
\& \nodeg \& \& \nodeh \& \\
};
	\draw[ultra thick] (a) -- (b) node {};
        \draw[ultra thick] (a) -- (c) node {};
        \draw[ultra thick] (b) -- (d) node {};
        \draw[ultra thick] (b) -- (f) node {};
        \draw[ultra thick] (c) -- (e) node {};
        \draw[ultra thick] (d) -- (h) node {};
        \draw[ultra thick] (e) -- (g) node {};
\node[left = -0.5mm of e] () {{\footnotesize $1$}};
	\node[left = -0.5mm of g] () {{\footnotesize $2$}};
	\node[left = -0.5mm of c] () {{\footnotesize $3$}};
	\node[left = -0.5mm of a] () {{\footnotesize $4$}};
	\node[left = -0.5mm of f] () {{\footnotesize $5$}};
	\node[left = -0.5mm of b] () {{\footnotesize $6$}};
	\node[left = -0.5mm of h] () {{\footnotesize $7$}};
	\node[left = -0.5mm of d] () {{\footnotesize $8$}};
\end{tikzpicture}
 \vspace*{-2mm}
\caption{An example of an increasing binary tree on the vertex set $[8]$.
	To the left of each vertex, its order as per the inorder traversal
        is recorded.
}
 \vspace*{1cm}
\label{fig.bt}
\end{figure}
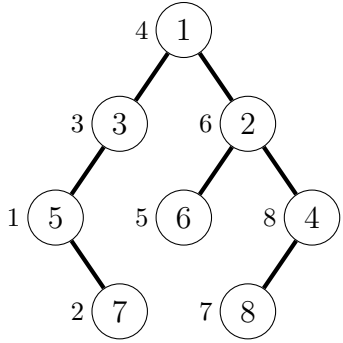

Figure~\ref{fig.bt} is an example of an increasing binary tree on the vertex set $[8]$.


\bigskip

{\bf Increasing restricted ternary trees.}
Our second family consists of increasing ternary trees 
on the vertex set $[n]$
such that middle children do not have any siblings.
We call these \textbfit{increasing restricted ternary trees} (RTs for short).
We will show (Corollary~\ref{cor.abt.numbers})
that the sequence of cardinalities of these trees
is generated by the T-fraction with
$\alpha$'s given by $1,1,2,2,3,3,4,4,5,5,\ldots$
and $\delta$'s given by $0,1,0,2,0,3,0,4,0,5,\ldots$:
that is, \reff{eq.alphadelta.quasiaffine}
with $x=y=u=v=b=d=1$ and $a=c=0$.
This sequence begins as
\be
   (a_n)_{n \ge 0}  \;=\;
   1, 1, 3, 11, 51, 295, 2055, 16715, 155355, 1624255, 18868575, \ldots
\label{eq.a=c=0}
\ee
and is \cite[A230008]{OEIS}.

The shifted sequence $(a_n)_{n\geq 1}$ was previously studied
by Kuba and Panholzer in \cite[Section~5.2, Example~5 and Remark~8]{Kuba_16}.
However, they considered a different combinatorial interpretation:
namely, {\em binary free multilabeled increasing trees}.
They showed that the exponential generating function
$F(t) = \sum_{n=1}^\infty a_n t^n/n!$ 
satisfies the ordinary differential equation
\be
F'(t) \;=\; F(t)^2 + 3 F(t) + 1\;.
\ee
It is not difficult to see,
using the general theory of increasing trees \cite{Bergeron_92},
that the exponential generating function
for increasing restricted ternary trees
satisfies this same differential equation.
Indeed, one can also show directly that, for $n \ge 1$,
binary free multilabeled increasing trees of size $n$
are in bijection with increasing restricted ternary trees of size $n$;
we do this by using the bijection in \cite[Theorem~10]{Kuba_16}
and then redrawing the black vertices in their bijection
to have middle children (see Section~\ref{subsubsec.free}).
We denote by $\RT_n$
the set of all increasing restricted ternary trees on the vertex set $[n]$.

\begin{figure}[t]
\centering
\newcommand{\nodea}{\node[draw,circle] (a) {$1$}
;}\newcommand{\nodeb}{\node[draw,circle] (b) {$2$}
;}\newcommand{\nodec}{\node[draw,circle] (c) {$3$}
;}\newcommand{\noded}{\node[draw,circle] (d) {$4$}
;}\newcommand{\nodee}{\node[draw,circle] (e) {$5$}
;}\newcommand{\nodef}{\node[draw,circle] (f) {$6$}
;}
\begin{tikzpicture}
\matrix[column sep=.1cm, row sep=0.5cm,ampersand replacement=\&]{
\& \& \nodea \& \\
\& \nodeb \& \& \nodec \\
\& \noded  \& \nodee \& \\
\& \& \& \nodef\\
};
        \draw[ultra thick] (a) -- (b) node {};
        \draw[ultra thick] (a) -- (c) node {};
        \draw[ultra thick] (b) -- (d) node {};
        \draw[ultra thick] (c) -- (e) node {};
        \draw[ultra thick] (e) -- (f) node {};
\end{tikzpicture}
   \vspace*{-2mm}
\caption{
   An example of an increasing restricted ternary tree on the vertex set $[6]$.
}
   \vspace*{1cm}
\label{fig.rt}
\end{figure}
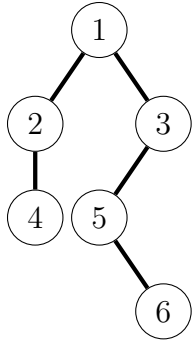

Figure~\ref{fig.rt} is an example of an increasing restricted ternary tree on the vertex set $[6]$.

\bigskip

{\bf Increasing interval-labeled restricted ternary trees.}
Our third family is defined as follows:
An \textbfit{increasing interval-labeled restricted ternary tree}
(IRT for short)
on the label set $[0,n]$
is a vertex-labeled tree satisfying the following rules:
\begin{itemize}
   \item The underlying tree is a ternary rooted tree
      in which middle children do not have any siblings.
   \item Each vertex $v$ in the tree is assigned an interval of integer labels
      $L_v = [i_1,i_2] \subseteq [0,n]$ such that 
      for every $i\in [0,n]$ there exists exactly one vertex $v$
      with $i\in L_v$.
      Thus, the collection of vertex labels forms an interval-partition
      of $[0,n]$.
   \item If a vertex $v$ has a middle child, then $|L_v| =1$,
      i.e., the vertex $v$ gets a single label.
   \item The labels are increasingly assigned, i.e., 
       for every pair of vertices $v,w$ such that $v$ is the parent of $w$,
       we impose that $\max L_v < \min L_w$. 
   \item The root (which is the vertex containing $0$ in its label set)
       can only have a left child, not a middle or right child.
\end{itemize}
We stress that $n+1$ is the total number of integer {\em labels}\/;
the number of {\em vertices}\/ can be anything from 1 to $n+1$.
We denote by $\IRT_n$ the set of all IRTs on the label set $[0,n]$.
In Figure~\ref{fig_IRT_n=0,1,2} we show all IRTs
on the label set $[0,n]$ with $n=0,1,2$,
and in Figure~\ref{fig.IRT.big} we show an example of an IRT
on the label set $[0,n]$ with $n=16$.

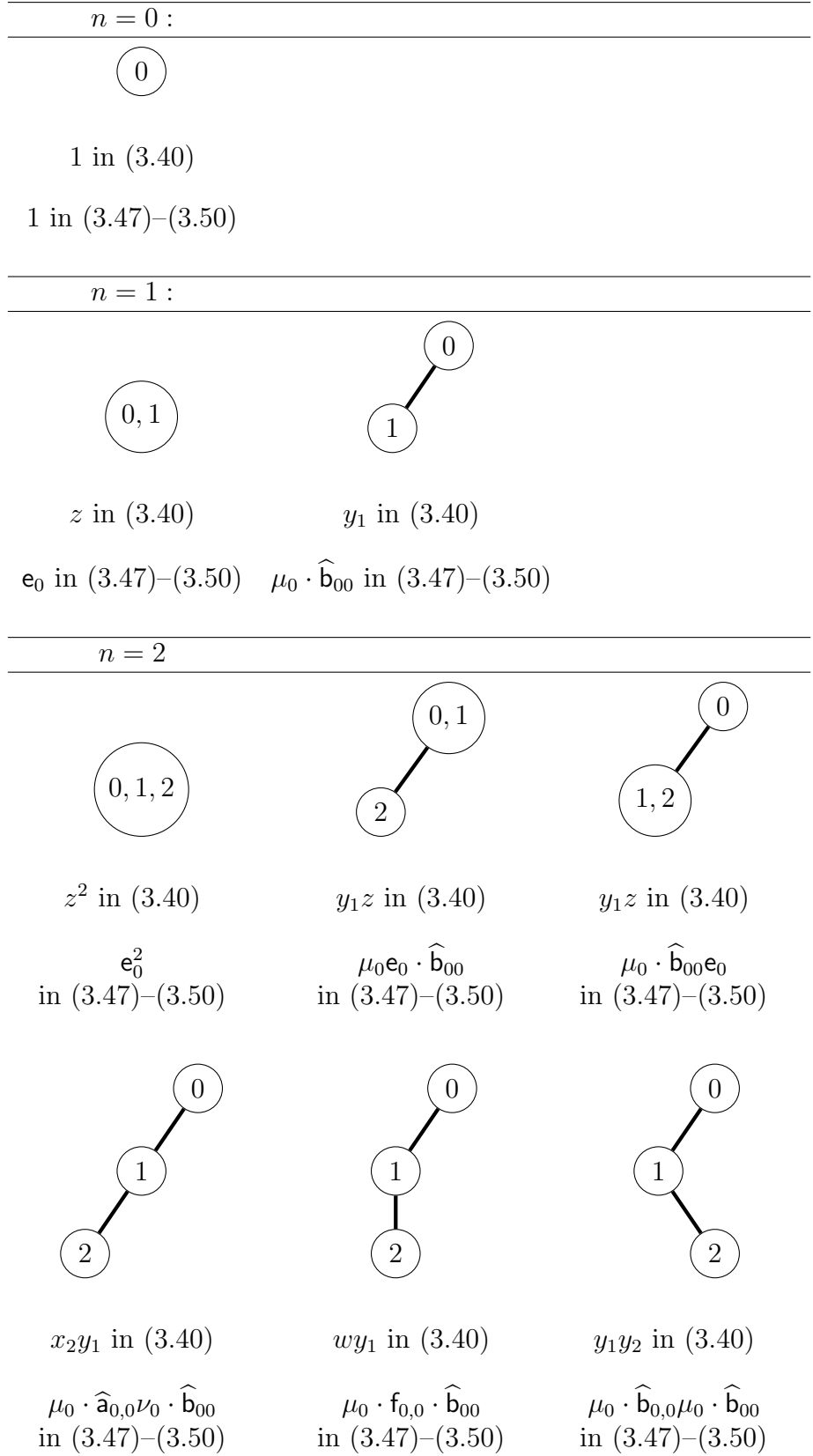
\begin{figure}[p]
\centering
\begin{tabular}{c c c c c}
\hline
	$n=0:$  & & &\\
\hline
{
\newcommand{\nodea}{\node[draw,circle] (a) {$0$};}
\begin{tikzpicture}
\matrix[column sep=.1cm, row sep=0.5cm,ampersand replacement=\&]{
 \nodea  \\
};
\end{tikzpicture}}  & & & \\[4mm]
 $1$ in  \eqref{eq.Pn.sbt} &&&\\[4mm]
 $1$ in \eqref{eq.weight.master.sbt.I}--\eqref{eq.weight.master.sbt} &&&\\[6mm]
\hline
$n=1:$  & & &\\
\hline
 {
\newcommand{\nodea}{\node[draw,circle] (a) {$0,1$};}
\begin{tikzpicture}
\matrix[column sep=.1cm, row sep=0.5cm,ampersand replacement=\&]{
 \nodea  \\
};
\end{tikzpicture}}
& 
{
\newcommand{\nodea}{\node[draw,circle] (a) {$0$}
;}\newcommand{\nodeb}{\node[draw,circle] (b) {$1$}
;}
\begin{tikzpicture}
\matrix[column sep=.1cm, row sep=0.5cm,ampersand replacement=\&]{
	\& \nodea \\
	\nodeb \& \\
};
\draw[ultra thick] (a) -- (b) node {};
\end{tikzpicture}}
& &\\[4mm]
$z$ in  \eqref{eq.Pn.sbt} 
& $y_1$ in  \eqref{eq.Pn.sbt} 
&&\\[4mm]
$\sfe_0$ in \eqref{eq.weight.master.sbt.I}--\eqref{eq.weight.master.sbt}
&
$\mu_0 \cdot  \sfbhat_{00}$  in \eqref{eq.weight.master.sbt.I}--\eqref{eq.weight.master.sbt}
&&\\[6mm]
\hline
$n=2$  & & &\\
\hline
{
\newcommand{\nodea}{\node[draw,circle] (a) {$0,1,2$};}
\begin{tikzpicture}
\matrix[column sep=.1cm, row sep=0.5cm,ampersand replacement=\&]{
 \nodea  \\
};
\end{tikzpicture}}
& 
{
\newcommand{\nodea}{\node[draw,circle] (a) {$0,1$}
;}\newcommand{\nodeb}{\node[draw,circle] (b) {$2$}
;}
\begin{tikzpicture}
\matrix[column sep=.1cm, row sep=0.5cm,ampersand replacement=\&]{
        \& \nodea \\
        \nodeb \& \\
};
\draw[ultra thick] (a) -- (b) node {};
\end{tikzpicture}}
& 
{
\newcommand{\nodea}{\node[draw,circle] (a) {$0$}
;}\newcommand{\nodeb}{\node[draw,circle] (b) {$1,2$}
;}
\begin{tikzpicture}
\matrix[column sep=.1cm, row sep=0.5cm,ampersand replacement=\&]{
        \& \nodea \\
        \nodeb \& \\
};
\draw[ultra thick] (a) -- (b) node {};
\end{tikzpicture}}
&\\[4mm]
$z^2$ in  \eqref{eq.Pn.sbt} 
& $y_1 z$ in  \eqref{eq.Pn.sbt} 
& $y_1 z$ in  \eqref{eq.Pn.sbt} 
\\[4mm]
$\sfe_0^2$ 
&
$\mu_0 \sfe_0 \cdot \sfbhat_{00}$ 
&
$\mu_0 \cdot \sfbhat_{00} \sfe_0$
&\\
 in \eqref{eq.weight.master.sbt.I}--\eqref{eq.weight.master.sbt}
&
 in \eqref{eq.weight.master.sbt.I}--\eqref{eq.weight.master.sbt}
&
 in \eqref{eq.weight.master.sbt.I}--\eqref{eq.weight.master.sbt}
\\[6mm]
{
\newcommand{\nodea}{\node[draw,circle] (a) {$0$}
;}\newcommand{\nodeb}{\node[draw,circle] (b) {$1$}
;}\newcommand{\nodec}{\node[draw,circle] (c) {$2$}
;}
\begin{tikzpicture}
\matrix[column sep=.1cm, row sep=0.5cm,ampersand replacement=\&]{
        \& \& \nodea \\
	\& \nodeb \& \\
	\nodec \&  \& \\
};
\draw[ultra thick] (a) -- (b) node {};
\draw[ultra thick] (b) -- (c) node {};
\end{tikzpicture}}
& 
{
\newcommand{\nodea}{\node[draw,circle] (a) {$0$}
;}\newcommand{\nodeb}{\node[draw,circle] (b) {$1$}
;}\newcommand{\nodec}{\node[draw,circle] (c) {$2$}
;}
\begin{tikzpicture}
\matrix[column sep=.1cm, row sep=0.5cm,ampersand replacement=\&]{
        \& \& \nodea \\
        \& \nodeb \& \\
	\& \nodec   \& \\
};
\draw[ultra thick] (a) -- (b) node {};
\draw[ultra thick] (b) -- (c) node {};
\end{tikzpicture}}
&
{
\newcommand{\nodea}{\node[draw,circle] (a) {$0$}
;}\newcommand{\nodeb}{\node[draw,circle] (b) {$1$}
;}\newcommand{\nodec}{\node[draw,circle] (c) {$2$}
;}
\begin{tikzpicture}
\matrix[column sep=.1cm, row sep=0.5cm,ampersand replacement=\&]{
        \& \& \nodea \\
        \& \nodeb \& \\
	\& \& \nodec  \\
};
\draw[ultra thick] (a) -- (b) node {};
\draw[ultra thick] (b) -- (c) node {};
\end{tikzpicture}}
\\[4mm]
 $x_2 y_1$ in  \eqref{eq.Pn.sbt} 
& $w y_1$ in  \eqref{eq.Pn.sbt} 
& $y_1 y_2$ in  \eqref{eq.Pn.sbt} 
\\[4mm]
$\mu_0\cdot \sfahat_{0,0} \nu_{0} \cdot \sfbhat_{00}$
&
$\mu_0 \cdot \sff_{0,0} \cdot \sfbhat_{00}$ 
&
$\mu_0\cdot \sfbhat_{0,0} \mu_{0} \cdot \sfbhat_{00}$ 
&\\
 in \eqref{eq.weight.master.sbt.I}--\eqref{eq.weight.master.sbt}
&
 in \eqref{eq.weight.master.sbt.I}--\eqref{eq.weight.master.sbt}
&
 in \eqref{eq.weight.master.sbt.I}--\eqref{eq.weight.master.sbt}
\\
\end{tabular}
   \vspace{3mm}
\caption{
   All IRTs on label set $[0,n]$ with $n=0,1,2$,
   along with their respective weights as per equations~\reff{eq.Pn.sbt}
   and~\eqref{eq.weight.master.sbt.I}--\eqref{eq.weight.master.sbt}.
}
\label{fig_IRT_n=0,1,2}
\end{figure}

\begin{figure}[t]
\centering
\newcommand{\nodea}{\node[draw,circle] (a) {$0,1$}
;}\newcommand{\nodeb}{\node[draw,circle] (b) {$2,3$}
;}\newcommand{\nodec}{\node[draw,circle] (c) {$4$}
;}\newcommand{\noded}{\node[draw,circle] (d) {$5$}
;}\newcommand{\nodee}{\node[draw,circle] (e) {$6$}
;}\newcommand{\nodef}{\node[draw,circle] (f) {$7$}
;}\newcommand{\nodeg}{\node[draw,circle] (g) {$8$}
;}\newcommand{\nodeh}{\node[draw,circle] (h) {$9$}
;}\newcommand{\nodei}{\node[draw,circle] (i) {$10$}
;}\newcommand{\nodej}{\node[draw,circle] (j) {$11$}
;}\newcommand{\nodek}{\node[draw,circle] (k) {$12,13$}
;}\newcommand{\nodel}{\node[draw,circle] (l) {$14,15$}
;}\newcommand{\nodem}{\node[draw,circle] (m) {$16$}
;}
\begin{tikzpicture}
\matrix[column sep=.1cm, row sep=0.5cm,ampersand replacement=\&]{
\&\& \& \& \nodea  \\
\&\&   \nodeb \& \&  \\
\& \nodec \& \& \noded \& \\
\nodeg \& \& \nodef \& \& \nodee\\
\& \nodeh  \& \nodei\& \& \nodej\\
\& \nodel \& \& \& \nodek\\
	\& \& \& \&  \& \nodem\\
};
        \draw[ultra thick] (a) -- (b) node {};
        \draw[ultra thick] (b) -- (c) node {};
        \draw[ultra thick] (b) -- (d) node {};
        \draw[ultra thick] (c) -- (g) node {};
        \draw[ultra thick] (d) -- (f) node {};
        \draw[ultra thick] (d) -- (e) node {};
        \draw[ultra thick] (e) -- (j) node {};
        \draw[ultra thick] (f) -- (i) node {};
        \draw[ultra thick] (g) -- (h) node {};
        \draw[ultra thick] (h) -- (l) node {};
        \draw[ultra thick] (j) -- (k) node {};
        \draw[ultra thick] (k) -- (m) node {};

	\node[left = -1.2mm of a] () {{\footnotesize 1:13}};
	\node[left = -1.2mm of b] () {{\footnotesize 2:5}};
	\node[left = -1.2mm of c] () {{\footnotesize 3:4}};
	\node[left = -1.2mm of d] () {{\footnotesize 7:8}};
	\node[left = -1.2mm of e] () {{\footnotesize 10:9}};
	\node[left = -1.2mm of f] () {{\footnotesize 8:6}};
	\node[left = -1.2mm of g] () {{\footnotesize 4:1}};
	\node[left = -1.2mm of h] () {{\footnotesize 5:2}};
	\node[left = -1.2mm of i] () {{\footnotesize 9:7}};
	\node[left = -1.2mm of j] () {{\footnotesize 11:10}};
	\node[left = -1.2mm of k] () {{\footnotesize 12:11}};
	\node[left = -1.2mm of l] () {{\footnotesize 6:3}};
	\node[left = -1.2mm of m] () {{\footnotesize 13:12}};
\end{tikzpicture}
   \vspace*{-2mm}
\caption{An example of an increasing interval-labeled restricted ternary tree on the label set $[0,16]$.
To the left of each vertex are two numbers $a$:$b$,
where $a$ is the order of the vertex as per the preorder traversal,
and $b$ is the order of the vertex as per traversal order $\bfA'$ 
defined as follows:
the left child if any, then the root, then the middle child if any,
and then the right child if any, all implemented recursively.
}
   \vspace*{1cm}
\label{fig.IRT.big}
\end{figure}
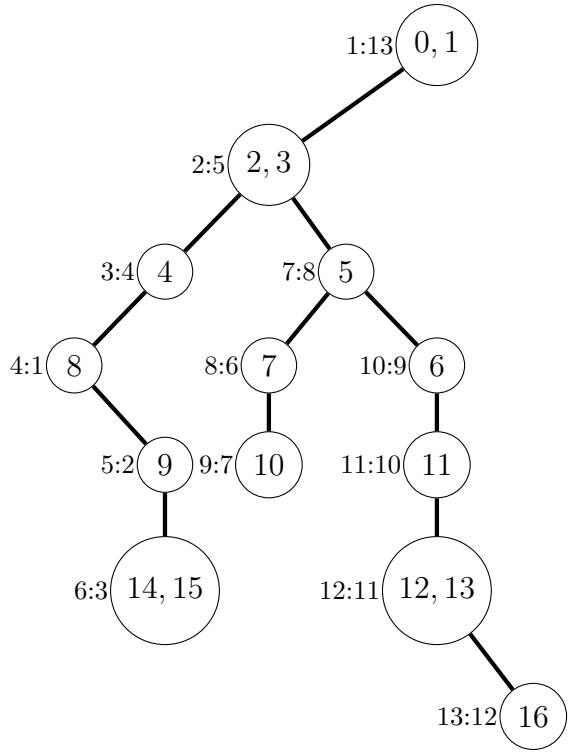

Clearly, if we impose the further condition
that all label sets have cardinality~1,
and we then remove the root
(which in this case necessarily has label set $\{0\}$),
then what remains is an increasing restricted ternary tree
on the vertex set $[n]$
(and conversely).
These trees therefore generalize RTs
by allowing the vertex labels (including that of the root)
to be intervals of cardinality $> 1$.

We will show (Corollary~\ref{cor.sbt.numbers})
that the sequence of cardinalities of these trees
is generated by the \hbox{T-fraction}
with $\alpha$'s given by $1,1,2,2,3,3,4,4,5,5,\ldots$
and $\delta$'s given by
$1,1,1,2,1,3,1,4,1,5,\ldots$:
that is, \reff{eq.alphadelta.quasiaffine}
with $x=y=u=v=a=b=d=1$ and $c=0$.
As already mentioned, this sequence is \reff{eq.c=0}
and is not at present in the OEIS.

Our motivation behind defining this class of trees
is a simple identity on T-fractions (Proposition~\ref{prop.transformation}).
With this identity in mind, one can view our definition
of interval-labeled restricted ternary trees
as a deformation of the definition of restricted ternary trees;
and we can in fact use Proposition~\ref{prop.transformation}
to \emph{prove} some of our results on IRTs as corollaries of those for RTs
(Section~\ref{subsec.algebraic.proofs.IRT}).

\bigskip

Finally, we pose an open problem:

\begin{openproblem}
\rm
Find a ``nice'' combinatorial interpretation for the T-fraction with $\alpha$'s
     given by 
     $1,1,2,2,3,3,4,4,5,5,\ldots$
     and $\delta$'s given by
     $1,1,2,2,3,3,4,4,5,5,\ldots$:
     that is,~\reff{eq.alphadelta.quasiaffine}
with $x=y=u=v=a=b=c=d=1$.
\end{openproblem}







\subsection{Outline of this paper}
   \label{subsec.intro.outline}

The remainder of this paper is structured as follows.
We begin, in Section~\ref{sec.preliminaries},
by introducing some necessary preliminaries.
Then, in Section~\ref{sec.results}, we state our main results.
After that we proceed to the proofs:
the preliminaries for the bijective proofs
are given in Section~\ref{sec.prelimproofs},
and the bijective proofs are presented in Section~\ref{sec.bijections.proofs}.
Also, for some of our less powerful results
we have very simple algebraic proofs:
these are presented in Section~\ref{sec.algebraic.proofs}.
Finally, in Section~\ref{sec.other.interpretations},
we will provide some other combinatorial interpretations of our T-fractions;
this section includes a conjecture
on the equidistribution of some vincular patterns on permutations
(Conjecture~\ref{conj.trivariate}).

This work would not have been possible without the existence
of the On-Line Encyclopedia of Integer Sequences~\cite{OEIS}.
Indeed, we began this project by doing an exhaustive search of the OEIS
for sequences with quasi-affine T-fraction coefficients of period~2,
as in \reff{eq.alphadelta.quasiaffine};
we will explain our procedure and results in Appendix~\ref{sec.OEIS}.
In Appendix~\ref{sec.cfg.operator}
we show how the sequences enumerating our families of trees
can be obtained in a very simple way using context-free (Chen) grammars.
We used these formulae to help guess our families of trees,
and also to help check that the sequences obtained from
the derivative operators and the T-fractions match.

\section{Preliminaries}
\label{sec.preliminaries}

\subsection{Contraction and transformation formulae}
\label{subsec.prelimtransform}

The formulae for even and odd contraction of an S-fraction
to an equivalent J-fraction are well known:
see e.g.~\cite[Lemmas~1 and 2]{Dumont_94b} \cite[Lemma~1]{Dumont_95}
for very simple algebraic proofs,
and see \cite[pp.~V-31--V-32]{Viennot_83}
for enlightening combinatorial proofs
based on grouping pairs of steps in a Dyck path.
These formulae have also been extended
to suitable subclasses of T-fractions
\cite{Sokal_totalpos} \cite[Propositions~2.1 and 2.2]{Deb-Sokal_Genocchi}.
In this paper we will need only
odd contraction \cite[Proposition~2.2]{Deb-Sokal_Genocchi}:

\begin{proposition}[Odd contraction for T-fractions with $\delta_1=\delta_3=\delta_5 = \ldots = 0$]
   \label{prop.contraction_odd.Ttype}
\hfill\break
\vspace*{4mm}
\hspace*{-2mm}We have
\begin{eqnarray}
   & &  \hspace*{-15mm}
   \cfrac{1}{1 - \cfrac{\alpha_1 t}{1 - \delta_2 t - \cfrac{\alpha_2 t}{1 -  \cfrac{\alpha_3 t}{1- \delta_4 t - \cfrac{\alpha_4 t}{1 - \cdots}}}}}
   \;\;=\;  \nonumber \\[3mm]
   & &
   1 \:+\:
   \cfrac{\alpha_1 t}{1 - (\alpha_1 + \alpha_2 + \delta_2) t - \cfrac{\alpha_2 \alpha_3 t^2}{1 - (\alpha_3 + \alpha_4 + \delta_4) t - \cfrac{\alpha_4 \alpha_5 t^2}{1 - \cdots}}}
 \;\;.
   \qquad
 \label{eq.contraction_odd.Ttype}
\end{eqnarray}
That is, the T-fraction on the left-hand side of
\reff{eq.contraction_odd.Ttype}
equals $1$ plus $\alpha_1 t$ times the J-fraction with coefficients
\begin{subeqnarray}
   \gamma_n  & = &  \alpha_{2n+1} + \alpha_{2n+2} + \delta_{2n+2} \\
   \beta_n  & = &  \alpha_{2n} \alpha_{2n+1}
 \label{eq.contraction_odd.coeffs.Ttype}
\end{subeqnarray}
Here \reff{eq.contraction_odd.Ttype}/\reff{eq.contraction_odd.coeffs.Ttype}
holds as an identity in $\Z[\balpha,\bdelta_{\rm even}][[t]]$,
where $\bdelta_{\rm even} = (0,\delta_2,0,\delta_4,\ldots)$.
\end{proposition}

Both the algebraic and the combinatorial proofs of the contraction formulae
for S-fractions can be easily generalized \cite{Sokal_totalpos}
to prove the contraction formulae for T-fractions.

Please note that Proposition~\ref{prop.contraction_odd.Ttype}
allows the J-fraction on the right-hand side of \reff{eq.contraction_odd.Ttype}
to be converted to a T-fraction on the left-hand side in numerous ways:
$\alpha_1$ can be chosen arbitrarily;
then each $\beta_n$ can be factored arbitrarily into
$\alpha_{2n}$ and $\alpha_{2n+1}$;
and then $\delta_{2n+2}$ is {\em defined}\/ by
(\ref{eq.contraction_odd.coeffs.Ttype}a).
Of course, some of the resulting T-fractions may have
coefficients $\bdelta$ that are ``nice'' in one sense or another,
while others will not.

\bigskip

We now prove a useful transformation formula for T-fractions:

\begin{proposition}[T-fraction with no odd delta to generic T-fraction]
\label{prop.transformation}
\be
\cfrac{1/(1-\delta_1 t)}{1 -  \cfrac{\alpha_1 t/(1-\delta_1 t)}{1 - \delta_2 t - \cfrac{\alpha_2 t/(1- \delta_3 t)}{1 - \cfrac{\alpha_3 t/(1-\delta_3 t)}{1-\delta_4 t-\ldots}}}}
   \;\;=\;\;
   \cfrac{1}{1 - \delta_1 t - \cfrac{\alpha_1 t}{1 - \delta_2 t - \cfrac{\alpha_2 t}{1 - \delta_3 t - \cfrac{\alpha_3 t}{1 - \ldots}}}}
   \;\:.
   \quad
\label{eq.transformation}
\ee
\end{proposition}

\proof
This is obtained by repeated use of the identity
\be
\dfrac{1/f(t)}{1 - g(t)/f(t)}
\;=\;
\dfrac{1}{f(t) - g(t)}
\label{eq.comm.ring.identity}
\ee
at alternate levels
of the left-hand side of~\eqref{eq.transformation},
with $f(t) = 1-\delta_1 t, 1-\delta_3 t, \ldots\;$.
\qed

Let $T(t;\balpha,\bdelta)$ denote the T-fraction~\eqref{eq.Tfrac.def}.
Also, let $\bdelta_{\rm even} = \left.\bdelta\right|_{\delta_{2k-1}=0}$,
i.e., $\bdelta_{\rm even}$ denotes the sequence $\bdelta$
with all odd-order $\delta_i$ set to $0$.
Then, Proposition~\ref{prop.transformation}
is equivalent to the following:
\be
\dfrac{1}{1-\delta_1 t}
\;
T\left(t;
\left.\balpha\right|_{
\alpha_{2k-1} \mapsto \frac{\alpha_{2k-1}}{1 - \delta_{2k-1} t}
\, ,\,
\alpha_{2k} \mapsto \frac{\alpha_{2k}}{1 - \delta_{2k+1} t}
},
\bdelta_{\rm even}\right)
\;=\;
T(t;\balpha,\bdelta)\;.
\label{eq.prop.transformation.equiv}
\ee
This identity gives us a recipe to start
from a T-fraction with $\delta_i$ only at even orders
and then insert $\delta_i$ at odd orders as well.
This is what motivated our construction of
increasing interval-labeled restricted ternary trees.
At the end of Section~\ref{subsec.prelimproofs.1}
we will also give a combinatorial interpretation/proof
of Proposition~\ref{prop.transformation}.

%
%

\subsection{Types of labeled trees}

In this paper, all trees are rooted and finite;
we henceforth omit these two adjectives.
A \textbfit{labeled tree} with vertex set $V$
is a tree in which each vertex is assigned a {\em distinct}\/ label from~$V$,
such that each element of $V$ is the label for some vertex of the tree
(hence for exactly one vertex).
A \textbfit{multilabeled tree} with label set $\scrl$
is a tree in which each vertex $v$ is assigned a nonempty label set
$L_v \subseteq \scrl$,
such that each label $i \in \scrl$ belongs to exactly one of the sets $L_v$.
Specific types of multilabeled trees arise by restricting the subsets of
$\scrl$ that are allowed to be vertex label sets  $L_v$.
For instance, in~\textbfit{$\bm{k}$-labeled trees}
\cite[Section~2.6]{Kuba_21},
all vertex label sets $L_v$ must have cardinality~$k$;
in~\textbfit{interval-labeled trees},
the label set $\scrl$ is some set of integers (usually an interval),
and all vertex label sets $L_v$ must be intervals.

A labeled tree on a partially ordered vertex set $V$
is called \textbfit{increasing}
if the label of a child is always greater than the label of its parent.
In most applications the vertex set is totally ordered:
this is the case for ordinary (unilabeled) increasing trees,
in which the vertex set is some set of integers
(ordinarily either $[n]$ or $[0,n]$)
and each vertex has a single integer label.
However, for multilabeled trees the vertex set $V$
is some collection of nonempty subsets of the label set $\scrl$,
equipped with the partial order $A < B$
$\Longleftrightarrow$ $\max A < \min B$.
We remark that although the set of {\em all}\/ intervals in $[0,n]$
fails to be totally ordered when $n \ge 2$, the set of intervals occurring
as vertex labels for a {\em given}\/ interval-labeled restricted ternary tree
(which is, by definition, the vertex set of that tree)
forms an interval-partition of $[0,n]$ and hence is totally ordered.

\subsection{Tree statistics: Node types and label surplus}

We will use the notion of node types for $k$-ary trees,
as introduced by Kuba and Varvak \cite{Kuba_21}. 
%
Let $T$ be a $k$-ary rooted tree, and let $V(T)$ be its vertex set. 
For $v\in V(T)$,
we define the \textbfit{node type} of $v$ to be the string
$N(v,T)\eqdef n_{v,1}\cdots n_{v,k}\in \{0,1\}^{k}$
where $n_{v,j}=1$ if the $j$-th child of the vertex $v$ exists,
and 0 if it does not.
%
%
%
Let us mention the possible node types in some families of trees:
\begin{itemize}
    \item The possible node types in binary trees are $00, 01, 10,11$.
    \item The possible node types in ternary trees are $000, 001, 010,100,$ $011,101,110,111$.
    \item The possible node types in restricted ternary trees are $000,100,010,001,101$.
\end{itemize}

For a string $s\in \{0,1\}^{k}$ with $k \ge 1$,
we define $I_T(s)$ to be
\be
I_T(s)\;\eqdef\;
\#\{v\colon\: v\in V(T) \text{ and } N(v,T) = s \}\;.
\label{eq.def.node.count}
\ee
Thus, $I_T(s)$ counts the number of vertices in the tree $T$
with node type $s$.


In some cases we will want to give weights only to {\em non-root}\/ vertices.
We therefore define also
\be
I'_T(s)\;\eqdef\;
\#\{v\colon\: v\in V(T) \setminus \text{root} \text{ and } N(v,T) = s \} \;,
\label{eq.def.node.count.prime}
\ee
so that $I'_T(s)$ counts the number of {\em non-root}\/ vertices in the tree $T$
with node type $s$.

In this paper we also study some $k$-ary trees
in which the vertices are allowed to have multiple labels.
For this reason, given a multilabeled $k$-ary tree, we also introduce
$I_T(\varepsilon)$ (here $\varepsilon$ is the empty string) to denote
\begin{equation}
I_T(\varepsilon) \;\eqdef\; \#\{i\colon\: i \text{ is a label in }T\}
   \:-\: \#\{\text{$v\colon\: v$ is a vertex in $T$}\}\;.
\label{eq.def.node.count.varep}
\end{equation}
We call the quantity $I_T(\varepsilon)$
the \textbfit{label surplus of the tree $\bm{T}$}.
Using $L_v$ to denote the label set of the vertex $v$,
we have equivalently
\be
I_\varepsilon(T) 
\;=\; 
\sum_{v\in V(T)} (|L_v| - 1)
\label{eq.label.surplus.dist}
\ee
(since the label sets $L_v$ are disjoint).
We therefore call the quantity $|L_v| - 1$
the \textbfit{label surplus of the vertex $\bm{v}$}.


\subsection{Tree statistics: Crossings and nestings}
  \label{subsec.crossnest}

We now introduce some new statistics on trees.
Similar statistics were used implicitly in
\cite{Kuba_21,latpath_lah,latpath_SRTR}
and will be used explicitly in \cite{masterpoly_trees}.

Let $T$ be an increasing tree on a totally ordered vertex set $V$.
(This includes the case of increasing interval-labeled trees,
for which the vertex set is a totally ordered set of intervals in $[0,n]$.)
Let $v_0, v_1, \ldots, v_m$ be the vertices of $T$ in this total order
(here $v_0$ must be the root, $v_1$ must be a child of the root,
 and $v_m$ must be a leaf).
For any non-root vertex $v$, we denote by $p(v)$ the parent of $v$ in $T$;
if $v$ is the root, then $p(v)$ is undefined.
We then define the \textbfit{level} of a vertex $v$ as
\begin{equation}
\lev(v,T)
 \;\eqdef\;
\#\{w\colon\: p(w)< v < w\} \;.
\label{eq.def.level.vertex.sbt}
\end{equation}
It is clear that $\lev(v_0,T) = \lev(v_m,T) = 0$.
Note also that when $v = v_1$, which is a child of the root,
any $w$ contributing in \reff{eq.def.level.vertex.sbt}
must be another child of the root;
so if the root has only one child
--- as is the case for our
increasing interval-labeled restricted ternary trees ---
then $\lev(v_1,T) = 0$.

Now suppose further that the children of each vertex are linearly ordered.
(This is certainly the case for $k$-ary trees, which include binary,
 ternary and restricted ternary trees as special cases.)
We use this linear order to introduce two new statistics:
{\em croix}\/ and {\em nid}\/
(the French words for crossings and nestings, respectively).
However, to introduce these, we need to choose first
a \textbfit{tree-traversal algorithm} $\bfA$:
by this we mean a mapping from triplets $(T,V,<)$ of
increasing ordered trees $T$ on a totally ordered vertex set $(V,<)$
to new total orders $<_\bfA$ on $V$,
that satisfies the following \textbfit{consistency property}:
\begin{quote}
   For any $j \in V$, define $V|_j = \{k \in V \colon\: k \le j\}$.
   Then for every tree $T$ on $V$ and every $j \in V$,
   the total order $<_\bfA$ on $T$ restricted to $V|_j$
   is the same as the total order $<_\bfA$ on the tree $T \restrict V|_j$.
\end{quote}
In other words, the vertices of $T \restrict V|_j$
are traversed in the same relative order as they are traversed in $T$.
Some examples of consistent tree-traversal algorithms are:
\begin{itemize}
   \item {\bf Preorder traversal.}
       First the root, then the children in order, all implemented recursively.
   \item {\bf Postorder traversal.}
       The children in order, then the root, all implemented recursively.
   \item {\bf Inorder (= symmetric) traversal for binary trees}
       \cite[pp.~44--45]{Stanley_12}.
       The left child if any, then the root, then the right child if any,
       all implemented recursively.\footnote{
   Inorder traversal can be generalized in a variety of ways to
   non-binary ordered trees.
   For instance, fix your favorite integer $c \ge 0$;
   traverse the first $c$ children
   (or all of the children if there are fewer than $c$ of them),
   then the root, then the remaining children, implemented recursively.
   Or alternatively:
   traverse all but the last $c$ children (if there are more than $c$ of them),
   then the root,
   then the last $c$ children
   (or all of the children if there are fewer than $c$ of them),
   implemented recursively.
}
\end{itemize}

Now choose a consistent tree-traversal algorithm $\bfA$.
We refine the definition~\reff{eq.def.level.vertex.sbt} of level as follows.
Consider a pair of vertices $v,w$ in the tree $T$ that satisfy $p(w)<v<w$.
We say that they form a
\begin{itemize}
    \item croix if $v<_\bfA w$,
    \item nid if $w<_\bfA v$.
\end{itemize}
We then define the
{\bf index-refined crossing and nestings statistics}
\begin{subeqnarray}
    \croix_\bfA(v,T) \;=\; \{w\colon\: p(w) < v < w \text{ and } v<_\bfA w \}\\
    \nid_\bfA(v,T) \;=\; \{w\colon\: p(w) < v < w \text{ and } w<_\bfA v \}
    \label{eq.def.croix.nid.vertex}
\end{subeqnarray}
It is clear from
definitions~\eqref{eq.def.level.vertex.sbt}/\eqref{eq.def.croix.nid.vertex}
that
\begin{equation}
    \lev(v,T) \;=\; \croix_\bfA(v,T) + \nid_\bfA(v,T) \;.
\label{eq.croix.plus.nid}
\end{equation}

In the remainder of this paper we will lighten the notation
by deleting the subscript $\bfA$ from $\croix$ and $\nid$,
but we stress that meaning of $\croix$ and $\nid$
{\em depends on the choice of tree-traversal algorithm}\/.
Nevertheless, the resulting ``master'' continued fractions
will be identical for all consistent tree-traversal algorithms.
Therefore, we will have proven that the families
$(\croix_\bfA,\nid_\bfA)$ and $(\croix_{\bfA'},\nid_{\bfA'})$
are {\em equidistributed}\/,
for any pair of consistent tree-traversal algorithms $\bfA$ and $\bfA'$.

\medskip

\begin{example}
\rm
Let $T$ be the increasing binary tree shown in Figure~\ref{fig.bt}.
It has three vertices of node type 00, two vertices of node type 11, 
two vertices of node type 10, and one vertex of node type 01.

Let $\bfA$ be the inorder traversal.
The vertices are traversed as follows:
$5 <_{\bfA} 7 <_{\bfA} 3 <_{\bfA} 1 <_{\bfA} 6 <_{\bfA} 2 <_{\bfA} 8 <_{\bfA} 4$.

We now note down
the node type $N(v,T)$ and 
the statistics $\lev(v,T)$, $\nid_{\bfA}(v,T)$ and $\croix_{\bfA}(v,T)$:


\begin{equation}
\begin{tabular}{c|c|c|c|c}
        $v$ &
	$N(v,T)$
	&
	$\lev(v,T)$
	&
        $\nid(v,T)$
        &
        $\croix(v,T)$\\
        \hline
	1 & 11 & 0 & 0 & 0\\
	2 & 11 & 1 & 1 & 0\\
	3 & 10 & 2 & 0 & 2\\
	4 & 10 & 2 & 2 & 0\\
	5 & 01 & 2 & 0 & 2\\
	6 & 00 & 2 & 1 & 1\\
	7 & 00 & 1 & 0 & 1\\
	8 & 00 & 0 & 0 & 0
\end{tabular}
  \label{eq.tab.bt.croix.nid}
\end{equation}
\myendremark
\end{example}

\smallskip

\begin{example}
\rm
Let $T$ be the IRT shown in Figure~\ref{fig.IRT.big}.
The total order on the vertices is 
      $\{0,1\} 
      < 
      \{2,3\}
      <
      \{4\}
      <
      \{5\}
      <
      \{6\}
      <
      \{7\}
      <
      \{8\}
      <
      \{9\}
      <
      \{10\}
      <
      \{11\}
      <
      \{12,13\}
      <
      \{14,15\}
      <\{16\}$.
There are three vertices of node type $000$,
two vertices each of node types $101,100,001$,
and four vertices of node type $010$.

Let $\bfA$ be the preorder traversal.
The vertices as traversed as follows:
	$\{0,1\}
	<_{\bfA}
	\{2,3\}
	<_{\bfA}
	\{4\}
	<_{\bfA}
	\{8\}
	<_{\bfA}
	\{9\}
	<_{\bfA}
	\{14,15\}
	<_{\bfA}
	\{5\}
	<_{\bfA}
	\{7\}
	<_{\bfA}
	\{10\}
	<_{\bfA}
	\{6\}
	<_{\bfA}
	\{11\}
	<_{\bfA}
	\{12,13\}
	<_{\bfA}
	\{16\}.
	$

Next let $\bfA'$ be the following traversal order: 
the left child if any, then the root, then the middle child if any,
and then the right child if any, all implemented recursively.
As per this order, the vertices are traversed as follows:
      $\{8\} 
      <_{\bfA'} 
      \{9\}
      <_{\bfA'}
      \{14,15\}
      <_{\bfA'}
      \{4\}
      <_{\bfA'}
      \{2,3\}
      <_{\bfA'}
      \{7\}
      <_{\bfA'}
      \{10\}
      <_{\bfA'}
      \{5\}
      <_{\bfA'}
      \{6\}
      <_{\bfA'}
      \{11\}
      <_{\bfA'}
      \{12,13\}
      <_{\bfA'}
      \{16\}
      <_{\bfA'}
      \{0,1\}$.

We now note down the 
node type $N(v,T)$
and the statistics $\lev(v,T)$, $\nid_{\bfA}(v,T)$ and $\croix_{\bfA}(v,T)$,
$\nid_{\bfA'}(v,T)$ and $\croix_{\bfA'}(v,T)$:

\begin{equation}
\begin{tabular}{c|c|c|c|c|c|c}
        $v$ &
	$N(v,T)$
	&
        $\lev(v,T)$
	&
	$\nid_{\bfA}(v,T)$
	&
	$\croix_{\bfA}(v,T)$	
	&
	$\nid_{\bfA'}(v,T)$
        &
	$\croix_{\bfA'}(v,T)$\\
        \hline
	\{0,1\} & 100 & 0 & 0 & 0 & 0 & 0  \\
	\{2,3\} & 101 & 0 & 0 & 0 & 0 & 0 \\
	\{4\} & 100 & 1 & 0 & 1 & 0 & 1 \\
	\{5\} & 101 & 1 & 1 & 0 & 1 & 0\\
	\{6\} & 010 & 2 & 2 & 0 & 2 & 0\\
	\{7\} & 010 & 2 & 1 & 1 & 1 & 1\\
	\{8\} & 001 & 2 & 0 & 2 & 0 & 2\\
	\{9\} & 010 & 2 & 0 & 2& 0 & 2\\
	\{10\} & 000 & 2 & 1 & 1 & 1 & 1\\
	\{11\} & 010 & 1 & 1 & 0 & 1 & 0\\
	\{12,13\} & 001 & 1 & 1 & 0 & 1 & 0\\
	\{14,15\} & 000 & 1 & 0 & 1 & 0 & 1\\
	\{16\} & 000 & 0 & 0 & 0 & 0 & 0\\
\end{tabular}
  \label{eq.tab.IRT.croix.nid}
\end{equation}
%
%
\myendremark
\end{example}

\section{Statement of results}
\label{sec.results}

All our continued fractions will come in two levels of generality:
``simple'' continued fractions, with finitely many indeterminates
that count node types and label surplus
[cf.\ \reff{eq.def.node.count}/\reff{eq.def.node.count.varep}];
and then ``master'' continued fractions,
with infinitely many indeterminates that count the pair $(\croix,\nid)$
at each vertex
[cf.\ \reff{eq.def.croix.nid.vertex}].

For pedagogical clarity, we begin by presenting these continued fractions
for increasing binary trees;
some (but not all) of these formulae are well known
and go back to Flajolet \cite{Flajolet_80}.
After this, we present our new results in increasing order of generality:
first increasing restricted ternary trees,
and then increasing interval-labeled restricted ternary trees.

\subsection{Increasing binary trees}
\label{subsec.bt}

Let $\mathcal{B}_n$ denote the set of increasing binary trees
on the vertex set $[n]$.
It is well known that this set has cardinality $|\mathcal{B}_n| = n!$:
see e.g.~\cite[pp.~44--45]{Stanley_12} for a bijective proof.
Here we define some polynomials that refine this counting,
and provide continued fractions for their ordinary generating functions.

\subsubsection{Simple J-fraction and T-fraction}

Consider first the polynomial in four variables
that enumerates increasing binary trees according to the node types:
for $n\geq 0$, we define
\begin{equation}
    P_n(x_1, x_2, y_1, y_2) 
    \;\eqdef\;
    \sum_{T\in \mathcal{B}_n} 
         x_1^{I_T(11)} y_1^{I_T(00)} x_2^{I_T(10)} y_2^{I_T(01)}
    \;,
\label{eq.Pn.bt}
\end{equation}
where we recall that $I_T(s)$ for a string $s$
was defined in~\eqref{eq.def.node.count}
as the number of vertices with node type $s$.
Thus, the variables $x_1, y_1, x_2, y_2$
are associated to the node types $11$, $00$, $10$, $01$, respectively.
In particular, we have $x_1$ or $x_2$ when a vertex has a left child,
and $y_1$ or $y_2$ when it does not;
we have a subscript 1 (resp.~2)
when the vertex has even (resp.~odd) out-degree.
By convention $P_0 = 1$ (corresponding to the empty tree);
for $n\geq 1$, the polynomial $P_n$ always has a factor $y_1$
since the vertex $n$ is always a leaf.

The polynomials~\reff{eq.Pn.bt} have the following beautiful J-fraction,
which is essentially~\cite[Theorem~3A]{Flajolet_80}
restated in terms of increasing binary trees
using the correspondence in~\cite[pp.~44--45]{Stanley_12}:

\begin{theorem}[{\cite[Theorem~3A]{Flajolet_80}+\cite[pp.~44--45]{Stanley_12}}]
\label{thm.Simple.J.bt}
The ordinary generating function of the polynomials
$P_{n+1}(x_1, x_2, y_1, y_2)$
has the J-fraction
\be
\sum_{n=0}^\infty P_{n+1}(x_1, x_2, y_1, y_2)\, t^n
\;=\;
\Scale[0.9]{
\cfrac{y_1}{1-(x_2 +y_2) t - \cfrac{2 x_1 y_1 \, t^2}{1-2(x_2+y_2) t - \cfrac{6 x_1 y_1 t^2}{1- \ldots}}}
}
\ee
with coefficients
\begin{eqnarray}
\gamma_n \;=\; (n+1)\,(x_2+y_2) \,,
\qquad 
\beta_n \;=\; n(n+1)\, x_1 y_1\;.
\label{eq.def.simple.bt.J}
\end{eqnarray}
\end{theorem}

When $x_1 = x_2 = y_1 = y_2 = 1$,
this is a J-fraction for the sequence $((n+1)!)_{n \ge 0}$.
It arises by even contraction of the S-fraction
with $\alpha_{2k-1} = k+1$, $\alpha_{2k} = k$.

\bigskip

We now obtain a T-fraction from Theorem~\ref{thm.Simple.J.bt}
by using odd contraction (Proposition~\ref{prop.contraction_odd.Ttype}):

\begin{theorem}
\label{thm.simple.bt.T}
The ordinary generating function of the polynomials
$P_{n}(x_1, x_2, y_1, y_2)$
has the T-fraction
\begin{eqnarray}
&&
\hspace*{-15mm}
\sum_{n=0}^\infty P_{n}(x_1, x_2, y_1, y_2)\, t^n\;\nonumber\\[2mm]
&&
\hspace*{-15mm}=\;
\Scale[0.9]{
\cfrac{1}{1- \cfrac{y_1 t}{1-(x_2+y_2 - x_1-y_1 ) t - \cfrac{ x_1 t}{1- \cfrac{2 y_1 t}{1  - 2(x_2+y_2 - x_1-y_1 ) t - \cfrac{2 x_1 t}{1-\ldots}}}}}
}
   \label{eq.thm.simple.bt.T}
\end{eqnarray}
with coefficients
\begin{subeqnarray}
    \alpha_{2k-1} & = & k y_1\\[2mm]
    \alpha_{2k} & = & k x_1\\[2mm]
    \delta_{2k-1} & = & 0\\[2mm]
    \delta_{2k} & = & k\, (x_2 + y_2  - x_1-y_1)
\label{eq.def.weights.simple.bt.T}
\end{subeqnarray}
\end{theorem}

If we now choose some specialization of $x_1,x_2,y_1,y_2$
that satisfies $x_1 + y_1 = x_2 + y_2$,
then from (\ref{eq.def.weights.simple.bt.T}d) we get $\delta_{2k} = 0$,
so that the T-fraction \reff{eq.thm.simple.bt.T} becomes an S-fraction.
In particular this happens if we take $x_1 = x_2$ and $y_1 = y_2$,
or alternatively if we take $x_1 = y_2$ and $y_1 = x_2$
(among many other possibilities).
Let us show the first of these:

\begin{corollary}
\label{cor.simple.bt.S}
The ordinary generating function of the polynomials~\reff{eq.Pn.abt} 
under the substitutions $x_1 = x_2 = x$ and $y_1 = y_2 = y$
has the S-fraction
\begin{eqnarray}
\sum_{n=0}^\infty P_{n}(x, x, y, y)\, t^n
\;=\;
\Scale[0.9]{
\cfrac{1}{1- \cfrac{y t}{1- \cfrac{ x t}{1- \cfrac{2 y t}{1  - \cfrac{2 x t}{1-  - \cfrac{3 y t}{1-\ldots}}}}}}
}
\end{eqnarray}
with coefficients
\begin{eqnarray}
    \alpha_{2k-1} \;=\; k y, \qquad    
    \alpha_{2k} \;=\; k x\;.
\label{eq.def.weights.simple.bt.S}
\end{eqnarray}
\end{corollary}

\noindent
This is the S-fraction for the homogenized Eulerian polynomials
\cite[section~79]{Stieltjes_1894} \cite[Section~2.2]{Sokal-Zeng_masterpoly}.
When $x=y=1$, it gives Euler's \cite[section~21]{Euler_1760}
S-fraction for the sequence $(n!)_{n \ge 0}$.

\subsubsection{Master J-fraction and T-fraction}

Fix a consistent tree-traversal algorithm $\bfA$.
Now let
$\bsfa = (\sfa_{\ell,\ell'})_{\ell,\ell'\geq 0}$,
$\bsfb = (\sfb_{\ell,\ell'})_{\ell,\ell'\geq 0}$,
$\bsfc = (\sfc_{\ell,\ell'})_{\ell,\ell'\geq 0}$,
$\bsfd = (\sfd_{\ell,\ell'})_{\ell,\ell'\geq 0}$
be infinite sets of indeterminates,
and define polynomials $Q_n(\bsfa, \bsfb,\bsfc,\bsfd)$ by
\begin{eqnarray}
Q_n(\bsfa, \bsfb,\bsfc,\bsfd)
&\eqdef&
\sum_{T\in \mathcal{B}_n}
\prod_{N(v,T) =11 }
\sfa_{\croix(v,T),\nid(v,T)}
\prod_{N(v,T) =00 }
\sfb_{\croix(v,T),\nid(v,T)} \times
\nonumber\\
&& \qquad\, \prod_{N(v,T) =10 }
\sfc_{\croix(v,T),\nid(v,T)}
\prod_{ N(v,T) =01 }
\sfd_{\croix(v,T),\nid(v,T)}
\;.
\label{eq.def.Qn.bt}
\end{eqnarray}
Since the vertex $n$ is always a leaf,
the polynomials~\reff{eq.def.Qn.bt} have a common factor $\sfb_{00}$
when $n \ge 1$.

We have the following master J-fraction,
which is implicit in \cite{Flajolet_80}
when the tree-traversal algorithm is taken to be
inorder (= symmetric) traversal;
but we write it out explicitly and allow {\em any}\/
consistent tree-traversal algorithm.

\begin{theorem}
\label{thm.bt.master}
The ordinary generating function of the polynomials
$Q_{n+1}(\bsfa, \bsfb,\bsfc,\bsfd)$
has the J-fraction
\begin{eqnarray}
&& 
\hspace*{-15mm}
\sum_{n=0}^\infty Q_{n+1}(\bsfa, \bsfb,\bsfc,\bsfd) \, t^n
\;=\; 
\nonumber \\[2mm]
&&
\hspace*{-15mm}
\cfrac{\sfb_{00}}{1 - (\sfc_{00}+\sfd_{00})\, t - \cfrac{\sfa_{00}(\sfb_{01} + \sfb_{10}) \, t^2}{1 - (\sfc_{01}+\sfc_{10}+\sfd_{01}+\sfd_{10})\, t - \cfrac{(\sfa_{01}+\sfa_{10})(\sfb_{02} +\sfb_{11} + \sfb_{20})\, t^2}{ 1 - \ldots}}}
\label{eq.bt.master.cf}
\end{eqnarray}
with coefficients
\begin{subeqnarray}
    \beta_{n} & = & 
     \biggl(
     \sum_{\xi=0}^{n-1} \sfa_{\xi, n-1-\xi}
     \biggr)
     \biggl(
     \sum_{\xi=0}^{n} \sfb_{\xi, n-\xi}
     \biggr)
    \\[2mm] 
    \gamma_{n} & = & \sum_{\xi=0}^n \sfc_{\xi,n-\xi} \, + \, 
    \sum_{\xi=0}^n \sfd_{\xi,n-\xi}
\label{eq.bt.master.weights}
\end{subeqnarray}

\end{theorem}
Here we will derive Theorem~\ref{thm.bt.master}
as a special case of Theorem~\ref{thm.abt.master.J}.

\bigskip

We can now obtain a T-fraction from Theorem~\ref{thm.bt.master}
by using odd contraction (Proposition~\ref{prop.contraction_odd.Ttype}):

\begin{theorem}
\label{thm.bt.master.T}
The ordinary generating function of the polynomials
$Q_n(\bsfa, \bsfb,\bsfc,\bsfd)$
has the T-fraction
\begin{eqnarray}
&& 
\hspace*{-10mm}
\sum_{n=0}^\infty Q_n(\bsfa, \bsfb,\bsfc,\bsfd) \, t^n
\;=\; 
\nonumber \\[2mm]
&&
\hspace*{-10mm}
\Scale[0.8]{\cfrac{1}{1 -  \cfrac{\sfb_{00} \, t}{1 - (\sfc_{00}+\sfd_{00}-\sfa_{00}-\sfb_{00})\, t - \cfrac{\sfa_{00}\, t}{ 1 - 
\cfrac{(\sfb_{01}+\sfb_{10})\,t }{1-
(\sfc_{01}+\sfc_{10}+\sfd_{01}+\sfd_{10}-\sfa_{01}-\sfa_{10}-\sfb_{01}-\sfb_{10})
 \,t -
\cfrac{(\sfa_{01}+\sfa_{10})\,t}{1-\ldots}
}}}}
}
\nonumber \\
\end{eqnarray}
with coefficients
\begin{subeqnarray}
   \alpha_{2k-1}  & = &  \sum_{\xi=0}^{k-1} \sfb_{\xi, k-1-\xi}  \\[1mm]
   \alpha_{2k}    & = &  \sum_{\xi=0}^{k-1} \sfa_{\xi, k-1-\xi}  \\[1mm]
   \delta_{2k-1}  & = &  0    \\[1mm]
   \delta_{2k}    & = &  
     \sum_{\xi=0}^{k-1} \sfc_{\xi,k-1-\xi} \,+\,
     \sum_{\xi=0}^{k-1} \sfd_{\xi,k-1-\xi} \,-\,
     \sum_{\xi=0}^{k-1} \sfa_{\xi,k-1-\xi} \,-\,  
     \sum_{\xi=0}^{k-1} \sfb_{\xi,k-1-\xi}
         \qquad
\label{eq.bt.master.weights.T}
\end{subeqnarray}
\end{theorem}

If we now choose some specialization of the variables $\bsfa,\bsfb,\bsfc,\bsfd$
that satisfies
\be
   \sum_{\xi=0}^{k-1} \sfa_{\xi,k-1-\xi} \,+\,
   \sum_{\xi=0}^{k-1} \sfb_{\xi,k-1-\xi}
   \;=\;
   \sum_{\xi=0}^{k-1} \sfc_{\xi,k-1-\xi} \,+\,
   \sum_{\xi=0}^{k-1} \sfd_{\xi,k-1-\xi}
   \;,
\ee
then obviously the weights (\ref{eq.bt.master.weights.T}c,d)
become $\bdelta = \bzero$, so that we have an S-fraction.
There are many ways that this can be done,
but the simplest is to take $\bsfc = \bsfa$ and $\bsfd = \bsfb$:
that is, we just consider the status of the left child
and ignore the status of the right child.
We then have:

\begin{corollary}
   \label{cor.bt.master.T}
The ordinary generating function of the polynomials
$Q_n(\bsfa, \bsfb,\bsfa,\bsfb)$
has the S-fraction
\begin{eqnarray}
\sum_{n=0}^\infty Q_n(\bsfa, \bsfb,\bsfa,\bsfb) \, t^n
\;=\; 
\cfrac{1}{1 -  \cfrac{\sfb_{00} \, t}{1 - \cfrac{\sfa_{00}\, t}{ 1 - 
\cfrac{(\sfb_{01}+\sfb_{10})\,t }{1- \cfrac{(\sfa_{01}+\sfa_{10})\,t}{1-\ldots}
}}}}
\end{eqnarray}
with coefficients
\begin{subeqnarray}
   \alpha_{2k-1}  & = &  \sum_{\xi=0}^{k-1} \sfb_{\xi, k-1-\xi}  \\[1mm]
   \alpha_{2k}    & = &  \sum_{\xi=0}^{k-1} \sfa_{\xi, k-1-\xi}
\label{eq.cor.bt.master.weights.T}
\end{subeqnarray}
\end{corollary}

\subsection{Increasing restricted ternary trees}

We now generalize the results of the preceding subsection
from increasing binary trees to increasing restricted ternary trees.
We write $\RT_n$ for the set of increasing restricted ternary trees
on the vertex set $[n]$.


\subsubsection{Simple J-fraction and T-fraction}

We first introduce a polynomial in five variables
that enumerates increasing restricted ternary trees
according to the node types:
for $n\geq 0$, we define
\begin{equation}
    P_n(x_1, x_2, y_1, y_2,w) 
    \;\eqdef\;
    \sum_{T\in \RT_n} 
x_1^{I_T(101)} y_1^{I_T(000)} x_2^{I_T(100)} y_2^{I_T(001)} w^{I_T(010)}\;.
\label{eq.Pn.abt}
\end{equation}
Thus, the variables $x_1, y_1, x_2, y_2, w$
are associated to the node types
$101$, $000$, $100$, $001$, $010$, respectively.
In particular, for vertices without a middle child
we have $x_1$ or $x_2$ when a vertex has a left child,
and $y_1$ or $y_2$ when it does not;
we have a subscript 1 (resp.~2)
when the vertex has even (resp.~odd) out-degree;
and finally, we have $w$ for a vertex with a middle child.
By convention $P_0 = 1$ (corresponding to the empty tree);
for $n\geq 1$, the polynomial $P_n$ always has a factor $y_1$
since the vertex $n$ is always a leaf.
When $w=0$, \reff{eq.Pn.abt} reduces to \reff{eq.Pn.bt}.

The polynomials~\reff{eq.Pn.abt} have a beautiful J-fraction:

\begin{theorem}
   \label{thm.Simple.J.abt}
The ordinary generating function of the polynomials
$P_{n+1}(x_1, x_2, y_1, y_2, w)$
has the J-fraction
\be
\sum_{n=0}^\infty P_{n+1}(x_1, x_2, y_1, y_2, w)\, t^n
\;=\;
\Scale[0.9]{
\cfrac{y_1}{1-(x_2 +y_2 +w) t - \cfrac{2 x_1 y_1 \, t^2}{1-2(x_2+y_2+w) t - \cfrac{6 x_1 y_1 t^2}{1- \ldots}}}
}
\ee
with coefficients
\begin{eqnarray}
\gamma_n \;=\; (n+1)\,(x_2+y_2+w) \,,
\qquad 
\beta_n \;=\; n(n+1)\, x_1 y_1\;.
\label{eq.def.simple.abt.J}
\end{eqnarray}
\end{theorem}

We will prove Theorem~\ref{thm.Simple.J.abt} 
in Section~\ref{subsec.bijection.J.abt},
as a special case of the more general ``master'' J-fraction;
we will also give a very simple alternate proof
in Section~\ref{sec.algebraic.proofs}.
Specializing Theorem~\ref{thm.Simple.J.abt} to $w=0$
yields Theorem~\ref{thm.Simple.J.bt}.

Specializing Theorem~\ref{thm.Simple.J.abt}
to $x_1 = x_2 = y_1 = y_2 = w = 1$,
we deduce that the sequence of cardinalities
$(|\RT_{n+1}|)_{n \geq 0}$
has a nice J-fraction:

\begin{corollary}
The ordinary generating function of the sequence
$(|\RT_{n+1}|)_{n \geq 0}$ has the J-fraction
\begin{equation}
\sum_{n=0}^\infty |\RT_{n+1}| \, t^n 
\;=\;
\cfrac{1}{1 - 3  t -  \cfrac{2 t^2}{1 - 6 t   \cfrac{6  t}{ 1 - 9 t -   \cfrac{12 t^2}{1- \ldots}} }}
\end{equation}
with coefficients
\be
\gamma_n \;=\; 3(n+1) \,,
\qquad
\beta_n \;=\; n(n+1) \;.
\ee
\end{corollary}

\bigskip

We now obtain a T-fraction from Theorem~\ref{thm.Simple.J.abt}
by using odd contraction (Proposition~\ref{prop.contraction_odd.Ttype}):

\begin{theorem}
\label{thm.simple.abt.T}
The ordinary generating function of the polynomials
$P_{n}(x_1, x_2, y_1, y_2, w)$
has the T-fraction
\begin{eqnarray}
&&
\hspace*{-15mm}
\sum_{n=0}^\infty P_{n}(x_1, x_2, y_1, y_2, w)\, t^n\;\nonumber\\[2mm]
&&
\hspace*{-15mm}=\;
\Scale[0.9]{
\cfrac{1}{1- \cfrac{y_1 t}{1-(x_2+y_2+w - x_1-y_1 ) t - \cfrac{ x_1 t}{1- \cfrac{2 y_1 t}{1  - 2(x_2+y_2+w - x_1-y_1 ) t - \cfrac{2 x_1 t}{1-\ldots}}}}}
}
\end{eqnarray}
with coefficients
\begin{subeqnarray}
    \alpha_{2k-1} & = & k y_1\\[2mm]
    \alpha_{2k} & = & k x_1\\[2mm]
    \delta_{2k-1} & = & 0\\[2mm]
    \delta_{2k} & = & k\, (x_2 + y_2 +w - x_1-y_1)
\label{eq.def.weights.simple.abt.T}
\end{subeqnarray}
\end{theorem}

\noindent
Indeed, it is straightforward to check that the weights
\reff{eq.def.simple.abt.J} and \reff{eq.def.weights.simple.abt.T}
satisfy~\reff{eq.contraction_odd.coeffs.Ttype}.
Specializing Theorem~\ref{thm.simple.abt.T} to $w=0$ yields
Theorem~\ref{thm.simple.bt.T}.

If we now choose some specialization of $x_1,x_2,y_1,y_2$
that satisfies $x_1 + y_1 = x_2 + y_2$,
then the weight (\ref{eq.def.weights.simple.abt.T}d) simplifies to
$\delta_{2k} = kw$.
In particular, when $x_1 = x_2 = y_1 = y_2 = w = 1$, we obtain:

\begin{corollary}
   \label{cor.abt.numbers}
The ordinary generating function of the sequence
$(|\RT_n|)_{n \ge 0}$ has the T-fraction
\begin{equation}
\sum_{n=0}^\infty |\RT_n| \, t^n 
\;=\;
\cfrac{1}{1 -  \cfrac{ t}{1 - t - \cfrac{t}{ 1 - \cfrac{2t}{1 - 
2  t - \cfrac{2  t}{1 -  \cfrac{3 t}{1 - 3 t- \cfrac{3 t}{1 - \ldots}} }}}}}
\end{equation}
with coefficients
\be
\alpha_{2k-1} \;=\; \alpha_{2k} \;=\; k \,,
\qquad
\delta_{2k-1} \;=\; 0 \,,
\qquad
\delta_{2k} \;=\; k \;.
\ee
\end{corollary}

\noindent
This is the sequence \reff{eq.a=c=0}  \cite[A230008]{OEIS}.

\subsubsection{Master J-fraction and T-fraction}

Fix a consistent tree-traversal algorithm $\bfA$.
Now let
$\bsfa = (\sfa_{\ell,\ell'})_{\ell,\ell'\geq 0}$,
$\bsfb = (\sfb_{\ell,\ell'})_{\ell,\ell'\geq 0}$,
$\bsfc = (\sfc_{\ell,\ell'})_{\ell,\ell'\geq 0}$,
$\bsfd = (\sfd_{\ell,\ell'})_{\ell,\ell'\geq 0}$,
$\bsff = (\sff_{\ell,\ell'})_{\ell,\ell'\geq 0}$
be infinite sets of indeterminates,
and define polynomials $Q_n(\bsfa, \bsfb,\bsfc,\bsfd,\bsff)$ by
\begin{eqnarray}
Q_n(\bsfa, \bsfb,\bsfc,\bsfd,\bsff)
&=&
\sum_{T\in \RT_n} \:
\prod_{N(v,T) =101 }
\sfa_{\croix(v,T),\nid(v,T)}
\prod_{N(v,T) =000 }
\sfb_{\croix(v,T),\nid(v,T)} \times
\nonumber\\
&& \qquad\quad\, \prod_{N(v,T) =100 }
\sfc_{\croix(v,T),\nid(v,T)}
\prod_{ N(v,T) =001 }
\sfd_{\croix(v,T),\nid(v,T)}\times
\nonumber\\
&& \qquad\quad\, \prod_{N(v,T) =010 }
\sff_{\croix(v,T),\nid(v,T)}
\;.
\label{eq.def.Qn.abt}
\end{eqnarray}
Of course $P_0 = 1$;
for $n\geq 1$, the polynomial $Q_n$ always has a factor $\sfb_{00}$
since the vertex $n$ is always a leaf.

We have the following master J-fraction:

\begin{theorem}
   \label{thm.abt.master.J}
The ordinary generating function of the polynomials
$Q_{n+1}(\bsfa, \bsfb,\bsfc,\bsfd,\bsff)$
has the J-fraction
\begin{eqnarray}
&& 
\hspace*{-15mm}
\sum_{n=0}^\infty Q_{n+1}(\bsfa, \bsfb,\bsfc,\bsfd,\bsff) \, t^n
\;=\; 
\nonumber \\[2mm]
&&
\hspace*{-15mm}
\Scale[0.8]{
\cfrac{\sfb_{00}}{1 - (\sfc_{00}+\sfd_{00}+\sff_{00})\, t - \cfrac{\sfa_{00}(\sfb_{01} + \sfb_{10}) \, t^2}{1 - (\sfc_{01}+\sfc_{10}+\sfd_{01}+\sfd_{10}+\sff_{01}+\sff_{10})\, t - \cfrac{(\sfa_{01}+\sfa_{10})(\sfb_{02} +\sfb_{11} + \sfb_{20})\, t^2}{ 1 - \ldots}}}
}
\label{eq.thm.abt.master.J}
\end{eqnarray}
where the coefficients are defined as follows:
\begin{subeqnarray}
    \beta_{n} & = & 
     \biggl(
     \sum_{\xi=0}^{n-1} \sfa_{\xi, n-1-\xi}
     \biggr)
     \biggl(
     \sum_{\xi=0}^{n} \sfb_{\xi, n-\xi}
     \biggr)
    \\[2mm] 
    \gamma_{n} & = & \sum_{\xi=0}^n \sfc_{\xi,n-\xi} \, + \, 
    \sum_{\xi=0}^n \sfd_{\xi,n-\xi}
    \, + \, 
    \sum_{\xi=0}^n \sff_{\xi,n-\xi}
\label{eq.abt.master.J.weights}
\end{subeqnarray}
\end{theorem}
We will prove Theorem~\ref{thm.abt.master.J}
in Section~\ref{subsec.bijection.J.abt},
by bijection onto a suitable class of labeled Motzkin paths.
When $\sff_{\ell,\ell'} = 0$ for all $\ell,\ell'\geq 0$,
this yields Theorem~\ref{thm.bt.master}.

Once again we can obtain a T-fraction
by using odd contraction (Proposition~\ref{prop.contraction_odd.Ttype}):

\begin{theorem}
\label{thm.abt.master.T}
The ordinary generating function of the polynomials
$Q_n(\bsfa, \bsfb,\bsfc,\bsfd,\bsff)$
has the T-fraction
\begin{eqnarray}
&& 
\hspace*{-10mm}
\sum_{n=0}^\infty Q_n(\bsfa, \bsfb,\bsfc,\bsfd,\bsff) \, t^n
\;=\; 
\nonumber \\[2mm]
&&
\hspace*{-10mm}
\Scale[0.77]{\cfrac{1}{1 -  \cfrac{\sfb_{00} \, t}{1 - 
(\sfc_{00}\!+\!\sfd_{00}\!+\!\sff_{00}\!-\!\sfa_{00}\!-\!\sfb_{00})\, t
 - \cfrac{\sfa_{00}\, t}{ 1 - 
\cfrac{(\sfb_{01}\!+\!\sfb_{10})\,t }{1-
(\sfc_{01}\!+\!\sfc_{10}\!+\!\sfd_{01}\!+\!\sfd_{10}\!+\!\sff_{01}\!+\!\sff_{10}\!-\!\sfa_{01}\!-\!\sfa_{10}\!-\!\sfb_{01}\!-\!\sfb_{10})
 \,t -
\cfrac{(\sfa_{01}\!+\!\sfa_{10})\,t}{1-\ldots}
}}}}
}
\nonumber \\[2mm]
\end{eqnarray}
with coefficients
\begin{subeqnarray}
   \alpha_{2k-1}  & = &  \sum_{\xi=0}^{k-1} \sfb_{\xi, k-1-\xi}  \\[1mm]
   \alpha_{2k}    & = &  \sum_{\xi=0}^{k-1} \sfa_{\xi, k-1-\xi}  \\[1mm]
   \delta_{2k-1}  & = &  0    \\[1mm]
   \delta_{2k}    & = &  
     \sum_{\xi=0}^{k-1} \sfc_{\xi,k-1-\xi} \,+\,
     \sum_{\xi=0}^{k-1} \sfd_{\xi,k-1-\xi} \,+\,
     \sum_{\xi=0}^{k-1} \sff_{\xi,k-1-\xi} \,-\,
     \sum_{\xi=0}^{k-1} \sfa_{\xi,k-1-\xi} \,-\,  
     \sum_{\xi=0}^{k-1} \sfb_{\xi,k-1-\xi}
         \nonumber \\
\label{eq.abt.master.weights.T}
\end{subeqnarray}
\end{theorem}

If we now choose some specialization of the variables $\bsfa,\bsfb,\bsfc,\bsfd$
that satisfies
\be
   \sum_{\xi=0}^{k-1} \sfa_{\xi,k-1-\xi} \,+\,
   \sum_{\xi=0}^{k-1} \sfb_{\xi,k-1-\xi}
   \;=\;
   \sum_{\xi=0}^{k-1} \sfc_{\xi,k-1-\xi} \,+\,
   \sum_{\xi=0}^{k-1} \sfd_{\xi,k-1-\xi}
   \;,
\ee
then obviously the weight (\ref{eq.abt.master.weights.T}d) simplifies to
\be
   \delta_{2k}  \;=\;  \sum_{\xi=0}^{k-1} \sff_{\xi,k-1-\xi}
   \;.
\ee
There are many ways that this can be done,
but the simplest is to take $\bsfc = \bsfa$ and $\bsfd = \bsfb$:
that is, if the vertex does not have a middle child,
we just consider the status of the left child
and ignore the status of the right child.
We then have:

\begin{corollary}
\label{cor.abt.master.T}
The ordinary generating function of the polynomials
$Q_n(\bsfa, \bsfb,\bsfa,\bsfb,\bsff)$
has the T-fraction
\begin{eqnarray}
\sum_{n=0}^\infty Q_n(\bsfa, \bsfb,\bsfa,\bsfb,\bsff) \, t^n
\;=\; 
\cfrac{1}{1 -  \cfrac{\sfb_{00} \, t}{1 - \sff_{00}\, t - \cfrac{\sfa_{00}\, t}{ 1 - 
\cfrac{(\sfb_{01}+\sfb_{10})\,t }{1-
(\sff_{01}+\sff_{10}) \,t -
\cfrac{(\sfa_{01}+\sfa_{10})\,t}{1-\ldots}
}}}}
\end{eqnarray}
with coefficients
\begin{subeqnarray}
   \alpha_{2k-1}  & = &  \sum_{\xi=0}^{k-1} \sfb_{\xi, k-1-\xi}  \\[1mm]
   \alpha_{2k}    & = &  \sum_{\xi=0}^{k-1} \sfa_{\xi, k-1-\xi}  \\[1mm]
   \delta_{2k-1}  & = &  0    \\[1mm]
   \delta_{2k}    & = &  
     \sum_{\xi=0}^{k-1} \sff_{\xi,k-1-\xi} \nonumber \\
\label{eq.cor.abt.master.weights.T}
\end{subeqnarray}
\end{corollary}

\subsubsection{A more general master T-fraction}
   \label{subsubsec.moregeneral}

But there is a more general way than Corollary~\ref{cor.abt.master.T}
to obtain a T-fraction from the J-fraction of Theorem~\ref{thm.abt.master.J},
which we introduce now because it foreshadows in simpler form
what we will do in Section~\ref{subsubsec.IRT.master}
for interval-labeled trees.
Rather than take $\bsfc = \bsfa$ and $\bsfd = \bsfb$,
we introduce indeterminates
$\bsfahat = (\sfahat_{\ell,\ell'})_{\ell,\ell'\geq 0}$
and $\bsfbhat = (\sfbhat_{\ell,\ell'})_{\ell,\ell'\geq 0}$
with two subscripts,
and indeterminates
$\bmu = (\mu_{\ell})_{\ell\geq 0}$ and $\bnu = (\nu_{\ell})_{\ell\geq 0}$
with one subscript,
and then specialize the formulae of the preceding subsection to
\begin{subeqnarray}
    \sfa_{\ell,\ell'} & = & \sfahat_{\ell,\ell'} \, \mu_{\ell+\ell'+1}\\[2mm]
    \sfb_{\ell,\ell'} & = & \sfbhat_{\ell,\ell'} \, \nu_{\ell+\ell'-1}\\[2mm]
    \sfc_{\ell,\ell'} & = & \sfahat_{\ell,\ell'} \, \nu_{\ell+\ell'}\\[2mm]
    \sfd_{\ell,\ell'} & = & \sfbhat_{\ell,\ell'} \, \mu_{\ell+\ell'}
\label{eq.specializations.J.to.T.abt}
\end{subeqnarray}
That is, the weights $\bsfahat$ and $\bsfbhat$
concern the status of the left child,
taking full account of both croix and nid;
while the weights $\bmu$ and $\bnu$
concern the status of the right child,
but taking account only of lev = croix + nid.
Another way of saying this is that we have assigned vertex weights
\begin{equation}
\wt(v)
\;=\;
\begin{cases}
   \sfahat_{\croix(v,T),\nid(v,T)}\, \mu_{\lev(v,T)+1}
      \quad\quad \text{if $N(v,T)=101$}\\[2mm]
   \sfbhat_{\croix(v,T),\nid(v,T)}\, \nu_{\lev(v,T)-1}
      \qquad\, \text{if $N(v,T)=000$}\\[2mm]
   \sfahat_{\croix(v,T),\nid(v,T)}\, \nu_{\lev(v,T)}
      \qquad\,\quad\text{if $N(v,T)=100$}\\[2mm]
   \sfbhat_{\croix(v,T),\nid(v,T)}\, \mu_{\lev(v,T)}
      \qquad\quad \text{if $N(v,T)=001$}\\[2mm]
   \sff_{\croix(v,T),\nid(v,T)}
      \qquad\quad\qquad\quad\:\:\: \text{if $N(v,T)=010$}
\end{cases}
\label{eq.weight.master.abt.T}
\end{equation}
When $\bmu = \bnu = \bone$,
this reduces to the previously considered case
$\bsfc = \bsfa$ and $\bsfd = \bsfb$.

Using the weights~\reff{eq.weight.master.abt.T},
we then define polynomials $Q_n^\star(\bsfahat,\bsfbhat,\bmu,\bnu,\bsff)$ by
\begin{subeqnarray}
Q_0^\star(\bsfahat,\bsfbhat,\bmu,\bnu,\bsff)
&=& 
1
   \slabel{eq.def.Qn.abt.T.a}
   \\[2mm]
Q_n^\star(\bsfahat,\bsfbhat,\bmu,\bnu,\bsff)
&=& 
\mu_0 \, \sfbhat_{00}
\sum_{T\in \RT_n}
\prod_{v=1}^{n-1}
\wt(v) \qquad \hbox{for $n\geq 1$}
   \slabel{eq.def.Qn.abt.T.b}
   \label{eq.def.Qn.abt.T}
\end{subeqnarray}
Note that vertex $n$ is not given a weight;
instead we have the prefactor $\mu_0 \, \sfbhat_{00}$.
Note also that any leaf $v \le n-1$ must have $\lev(v,T) \ge 1$,
since the parent of vertex $v+1$ must be $< v$;
this means that the subscript on $\nu$ is always $\ge 0$,
as it should be.

We then have the following theorem:

\begin{theorem}
   \label{thm.abt.master.1}
The ordinary generating function of the polynomials
$Q_n^\star(\bsfahat,\bsfbhat,\bmu,\bnu,\bsff)$
has the T-fraction
\be
\sum_{n=0}^\infty Q_n^\star(\bsfahat,\bsfbhat,\bmu,\bnu,\bsff) \, t^n
\;=\;
\cfrac{1}{1 - \cfrac{\mu_0\,\sfbhat_{00} t}{1 - \sff_{00} t - \cfrac{\nu_0\,\sfahat_{00} t}{ 1 - \cfrac{\mu_1(\sfbhat_{01}+\sfbhat_{10}) t}{1 - (\sff_{01}+\sff_{10}) t - \cfrac{\nu_1(\sfahat_{01}+\sfahat_{10}) t}{1-\ldots}}}}}
\ee
with coefficients
\begin{subeqnarray}
    \alpha_{2k-1} & = &
       \mu_{k-1} \left(\sum_{\xi =0}^{k-1} \sfbhat_{\xi, k-1-\xi}\right)
       \slabel{eq.thm.abt.master.1.weights.a}  \\[2mm] 
    \alpha_{2k} & = & \nu_{k-1} \left(\sum_{\xi =0}^{k-1} \sfahat_{\xi, k-1-\xi} \right) \\[2mm]
    \delta_{2k-1} & = & 0 \qquad\qquad\qquad\qquad\\[2mm]
    \delta_{2k} & = & \sum_{\xi =0}^{k-1} \sff_{\xi, k-1-\xi}
 \label{eq.thm.abt.master.1.weights}
\end{subeqnarray}
\end{theorem}

\noindent
When $\bmu = \bnu = \bone$,
this reduces to Corollary~\ref{cor.abt.master.T}.

\bigskip

\noindent{\sc Proof of Theorem~\ref{thm.abt.master.1},
   assuming Theorem~\ref{thm.abt.master.J}.\ }
We first define polynomials 
$\widehat{Q}_n(\bsfahat,\bsfbhat,\bmu,\bnu,\bsff)$
by making the substitutions \reff{eq.specializations.J.to.T.abt}
in the polynomials~\reff{eq.def.Qn.abt}.
Since for $n \ge 1$ this gives the vertex $n$ a weight
$\sfb_{00} = \sfbhat_{00} \nu_{-1}$,
we have
\be
Q^\star_n(\bsfahat,\bsfbhat,\bmu,\bnu,\bsff)
\;=\; 
\mu_0 \, (\nu_{-1})^{-1}\,
\widehat{Q}_n(\bsfahat,\bsfbhat,\bmu,\bnu,\bsff)
\quad\hbox{for $n \ge 1$}
\;.
\ee
Now, making the substitutions~\reff{eq.specializations.J.to.T.abt}
in Theorem~\ref{thm.abt.master.J} shows that
$\sum\limits_{n=0}^\infty \widehat{Q}_{n+1}(\bsfahat,\bsfbhat,\bmu,\bnu,\bsff)$
is given by the J-fraction \reff{eq.thm.abt.master.J}
with $\sfb_{00} = \sfbhat_{00} \nu_{-1}$.
Therefore,
$\sum\limits_{n=0}^\infty Q^\star_{n+1}(\bsfahat,\bsfbhat,\bmu,\bnu,\bsff)$
is given by the same J-fraction
but with the top coefficient $\sfb_{00}$
replaced by $\mu_0 \sfbhat_{00}$.
Applying odd contraction (Proposition~\ref{prop.contraction_odd.Ttype})
to this latter J-fraction,
with $\balpha$ and $\bdelta$ defined by \reff{eq.thm.abt.master.1.weights},
proves Theorem~\ref{thm.abt.master.1}:
the key point is that \reff{eq.thm.abt.master.1.weights.a}
gives $\alpha_1 = \mu_0 \sfbhat_{00}$.
\qed

We will see in Section~\ref{subsubsec.IRT.master}
that Theorem~\ref{thm.abt.master.1}
is a special case of a more general result (Theorem~\ref{thm.sbt.master})
for increasing interval-labeled restricted ternary trees
(for which we will provide a bijective proof
 in Section~\ref{subsec.sbt.bijection}).

\subsection{Increasing interval-labeled restricted ternary trees}

Finally, we state our results for
increasing interval-labeled restricted ternary trees,
which are counted by the sequence \reff{eq.c=0}.
We write $\IRT_n$ for the set of 
increasing interval-labeled restricted ternary trees
on the label set $[0,n]$.

\subsubsection{Simple T-fraction}

We first introduce a polynomial in six variables
that enumerates increasing interval-labeled restricted ternary trees
on the label set $[0,n]$
according to the node types and the label surplus:
\begin{equation}
    P_n(x_1, x_2, y_1, y_2,w, z) 
    \;\eqdef\;
    \sum_{T\in \IRT_n} 
x_1^{I'_T(101)} y_1^{I'_T(000)} x_2^{I'_T(100)} y_2^{I'_T(001)} w^{I'_T(010)} z^{I_T(\varepsilon)}\;.
\label{eq.Pn.sbt}
\end{equation}

Please note that we are here using the weights $I'_T(s)$
defined in \reff{eq.def.node.count.prime}:
that is, we are giving node-type weights only to {\em non-root}\/ vertices.
(Of course, using $I'_T(s)$ instead of $I_T(s)$ matters only
 for the node types 000 and 100,
 since these are the only possible node types for the root
 in an increasing interval-labeled restricted ternary tree.)
As before, $I_T(\varepsilon)$ is the label surplus of the tree $T$,
as defined in \reff{eq.def.node.count.varep}/\reff{eq.label.surplus.dist}.
Thus, the variables $x_1, y_1, x_2, y_2, w$
are associated to the node types $101$, $000$, $100$, $001$, $010$,
respectively,
while $z$ is associated to the label surplus.
Note that $P_n$ is homogeneous of degree $n$
in the six variables $x_1, x_2, y_1, y_2,w, z$.

The polynomials~\reff{eq.Pn.sbt} do not have a nice J-fraction;
and as far as we know they do not have a nice T-fraction either.
However, under the two specializations $x_1=x_2$ and
$y_1=y_2$, they have a beautiful T-fraction: 

\begin{theorem}
   \label{thm.sbt.simple}
The ordinary generating function of the polynomial \reff{eq.Pn.sbt} 
specialized to $x_1 = x_2 = x$ and $y_1 = y_2 = y$
has the T-fraction
\begin{eqnarray}
   & &
   \sum_{n=0}^\infty  P_n(x, x, y, y,w,z) \, t^n
           \nonumber\\
   & &
   \qquad \;=\;
\Scale[0.9]{\cfrac{1}{1 - z t -  \cfrac{y t}{1 - w t -  \cfrac{x  t}{ 1 - z t - \cfrac{2y t}{1 - 2w t - \cfrac{2x  t}{1 - z t - \cfrac{3y t}{1 - 3w t - \cfrac{3x t}{1 - \ldots}} }}}}}}
   \qquad\qquad
\end{eqnarray}
with coefficients:
\begin{subeqnarray}
    \alpha_{2k-1} & = & k y \\[2mm] 
    \alpha_{2k} & = & k x \\[2mm]
    \delta_{2k-1} & = & z  \\[2mm]
    \delta_{2k} & = & k w
\label{eq.sbt.simple.weights}
\end{subeqnarray}

\end{theorem}

We will prove Theorem~\ref{thm.sbt.simple}
in Section~\ref{subsec.sbt.bijection},
as a special case of the more general ``master'' T-fraction;
we will also give a very simple alternate proof
in Section~\ref{subsec.algebraic.proofs.IRT}.

%
%
%

\bigskip
\bigskip

Setting $x=y=w=z=1$, we obtain the following simple corollary:

\begin{corollary}
\label{cor.sbt.numbers}
The ordinary generating function of the sequence
$(a_n)_{n\geq 0}$ where $a_n = |\IRT_n|$
has the T-fraction
\be
\sum_{n=0}^\infty a_n t^n 
\;=\;
\Scale[0.95]{
\cfrac{1}{1 -  t -  \cfrac{ t}{1 -  t -  \cfrac{1  t}{ 1 -  t - \cfrac{2 t}{1 - 2 t - \cfrac{2  t}{1 -  t - \cfrac{3 t}{1 - 3 t - \cfrac{3 t}{1 - \ldots}} }}}}}
}
\ee
with coefficients
\be
\alpha_{2k-1} \;=\; \alpha_{2k} \;=\; \delta_{2k} \;=\; k \;, \qquad
  \delta_{2k-1} \;=\; 1\;.
\ee
\end{corollary}

This is the sequence shown in~\reff{eq.c=0}.

\medskip

In fact, we can further generalize Theorem~\ref{thm.sbt.simple}
to have only one specialization instead of two; we state this now:

\begin{theorem}
   \label{thm.sbt.simple.new}
The ordinary generating function of the polynomial \reff{eq.Pn.sbt}
specialized to $y_2 = x_1 + y_1 - x_2$
has the T-fraction
\begin{eqnarray}
   & &
   \sum_{n=0}^\infty  P_n(x_1, x_2, y_1, y_2,w,z) \, t^n
           \nonumber\\
   & &
   \qquad \;=\;
\Scale[0.9]{\cfrac{1}{1 - z t -  \cfrac{y_1 t}{1 - w t -  \cfrac{x_1  t}{ 1 - z t - \cfrac{2y_1 t}{1 - 2w t - \cfrac{2x_1  t}{1 - z t - \cfrac{3y_1 t}{1 - 3w t - \cfrac{3x_1 t}{1 - \ldots}} }}}}}}
   \qquad\qquad
\end{eqnarray}
with coefficients:
\begin{subeqnarray}
    \alpha_{2k-1} & = & k y_1 \\[2mm]
    \alpha_{2k} & = & k x_1 \\[2mm]
    \delta_{2k-1} & = & z  \\[2mm]
    \delta_{2k} & = & k w
\label{eq.sbt.simple.weights.new}
\end{subeqnarray}
\end{theorem}

Note that the result does not depend on the variable $x_2$. 
In Section~\ref{subsec.algebraic.proofs.IRT}
we will give an algebraic proof of Theorem~\ref{thm.sbt.simple.new},
based on Proposition~\ref{prop.transformation}.

\subsubsection{Master T-fraction}
   \label{subsubsec.IRT.master}

We now go farther, and introduce a polynomial in six infinite families
of indeterminates:
$\bsfahat = (\sfahat_{\ell,\ell'})_{\ell,\ell'\geq 0}$,
$\bsfbhat = (\sfbhat_{\ell,\ell'})_{\ell,\ell'\geq 0}$
and $\bsff = (\sff_{\ell,\ell'})_{\ell,\ell'\geq 0}$
with two subscripts,
and
$\bmu = (\mu_{\ell})_{\ell\geq 0}$,
$\bnu = (\nu_{\ell})_{\ell\geq 0}$
and $\bsfe = (\sfe_{\ell})_{\ell\geq 0}$
with one subscript.
The notation is thus the same as in Section~\ref{subsubsec.moregeneral},
together with the new indeterminates $\bsfe$.

Let $T\in \IRT_n$ be an increasing interval-labeled restricted ternary tree
on the label set $[0,n]$,
and let $v$ be a vertex of $T$ with label set $L_v = \{l,l+1, \ldots, l+j\}$;
here $j = |L_v| - 1$ is the label surplus of the vertex $v$.
To each vertex $v$ we assign a weight $\wt(v)$ as follows:
\begin{eqnarray}
   \hspace*{-4cm}
   & &
   \bullet\;
   \textrm{If $l=0$ and $l+j=n$, we assign weight} \;
   \wt(v) \,=\, \sfe_{0}^n \;.
      \label{eq.weight.master.sbt.I} \\[2mm]
   \hspace*{-4cm}
   & &
   \bullet\;
   \textrm{If $l=0$ and $l+j<n$, we assign weight} \;
   \wt(v) \,=\, \mu_0 \sfe_{0}^j \;.
      \label{eq.weight.master.sbt.II} \\[2mm]
   \hspace*{-4cm}
   & &
   \bullet\;
   \textrm{If $l>0$ and $l+j=n$, we assign weight} \;
   \wt(v) \,=\, \sfbhat_{00} \sfe_{0}^j \;.
      \label{eq.weight.master.sbt.III}  \\[2mm]
   \hspace*{-4cm}
   & &
   \bullet\;
   \textrm{If $l>0$ and $l+j<n$, we assign weight} \;
        \nonumber
\end{eqnarray}
\begin{equation}
\wt(v)
\;=\;
\begin{cases}
\sfahat_{\croix(v,T),\nid(v,T)}\, \mu_{\lev(v,T)+1}\, (\sfe_{\lev(v,T)+1})^j
\quad\; \text{if $N(v,T)=101$}\\[2mm]
\sfbhat_{\croix(v,T),\nid(v,T)}\, \nu_{\lev(v,T)-1}\, (\sfe_{\lev(v,T)})^j 
\qquad \text{if $N(v,T)=000$}\\[2mm]
\sfahat_{\croix(v,T),\nid(v,T)}\, \nu_{\lev(v,T)}\, (\sfe_{\lev(v,T)+1})^j
\qquad\, \text{if $N(v,T)=100$}\\[2mm]
\sfbhat_{\croix(v,T),\nid(v,T)}\, \mu_{\lev(v,T)}\, (\sfe_{\lev(v,T)})^j
\qquad\quad \text{if $N(v,T)=001$}\\[2mm]
\sff_{\croix(v,T),\nid(v,T)}
\qquad\qquad\qquad\qquad\qquad\quad \text{if $N(v,T)=010$}\end{cases}
\label{eq.weight.master.sbt}
\end{equation}
(The reasons for these weights will be seen in
 Section~\ref{subsec.sbt.bijection} in the context of the bijective proof.)

Thus, the case \reff{eq.weight.master.sbt.I}
corresponds to a {\em trivial tree}\/
(i.e., a tree consisting only of the root),
in which the root $v$ has label surplus $n$.
The case \reff{eq.weight.master.sbt.II}
corresponds to $v$ being the root of a nontrivial tree.
The case \reff{eq.weight.master.sbt.III}
corresponds to $v$ being the {\em final vertex}\/
(i.e., the vertex having $n$ in its label set)
of a nontrivial tree.
And finally, the case \reff{eq.weight.master.sbt}
corresponds to $v$ being neither the root nor the final vertex.

Notice that a vertex $v$ with no middle child
gets a letter $\sfahat$ if it has a left child, and $\sfbhat$ if not,
with subscripts indicating croix and nid;
it also gets a letter $\mu$ if it has a right child, and $\nu$ if not,
with a single subscript indicating lev = croix + nid.
A vertex $v$ with a middle child gets a letter $\sff$.
And finally, a vertex with $|L_v| > 1$
gets a letter $\sfe$ raised to the power $|L_v| - 1$
(which is the label surplus of $v$).

Note also that if $l+j < n$ and $N(v,T)=000$ (i.e., $v$ is a leaf),
then we necessarily have $\lev(v,T) \ge 1$,
since some vertex higher-numbered than $v$
must have a parent that is lower-numbered than $v$;
then the subscript on $\nu$ is $\ge 0$.

Now define polynomials $Q_n(\bsfahat,\bsfbhat,\bmu,\bnu,\bsfe,\bsff)$  as
\be
Q_n(\bsfahat,\bsfbhat,\bmu,\bnu,\bsfe,\bsff)
\;=\;
\sum_{T\in \IRT_n}\prod_{v\in V(T)} \wt(v)\;.
\label{eq.def.Qn.irtt}
\ee
For instance, the first few $Q_n$ (cf.\ Figure~\ref{fig_IRT_n=0,1,2}) are
\begin{subeqnarray}
   Q_0  & = &   1  \\[1mm]
   Q_1  & = &   \mu_0 \sfbhat_{00} \,+\, \sfe_0  \\[1mm]
   Q_2  & = &  (\mu_0 \sfbhat_{00} + \sfe_0)^2 \:+\:
               \mu_0 \sfbhat_{00} (\nu_0 \sfahat_{00} + \sff_{00})
\end{subeqnarray}




We have the following theorem:

\begin{theorem}
   \label{thm.sbt.master}
The ordinary generating function of the polynomials
$Q_n(\bsfahat,\bsfbhat,\bmu,\bnu,\bsfe,\bsff)$
has the T-fraction
\begin{eqnarray}
&&
\hspace*{-15mm}
\sum_{n=0}^\infty Q_n(\bsfahat,\bsfbhat,\bmu,\bnu,\bsfe,\bsff) \, t^n
\;=\;
\nonumber \\
&&
\hspace*{-10mm}
\cfrac{1}{1 - \sfe_{0} t - \cfrac{\mu_0\,\sfbhat_{00} t}{1 - \sff_{00} t - \cfrac{\nu_0\,\sfahat_{00} t}{ 1 - \sfe_{1} t - \cfrac{\mu_1(\sfbhat_{01}+\sfbhat_{10}) t}{1 - (\sff_{01}+\sff_{10}) t - \cfrac{\nu_1(\sfahat_{01}+\sfahat_{10}) t}{1-\ldots}}}}}
\end{eqnarray}
where the coefficients are defined as follows:
\begin{subeqnarray}
    \alpha_{2k-1} & = & \mu_{k-1} \,
       \left(\sum_{\xi =0}^{k-1} \sfbhat_{\xi, k-1-\xi}\right)\\[2mm] 
    \alpha_{2k}   & = & \nu_{k-1} \,
       \left(\sum_{\xi =0}^{k-1} \sfahat_{\xi, k-1-\xi} \right) \\[2mm]
    \delta_{2k-1} & = & \sfe_{k-1}   \\[2mm]
    \delta_{2k} & = & \sum_{\xi =0}^{k-1} \sff_{\xi, k-1-\xi}
\label{eq.sbt.master.weights}
\end{subeqnarray}
\end{theorem}

We will prove Theorem~\ref{thm.sbt.master}
in Section~\ref{subsec.sbt.bijection},
by bijection onto a suitable class of labeled Schr\"oder paths.

Note that when $\sfe_\ell = 0$ for all $\ell$,
the tree $T$ is forced to be single-labeled,
and we simply have the a restricted ternary tree $T'$
on the vertex set $[n]$ together with a root vertex~0
that has (when $n \ge 1$) the vertex~1 (which is the root of $T'$)
as its left child.
When $n=0$, the tree $T$ (which consists solely of the root~0)
gets a weight $Q_0 = 1$, which also equals $Q^\star_0 = 1$
from \reff{eq.def.Qn.abt.T.a}.
When $n \ge 1$, the root of $T$ (vertex~0) gets weight $\mu_0$
from \reff{eq.weight.master.sbt.II},
and the vertex~$n$ gets weight $\sfbhat_{00}$
from \reff{eq.weight.master.sbt.III};
all other vertices get the same weight \reff{eq.weight.master.sbt}
as in \reff{eq.weight.master.abt.T}.
This gives the same prefactor $\mu_0 \, \sfbhat_{00}$
as in \reff{eq.def.Qn.abt.T.b}.
It follows that
\be
   Q_n(\bsfahat,\bsfbhat,\bmu,\bnu,\bzero,\bsff)
   \;=\;
   Q_n^\star(\bsfahat,\bsfbhat,\bmu,\bnu,\bsff)
   \;,
\ee
so that Theorem~\ref{thm.abt.master.1}
is the special case $\bsfe = \bzero$ of Theorem~\ref{thm.sbt.master}.





\section{Preliminaries for the proofs}   \label{sec.prelimproofs}

Our proofs are based on Flajolet's \cite{Flajolet_80}
combinatorial interpretation of continued fractions
in terms of Dyck and Motzkin paths
and its generalization
\cite{Fusy_15,Oste_15,Josuat-Verges_18,Sokal_totalpos,Elvey-Price-Sokal_wardpoly}
to Schr\"oder paths,
together with some bijections mapping our tree models
to labeled Motzkin or Schr\"oder paths.
In Sections~\ref{subsec.prelimproofs.1} and \ref{subsec.prelimproofs.2}
we review briefly these two ingredients and fix our notation.
At the end of Section~\ref{subsec.prelimproofs.1}
we will also give a combinatorial interpretation/proof
of Proposition~\ref{prop.transformation}.

\subsection{Combinatorial interpretation of continued fractions}
   \label{subsec.prelimproofs.1}

Recall that a \textbfit{Motzkin path} of length $n \ge 0$
is a path $\omega = (\omega_0,\ldots,\omega_n)$
in the right quadrant $\N \times \N$,
starting at $\omega_0 = (0,0)$ and ending at $\omega_n = (n,0)$,
whose steps $s_j = \omega_j - \omega_{j-1}$
are $(1,1)$ [``rise'' or ``up step''], $(1,-1)$ [``fall'' or ``down step'']
or $(1,0)$ [``level step''].
We write $h_j$ for the \textbfit{height} of the Motzkin path at abscissa~$j$,
i.e.\ $\omega_j = (j,h_j)$;
note in particular that $h_0 = h_n = 0$.
We write $\scrm_n$ for the set of Motzkin paths of length~$n$,
and $\scrm =  \bigcup_{n=0}^\infty \scrm_n$.
A Motzkin path is called a \textbfit{Dyck path} if it has no level steps.
A Dyck path always has even length;
we write $\scrd_{2n}$ for the set of Dyck paths of length~$2n$,
and $\scrd = \bigcup_{n=0}^\infty \scrd_{2n}$.

Let ${\bf a} = (a_i)_{i \ge 0}$, ${\bf b} = (b_i)_{i \ge 1}$
and ${\bf c} = (c_i)_{i \ge 0}$ be indeterminates;
we will work in the ring $\Z[[{\bf a},{\bf b},{\bf c}]]$
of formal power series in these indeterminates.
To each Motzkin path $\omega$ we assign a weight
$W(\omega) \in \Z[{\bf a},{\bf b},{\bf c}]$
that is the product of the weights for the individual steps,
where a rise starting at height~$i$ gets weight~$a_i$,
a~fall starting at height~$i$ gets weight~$b_i$,
and a level step at height~$i$ gets weight~$c_i$.
Flajolet \cite{Flajolet_80} showed that
the generating function of Motzkin paths
can be expressed as a continued fraction:

\begin{theorem}[Flajolet's master theorem]
   \label{thm.flajolet}
We have
\be
   \sum_{\omega \in \scrm}  W(\omega)
   \;=\;
   \cfrac{1}{1 - c_0 - \cfrac{a_0 b_1}{1 - c_1 - \cfrac{a_1 b_2}{1- c_2 - \cfrac{a_2 b_3}{1- \cdots}}}}
 \label{eq.thm.flajolet}
\ee
as an identity in $\Z[[{\bf a},{\bf b},{\bf c}]]$.
\end{theorem}

In particular, if $a_{i-1} b_i = \beta_i t^2$ and $c_i = \gamma_i t$
(note that the parameter $t$ is conjugate to the length of the Motzkin path),
we have
\be
   \sum_{n=0}^\infty t^n \sum_{\omega \in \scrm_n}  W(\omega)
   \;=\;
   \cfrac{1}{1 - \gamma_0 t - \cfrac{\beta_1 t^2}{1 - \gamma_1 t - \cfrac{\beta_2 t^2}{1 - \cdots}}}
   \;\,,
 \label{eq.flajolet.motzkin}
\ee
so that the generating function of Motzkin paths with height-dependent weights
is given by the J-fraction \reff{def.Jtype}.
Similarly, if $a_{i-1} b_i = \alpha_i t$ and $c_i = 0$
(note that $t$ is now conjugate to the semi-length of the Dyck path), we have
\be
   \sum_{n=0}^\infty t^n \sum_{\omega \in \scrd_{2n}}  W(\omega)
   \;=\;
   \cfrac{1}{1 - \cfrac{\alpha_1 t}{1 - \cfrac{\alpha_2 t}{1 - \cdots}}}
   \;\,,
 \label{eq.flajolet.dyck}
\ee
so that the generating function of Dyck paths with height-dependent weights
is given by the S-fraction \reff{def.Stype}.

Let us now show how to handle Schr\"oder paths within this framework.
A \textbfit{Schr\"oder path} of length $2n$ ($n \ge 0$)
is a path $\omega = (\omega_0,\ldots,\omega_{2n})$
in the right quadrant $\N \times \N$,
starting at $\omega_0 = (0,0)$ and ending at $\omega_{2n} = (2n,0)$,
whose steps are $(1,1)$ [``rise'' or ``up step''],
$(1,-1)$ [``fall'' or ``down step'']
or $(2,0)$ [``long level step''].
We write $s_j$ for the step starting at abscissa $j-1$.
If the step $s_j$ is a rise or a fall,
we set $s_j = \omega_j - \omega_{j-1}$ as before.
If the step $s_j$ is a long level step,
we set $s_j = \omega_{j+1} - \omega_{j-1}$ and leave $\omega_j$ undefined;
furthermore, in this case there is no step $s_{j+1}$.
We write $h_j$ for the height of the Schr\"oder path at abscissa~$j$
whenever this is defined, i.e.\ $\omega_j = (j,h_j)$.
Please note that $\omega_{2n} = (2n,0)$ and $h_{2n} = 0$
are always well-defined,
because there cannot be a long level step starting at abscissa $2n-1$.
Note also that a long level step at even (resp.~odd) height
can occur only at an odd-numbered (resp.~even-numbered) step.
We write $\scrs_{2n}$ for the set of Schr\"oder paths of length~$2n$,
and $\scrs = \bigcup_{n=0}^\infty \scrs_{2n}$.

There is an obvious bijection between Schr\"oder paths and Motzkin paths:
namely, every long level step is mapped onto a level step.
If we apply Flajolet's master theorem with
$a_{i-1} b_i = \alpha_i t$ and $c_i = \delta_{i+1} t$
to the resulting Motzkin path
(note that $t$ is now conjugate to the semi-length
 of the underlying Schr\"oder path),
we obtain
\be
   \sum_{n=0}^\infty t^n \sum_{\omega \in \scrs_{2n}}  W(\omega)
   \;=\;
   \cfrac{1}{1 - \delta_1 t - \cfrac{\alpha_1 t}{1 - \delta_2 t - \cfrac{\alpha_2 t}{1 - \cdots}}}
   \;\,,
 \label{eq.flajolet.schroder}
\ee
so that the generating function of Schr\"oder paths
with height-dependent weights
is given by the T-fraction \reff{eq.Tfrac.def}.
More precisely, every rise gets a weight~1,
every fall starting at height~$i$ gets a weight $\alpha_i$,
and every long level step at height~$i$ gets a weight $\delta_{i+1}$.
This combinatorial interpretation of T-fractions in terms of Schr\"oder paths
was found recently by several authors
\cite{Fusy_15,Oste_15,Josuat-Verges_18,Sokal_totalpos}.

Since the up steps $i \to i+1$ and the down steps $i+1 \to i$
in a Schr\"oder path can be paired,
we may alternatively distribute the weights on rises and falls
by assigning a weight~1 to all rises and falls starting at even heights,
a weight $\alpha_{2k}$ to a rise starting at odd height $2k-1$,
and a weight $\alpha_{2k-1}$ to a fall starting at odd height $2k-1$;
a~long level step at height~$i$ gets a weight $\delta_{i+1}$ as before.

\bigskip

{\bf Remark.}
With these preliminaries in place,
we can now give a combinatorial interpretation/proof
of Proposition~\ref{prop.transformation}:

\combinatorialproofof{Proposition~\ref{prop.transformation}}
We will provide a bijective interpretation
of the identity~\reff{eq.prop.transformation.equiv},
which is equivalent to Proposition~\ref{prop.transformation}.
We use our alternative assignment of weights on Schr\"oder paths.

By the {\em initial long level steps}\/ of a Schr\"oder path $\omega$,
we shall refer to the long level steps occurring before the first rise, if any.
We will add the adjective {\em non-initial}\/
to refer to all other long level steps.
Also, by a {\em restricted Schr\"oder path}\/,
we mean a  Schr\"oder path with no long level steps at even heights.
Given a Schr\"oder path $\omega$, we define a
restricted Schr\"oder path $\omega^\flat$
by simply removing all the long level steps at even heights in $\omega$.
The map $\omega\mapsto\omega^\flat$ is clearly a surjection.

Notice that a non-initial long level step at an even height
is preceded either by another non-initial long level step at the same height
or by a rise or fall starting at an odd height.
With this observation, we can immediately describe the preimage of a restricted Schr\"oder path $\omega^\flat$
and hence prove
equation~\reff{eq.prop.transformation.equiv}:
Firstly, we may choose to insert any number of initial long level steps
to $\omega^\flat$;
this contributes the prefactor $\dfrac{1}{1-\delta_1 t}$ on the left-hand side.
For a rise from height~$2k-1$ to height~$2k$
(which has weight $\alpha_{2k}$),
we may choose to insert any number of 
non-initial long level steps at height~$2k$ after this rise;
this justifies the substitution
$\alpha_{2k} \mapsto \dfrac{\alpha_{2k}}{1-\delta_{2k+1} t}$.
Likewise, for a fall from height~$2k-1$ to height~$2k-2$
(which has weight $\alpha_{2k-1}$),
we may choose to insert any number of 
non-initial long level steps at height~$2k-2$ after this fall;
this justifies the substitution
$\alpha_{2k-1} \mapsto \dfrac{\alpha_{2k-1}}{1-\delta_{2k-1} t}$.
\qed

\subsection{Labeled Dyck, Motzkin and Schr\"oder paths}
   \label{subsec.prelimproofs.2}

Let $\bfscra = (\scra_h)_{h \ge 0}$, $\bfscrb = (\scrb_h)_{h \ge 1}$
and $\bfscrc = (\scrc_h)_{h \ge 0}$ be sequences of finite sets.
An
\textbfit{$(\bfscra,\bfscrb,\bfscrc)$-labeled Motzkin path of length $\bm{n}$}
is a pair $(\omega,\xi)$
where $\omega = (\omega_0,\ldots,\omega_n)$
is a Motzkin path of length $n$,
and $\xi = (\xi_1,\ldots,\xi_n)$ is a sequence satisfying
\be
   \xi_i  \:\in\:
   \begin{cases}
       \scra(h_{i-1})  & \textrm{if step $i$ is a rise (i.e.\ $h_i = h_{i-1} + 1$)}
              \\[1mm]
       \scrb(h_{i-1})  & \textrm{if step $i$ is a fall (i.e.\ $h_i = h_{i-1} - 1$)}
              \\[1mm]
       \scrc(h_{i-1})  & \textrm{if step $i$ is a level step (i.e.\ $h_i = h_{i-1}$)}
   \end{cases}
 \label{eq.xi.ineq}
\ee
where $h_{i-1}$ (resp.~$h_i$) is the height of the Motzkin path
before (resp.~after) step $i$.
[For typographical clarity
 we have here written $\scra(h)$ as a synonym for $\scra_h$, etc.]
We call $\xi_i$ the \textbfit{label} associated to step $i$.
We call the pair $(\omega,\xi)$
an \textbfit{$(\bfscra,\bfscrb)$-labeled Dyck path}
if $\omega$ is a Dyck path (in this case $\bfscrc$ plays no role).
We denote by $\scrm_n(\bfscra,\bfscrb,\bfscrc)$
the set of $(\bfscra,\bfscrb,\bfscrc)$-labeled Motzkin paths of length $n$,
and by $\scrd_{2n}(\bfscra,\bfscrb)$
the set of $(\bfscra,\bfscrb)$-labeled Dyck paths of length $2n$.

We define a \textbfit{$(\bfscra,\bfscrb,\bfscrc)$-labeled 
Schr\"oder path}
in an analogous way;
now the sets $\scrc_h$ refer to long level steps.
We denote by $\scrs_{2n}(\bfscra,\bfscrb,\bfscrc)$
the set of $(\bfscra,\bfscrb,\bfscrc)$-labeled Schr\"oder paths
of length $2n$.

Let us stress that the sets $\scra_h$, $\scrb_h$ and $\scrc_h$ are allowed
to be empty.
Whenever this happens, the path $\omega$ is forbidden to take a step
of the specified kind starting at the specified height.


\bigskip

{\bf Remark.}  What we have called an
$(\bfscra,\bfscrb,\bfscrc)$-labeled Motzkin path
is (up to small changes in notation)
called a {\em path diagramme}\/ by Flajolet \cite[p.~136]{Flajolet_80}
and a {\em history}\/ by Viennot \cite[p.~II-9]{Viennot_83}.
Often the label sets $\scra_h, \scrb_h, \scrc_h$ are intervals of integers,
e.g.\ $\scra_h = \{1,\ldots,A_h\}$ or $\{0,\ldots,A_h\}$;
in this case the triplet $({\bf A},{\bf B},{\bf C})$
of sequences of maximum values is called a {\em possibility function}\/.
On the other hand, it is sometimes useful to employ labels
that are {\em pairs}\/ of integers
(e.g.\ \cite[Section~6.2]{Sokal-Zeng_masterpoly}
 and \cite[Section~7]{Deb-Sokal_Genocchi}).
It therefore seems preferable to state the general theory
without any specific assumption about the nature of the label sets.
\myendremark

\bigskip

Following Flajolet \cite[Proposition~7A]{Flajolet_80},
we can state a ``master J-fraction'' for
$(\bfscra,\bfscrb,\bfscrc)$-labeled Motzkin paths.
Let ${\bf a} = (a_{h,\xi})_{h \ge 0 ,\; \xi \in \scra(h)}$,
${\bf b} = (b_{h,\xi})_{h \ge 1 ,\; \xi \in \scrb(h)}$
and ${\bf c} = (c_{h,\xi})_{h \ge 0 ,\; \xi \in \scrc(h)}$
be indeterminates;
we give an $(\bfscra,\bfscrb,\bfscrc)$-labeled Motzkin path $(\omega,\xi)$
a weight $W(\omega,\xi)$
that is the product of the weights for the individual steps,
where a rise starting at height~$h$ with label $\xi$ gets weight~$a_{h,\xi}$,
a~fall starting at height~$h$ with label $\xi$ gets weight~$b_{h,\xi}$,
and a level step at height~$h$ with label $\xi$ gets weight~$c_{h,\xi}$.
Then:

\begin{theorem}[Flajolet's master theorem for labeled Motzkin paths]
   \label{thm.flajolet_master_labeled_Motzkin}
We have
\be
   \sum_{n=0}^\infty t^n
   \!\!
   \sum_{(\omega,\xi) \in \scrm_n(\bfscra,\bfscrb,\bfscrc)} \!\!\!  W(\omega,\xi)
   \;=\;
   \cfrac{1}{1 - c_0 t - \cfrac{a_0 b_1 t^2}{1 - c_1 t - \cfrac{a_1 b_2 t^2}{1- c_2 t - \cfrac{a_2 b_3 t^2}{1- \cdots}}}}
\ee
as an identity in $\Z[{\bf a},{\bf b},{\bf c}][[t]]$, where
\be
   a_h  \;=\;  \sum_{\xi \in \scra(h)} a_{h,\xi}
   \;,\qquad
   b_h  \;=\;  \sum_{\xi \in \scrb(h)} b_{h,\xi}
   \;,\qquad
   c_h  \;=\;  \sum_{\xi \in \scrc(h)} c_{h,\xi}
   \;.
 \label{def.weights.akbkck}
\ee
\end{theorem}

\noindent
This is an immediate consequence of Theorem~\ref{thm.flajolet}
together with the definitions.

By specializing to ${\bf c} = \bzero$ and replacing $t^2$ by $t$,
we obtain the corresponding theorem
for $(\bfscra,\bfscrb)$-labeled Dyck paths:

\begin{corollary}[Flajolet's master theorem for labeled Dyck paths]
   \label{cor.flajolet_master_labeled_Dyck}
We have
\be
   \sum_{n=0}^\infty t^n
   \!\!
   \sum_{(\omega,\xi) \in \scrd_{2n}(\bfscra,\bfscrb)} \!\!\!  W(\omega,\xi)
   \;=\;
   \cfrac{1}{1 - \cfrac{a_0 b_1 t}{1 - \cfrac{a_1 b_2 t}{1- \cfrac{a_2 b_3 t}{1- \cdots}}}}
\ee
as an identity in $\Z[{\bf a},{\bf b}][[t]]$, where
$a_h$ and $b_h$ are defined by \reff{def.weights.akbkck}.
\end{corollary}

Similarly, for labeled Schr\"oder paths we have:

\begin{theorem}[Flajolet's master theorem for labeled Schr\"oder paths]
   \label{thm.flajolet_master_labeled_Schroder}
We have
\be
   \sum_{n=0}^\infty t^n
   \!\!
   \sum_{(\omega,\xi) \in \scrs_{2n}(\bfscra,\bfscrb,\bfscrc)} \!\!\!  W(\omega,\xi)
   \;=\;
   \cfrac{1}{1 - c_0 t - \cfrac{a_0 b_1 t}{1 - c_1 t - \cfrac{a_1 b_2 t}{1- c_2 t - \cfrac{a_2 b_3 t}{1- \cdots}}}}
\ee
as an identity in $\Z[{\bf a},{\bf b},{\bf c}][[t]]$, where
$a_h, b_h, c_h$ are defined by \reff{def.weights.akbkck},
with $c_{h,\xi}$ now referring to long level steps.
\end{theorem}

\section{Bijective proofs of Theorems~\ref{thm.abt.master.J} and \ref{thm.sbt.master}}
   \label{sec.bijections.proofs}

Continued fractions for increasing binary trees go back to
the celebrated bijection of 
Fran\c{c}on and Viennot \cite{Francon_79,Francon_78}:
though ordinarily understood as a bijection
from permutations to labeled Motzkin paths,
the Fran\c{c}on--Viennot bijection can also be understood,
by virtue of the standard bijection from
increasing binary trees to permutations \cite[pp.~44--45]{Stanley_12},
as a bijection from increasing binary trees to labeled Motzkin paths;
in this form it was first written down by Flajolet \cite{Flajolet_80}
(see also \cite[Ch.~4b]{Viennot_abjc1}).
This bijection was rediscovered by
Albert, Linton and Ru\v{s}kuc~\cite{Albert_05},
and in the permutation-patterns community
it is often referred to as {\em insertion encoding}\/.

This bijection has recently been generalized
by Kuba and Varvak \cite{Kuba_21} 
and P\'etr\'eolle, Sokal and Zhu \cite{latpath_SRTR}
to the setting of increasing trees with higher arity,
and equivalently to generalized Stirling permutations.
We will provide here yet another generalization.



\subsection{Bijection for increasing restricted ternary trees:
   Proof of Theorem~\ref{thm.abt.master.J}}
\label{subsec.bijection.J.abt}


In this section we will construct a bijection from
increasing restricted ternary trees on the vertex set $[n+1]$
to labeled Motzkin paths of length $n$, as follows: Given a 
tree $T\in\RT_{n+1}$, we first define the path $\omega$ and then define the labels $\xi_i$, which will lie in the sets
\begin{subeqnarray}
   \scra_h  & = &  \{0,\ldots, h\}  \\
   \scrb_h  & = &  \{0,\ldots, h\}  \quad\hbox{for $h\geq 1$}  \\
   \scrc_{h}        & = &  \{1,2,3\}\times\{0,\ldots, h \} 
 \label{def.abc.bt}
\end{subeqnarray}
A level step that has label $\xi_h \in \{i\}\times \{0,\ldots, h \}$
will be called a level step of type $i$ ($i = 1, 2, 3$).


\bigskip

\noindent{\bf Step 1: Definition of the Motzkin path.}

Given a restricted ternary tree $T\in \RT_{n+1}$, 
we classify the indices $i\in[n]$ according to their node type.
We then define a Motzkin path $\omega = (\omega_0,\ldots,\omega_n)$
starting at $\omega_0 = (0,0)$ and ending at $\omega_n = (n,0)$,
with steps $s_1,\ldots,s_n$ as follows:
\begin{itemize}
\item If $N(i,T) = 101$, $s_i$ is a rise.
\item If $N(i,T) = 000$, $s_i$ is a fall.
\item If $N(i,T) = 100$, $010$ or $001$,
   $s_i$ is a level step of type 1, 2 or 3, respectively.
\end{itemize}
These definitions can equivalently be written as
\be
   h_i - h_{i-1}  \;=\;  \deg(i) \,-\, 1   \;.
 \label{eq.step.deg}
\ee
Please note also that no step is assigned to vertex $n+1$,
which is anyway always a leaf.

It is clear that $\omega$ consists of $n$ steps
which are rises, falls and level steps. 
Thus, it remains to show that $\omega$ always stays on or above the $x$-axis
and that it ends at $(n,0)$. 
We will do this by obtaining a precise interpretation of the heights:

\begin{lemma} [Interpretation of the heights]
\label{lem.abt.height.motzkin}
The height of $\omega$ at position $i$ is given by
\be
    h_i \;=\; \lev(i+1,T)\;.
 \label{eq.lem.abt.height.motzkin}
\ee
In particular, $h_i \ge 0$ and $h_n = \lev(n+1,T) = 0$.
\end{lemma}
\proof 
We proceed by induction.
By definition, the path $\omega$ starts at height $h_0 = 0$.
On the other hand, the root is here vertex 1;
and was observed following \reff{eq.def.level.vertex.sbt},
we have $\lev({\rm root}, T) = 0$.
This proves the base case $i=0$ of \reff{eq.lem.abt.height.motzkin}.

Consider now $i \ge 1$, and assume that $h_{i-1} = \lev(i,T)$.
We will compare ${\lev(i+1,T)}$ with $\lev(i,T)$,
and we will show that
\be
   \lev(i+1,T) - \lev(i,T)  \;=\;  \deg(i)-1\;.
 \label{eq.diff.lev.deg.1.bt}
\ee
By \reff{eq.step.deg}, this will complete the proof of the inductive step.

Let us start from the definitions
\begin{subeqnarray}
   \lev(i,T)    & = &  \#\{j \in [n+1] \colon\, p(j) < i < j\}
      \\[2mm]
   \lev(i+1,T)  & = &  \#\{j \in [n+1] \colon\, p(j) < i+1 < j\}
      \slabel{eq.lev.i+1}
\end{subeqnarray}
We see that:
\begin{itemize}
   \item  A vertex $j$ contributes to both $\lev(i,T)$ and $\lev(i+1,T)$
       in case $p(j) < i < i+1 < j$.
   \item   A vertex $j$ contributes to $\lev(i,T)$ but not to $\lev(i+1,T)$
       in case $p(j) < i$ and $j = i+1$
       --- or in other words, $j = i+1$ is not a child of $i$.
       [Note that we must always have $p(i+1) \le i$.]
   \item   A vertex $j$ contributes to $\lev(i+1,T)$ but not to $\lev(i,T)$
       in case $p(j) = i$ and $j > i+1$
       --- or in other words, $j$ is a child of $i$ other than $i+1$.
\end{itemize}
It follows that
\begin{subeqnarray}
   & &
   \lev(i+1,T)
   \,-\,
   \lev(i,T)
         \nonumber \\[2mm]
   & &
   \qquad \;=\;
       \bigl( \deg(i) \,-\, {\rm I}[i+1 \hbox{ is a child of } i] \bigr)
       \:-\:
       {\rm I}[i+1 \hbox{ is not a child of } i]
       \qquad\qquad
              \\[2mm]
   & &
   \qquad \;=\;
       \deg(i) - 1  \;.
\end{subeqnarray}
(Here ${\rm I}[\hbox{\sl proposition}] = 1$ if {\sl proposition} is true,
and 0 if it is false.)
\qed

{\bf Remark.}
For a general increasing tree (not necessarily restricted ternary),
the steps \reff{eq.step.deg} define in general an upper-\L{}ukasiewicz path,
i.e.~$h_i - h_{i-1} \in \{-1,0,1,2,\ldots\}$.
The key identity \reff{eq.diff.lev.deg.1.bt}
continues to hold in this generality:
see \cite[Lemma~3.3]{latpath_lah} and \cite{masterpoly_trees}.
\myendremark

\bigskip

\noindent{\bf Step 2: Definition of the labels $\xi_i$.} 
Fix a consistent tree-traversal algorithm $\bfA$.
We will now describe the labels.
We assign labels to the steps according to 
the status of the corresponding vertices as follows:
\begin{equation}
\xi_{i} 
\;\eqdef\;
\#\{ j \colon\: p(j) < i < j \text{ and } j<_\bfA i \}
\label{eq.labels.abt.bt}
\end{equation}
where $<_\bfA$ is the total order on vertices given by
the tree-traversal algorithm $\bfA$.
In other words,
\be
   \xi_{i} \;=\; \nid(i,T)
 \label{eq.labels.abt.nid.bt}
\ee
as defined in \reff{eq.def.croix.nid.vertex}.
To verify that the inequalities~\eqref{def.abc.bt} are satisfied,
we need to check that 
$0\leq \xi_{i} \leq \,  h_{i-1}$,
where $h_{i-1}$ is the starting height of step $s_i$.
The lower bound is immediate from the definition of $\xi_{i}$.
For the upper bound, notice that from Lemma~\ref{lem.abt.height.motzkin} it follows
$h_{i-1}=\lev(i,T)$.
Thus, we need to show that $\xi_{i} \leq \lev(i,T)$,
which is clear from the definitions
\eqref{eq.labels.abt.bt}/\eqref{eq.def.level.vertex.sbt}.
In fact, we also have an interpretation of the difference in terms of the 
statistic croix:
\be
h_{i-1}- \xi_i 
\;=\; 
\lev(i,T) - \xi_i 
\;=\; 
\croix(i,T)
\label{eq.height.minus.nid.abt}
\ee
by \eqref{eq.croix.plus.nid}.

\bigskip

\noindent {\bf Step 3: Proof of bijection.}
We prove that the mapping $T\mapsto (\omega,\xi)$ is a bijection,
by constructing the inverse bijection.
To do this, we first use the path $\omega$
to identify the node types of the vertices.
We will then use the labels $\xi$ to glue the vertices together 
and construct our restricted ternary tree $T\in \RT_n$;
the details are as follows.

We first define a class of intermediate objects in our bijection:
a {\em slotted restricted ternary tree}\/
is an increasing restricted ternary tree
whose set of vertex labels is $[i]\cup \{\infty\}$
but we now allow the label $\infty$ to be assigned to multiple vertices
(all of which must be leaves).
(We think of the vertices labeled $\infty$ as placeholders
where new vertices may be inserted into the tree in the future.)

Given a tree $T\in \RT_{n+1}$,
we define the tree $T|_{i}$
to be the subtree of $T$ 
consisting of vertices $\{1,2,\ldots,i\}$ along with their children,
in which the vertices with labels $>i$ are relabeled to $\infty$.
Clearly, $T|_i$ is a slotted restricted ternary tree,
and $T|_{n+1} = T$.
By the {\em history}\/ of tree $T$,
we will mean the sequence of slotted restricted ternary trees
$T|_1 \to T|_2 \to \ldots \to T|_{n+1}$.
Notice that for $i<n+1$,
the number of vertices labeled $\infty$ in the tree $T|_i$ is
the number of vertices $\ge i+1$ in $T$ whose parents are $\le i$:
by \reff{eq.lev.i+1} this is $\lev(i+1,T)+1$. 

We are now ready to build the tree $T$ 
from the pair $(\omega,\xi)$ 
by successively reading the steps $s_i$ and the labels $\xi_{i}$,
using which we construct $T|_i$ from $T|_{i-1}$.
We select the $(\xi_{i}+1)$-th vertex labeled $\infty$
in the traversal algorithm $\bfA$ applied to $T|_{i-1}$,
and rename this vertex to $i$. 
And then we choose to add children (labeled $\infty$) to vertex $i$,
depending on the status of the step $s_i$, as follows:
\begin{itemize}
\item If step $s_i$ is a rise, 
we add a left child and a right child to $i$,
and both children get the label $\infty$.

\item If step $s_i$ is a fall,
we do not add any children.

\item If step $s_i$ is a level step of type 1,
we add a left child with label $\infty$.

\item If step $s_i$ is a level step of type 2,
we add a middle child with label $\infty$.

\item If step $s_i$ is a level step of type 3,
we add a right child with label $\infty$.

\end{itemize}

Finally, at the end of this process, the tree $T|_{n}$
will only contain a single vertex labeled $\infty$;
we rename this vertex to $n+1$,
thereby constructing $T = T|_{n+1}$.

\bigskip

\noindent{\bf Step 4: Computation of the weights.}
We can now compute the weights associated to 
the Motzkin path $\omega$ in Theorem~\ref{thm.flajolet_master_labeled_Motzkin},
which we recall are $a_{h,\xi}$ for a rise starting at 
height $h$ with label $\xi$, 
$b_{h,\xi}$ for a fall starting at 
height $h$ with label $\xi$, 
and $c_{h,\xi}$ for a level step starting at 
height $h$ with label $\xi$.
We do this by putting together the information collected in 
the description of the bijection.
We write out the contribution of step $i$ of
the path $\omega$ 
as per the different node types:
\begin{itemize}
    \item[(a)] Rise from height $h$ to height $h+1$:
        \begin{itemize}
            \item This step corresponds to a vertex $i$ 
            having node type $N(i,T) = 101$,
            so from \reff{eq.def.Qn.abt} it contributes the letter $\sfa$.
            \item From~\eqref{eq.labels.abt.bt}/\eqref{eq.labels.abt.nid.bt}
               we know that $\xi_{i} = \nid(i,T)$,
            and from~\eqref{eq.height.minus.nid.abt}
            we get $h-\xi_{i} = \croix(i,T)$.
        \end{itemize}
        From~\eqref{eq.def.Qn.abt},
        we get that the weight for this step is 
        \be
        a_{h,\xi} \;=\; 
        \sfa_{h-\xi, \xi}
        \ee

    \item[(b)] Fall from height $h$ to height $h-1$:
        \begin{itemize}
            \item This step corresponds to a vertex $i$ 
            having node type $N(i,T) = 000$,
            so from \reff{eq.def.Qn.abt} it contributes the letter $\sfb$.
            \item Once again we have $\xi_{i} = \nid(i,T)$
               and $h-\xi_{i} = \croix(i,T)$.
        \end{itemize}
        From~\eqref{eq.def.Qn.abt},
        we get that the weight for this step is 
        \be
        b_{h,\xi} \;=\; 
        \sfb_{h-\xi, \xi}
        \ee
    \item[(c)] Level step at height $h$:
        \begin{itemize}
            \item This step corresponds to a vertex $i$ 
            having node type $N(i,T) = 100$, $010$ or $001$,
            so from \reff{eq.def.Qn.abt} it contributes
            the letter $\sfc$, $\sff$ or $\sfd$, respectively.
            
            \item Once again we have $\xi_{i} = \nid(i,T)$
               and $h-\xi_{i} = \croix(i,T)$.
        \end{itemize}
        From~\eqref{eq.def.Qn.abt},
        we get that the weight for this step is 
        \be
        c_{h,\xi} \;=\; 
        \sfc_{h-\xi, \xi}
        \,+\,
        \sfd_{h-\xi, \xi}
        \,+\,
        \sff_{h-\xi, \xi}
        \ee
\end{itemize}
Putting this all together in Theorem~\ref{thm.flajolet_master_labeled_Motzkin},
we obtain a J-fraction with
\begin{eqnarray}
   \beta_h &=& 
   (\hbox{rise from $h-1$ to $h$})
    \,\times\,
    (\hbox{fall from $h$ to $h-1$})
    \nonumber\\
 &=& 
 \biggl(\sum_{\xi=0}^{h-1} \sfa_{h-1-\xi, \xi} \biggr)
 \biggl(\sum_{\xi=0}^{h} \sfb_{h-\xi, \xi} \biggr)
 \\[2mm]
 \gamma_{h} &=& 
 \hbox{ level step at height $h$}
 \nonumber\\
 &=& 
 \sum_{\xi=0}^{h} \sfc_{h-\xi,\xi}
 \,+\,
 \sum_{\xi=0}^{h} \sfd_{h-\xi,\xi}
 \,+\,
 \sum_{\xi=0}^{h} \sff_{h-\xi,\xi}
 \label{eq.computation.master.abt.J}
\end{eqnarray}

This completes the proof of Theorem~\ref{thm.abt.master.J}.
\qed

We can now deduce Theorem~\ref{thm.Simple.J.abt} 
as a corollary:

\proofof{Theorem~\ref{thm.Simple.J.abt}}
It suffices to make the substitutions
\be
\sfa_{\ell, \ell'}  \;=\;  x_1 \;,\quad
\sfb_{\ell,\ell'}  \;=\;  y_1 \;,\quad
\sfc_{\ell,\ell'}  \;=\;  x_2 \;,\quad
\sfd_{\ell, \ell'}  \;=\;  y_2 \;,\quad
\sff_{\ell, \ell'}  \;=\;  w
\label{eq.substitutions.simple.rtt}
\ee
in Theorem~\ref{thm.abt.master.J}.
\qed

%
%

\bigskip

\subsection{Bijection for increasing interval-labeled restricted ternary trees:
   Proof of Theorem~\ref{thm.sbt.master}}
\label{subsec.sbt.bijection}

In this section we will construct a bijection from
increasing interval-labeled restricted ternary trees on the label set $[0,n]$
to labeled Schr\"oder paths of length $2n$, as follows:
Given a tree $T\in\IRT_n$, we first define the path $\omega$
and then define the labels $\xi_i$, which will lie in the sets

%
\begin{subeqnarray}
   \scra_h        & = &  \{0\}  \qquad\qquad\qquad\;\:\:\hbox{for $h$ even}  \\
   \scra_h        & = &  \{0,\ldots, \lfloor h/2 \rfloor \}    \qquad\hbox{for $h$ odd}  \\
   \scrb_h        & = &  \{0\}  \qquad\qquad\quad\quad\;\:\:\hbox{for $h$ even}  \\
   \scrb_h        & = &  \{0,\ldots, \lfloor h/2 \rfloor \}    \qquad\hbox{for $h$ odd}  \\
   \scrc_{h}       & = &  \{0\} 
   \quad\qquad\qquad\quad\;\:\:\hbox{for $h$ even}\\
   \scrc_{h}        & = &  \{0,\ldots, \lfloor h/2 \rfloor \}        
   \qquad\hbox{for $h$ odd}
 \label{def.abc.sbt}
\end{subeqnarray}
Notice that only steps starting at an odd height 
may have a non-unique choice of label.
We will then interpret the heights and labels,
which will show that our path and labels are well-defined.
Finally, we will prove that the map $T \mapsto (\omega, \xi)$
is indeed a bijection, by describing the inverse bijection.

\bigskip

{\bf Step 1: Definition of the Schr\"oder path.}
Recall that in an increasing interval-labeled restricted ternary tree
$T\in\IRT_n$,
the vertex labels are disjoint intervals in $[0,n]$.
The vertices therefore have a natural total order,
obtained by comparing their label sets.
Let the vertices of $T$ be ${v_0 < v_1 < \ldots < v_m}$ in this total order;
note that $v_0$ is the root, $v_m$ is a leaf,
$0\in L_{v_0}$ and $n\in L_{v_m}$.
Let us call the tree \textbfit{trivial}
if it consists only of a root
(then $m=0$ and the root has label set $[0,n]$),
and \textbfit{nontrivial} otherwise
(then $m>0$ and the root has a label set $[0,j]$ with $j < n$).

We will now describe a Schr\"oder path $\omega$ of length $2n$;
to do this we will assign a segment $\omega(v_i)$ 
to every vertex $v_i$ in the tree,
and the path $\omega$ will be obtained by concatenating the segments:
$\omega = \omega(v_0)\omega(v_1)\cdots\omega(v_m)$.

Let $v$ be a vertex of $T$ with $L_v = \{l,l+1, \ldots, l+j\}$,
so that the label surplus of this vertex is $j = |L_v| - 1$.
We now describe the segment $\omega(v)$ 
as a word in the letters $\{\nearrow, \searrow, \longrightarrow \}$,
which represent rise, fall and long level step, respectively:
\begin{eqnarray}
   & &
   \hspace*{-1.3cm}
   \bullet\;
   \textrm{If $l=0$ and $j=n$ then $\omega(v) = (\longrightarrow)^n$.}
       \label{eq.RT.step.I} \\
   & &
   \hphantom{\bullet}\;
          \textrm{(Here $v$ is the root of a trivial tree.)}
       \nonumber \\[3mm]
   & &
   \hspace*{-1.3cm}
   \bullet\;
   \textrm{If $l=0$ and $j<n$ then $\omega(v) = (\longrightarrow)^j \nearrow$.}
       \label{eq.RT.step.II} \\
   & &
   \hphantom{\bullet}\;
          \textrm{(Here $v$ is the root of a nontrivial tree.)}
       \nonumber \\[3mm]
   & &
   \hspace*{-1.3cm}
   \bullet\;
   \textrm{If $l>0$ and $l+j =n$, then
          $\omega(v) = \searrow (\longrightarrow)^j$.}
   \hspace*{5cm}
       \label{eq.RT.step.III} \\
   & &
   \hphantom{\bullet}\;
          \textrm{(Here $v$ is the last vertex $v_m$ of a nontrivial tree.)}
       \nonumber
\end{eqnarray}
\begin{itemize}
    \item If $l>0$ and $l+j<n$ then 
    \begin{equation}
        \omega(v) 
        \;=\;
        \begin{cases}
            \nearrow (\longrightarrow)^j \nearrow \qquad \text{if $N(v,T) = 101$}\\[2mm]
            \searrow (\longrightarrow)^j \searrow \qquad \text{if $N(v,T) = 000$}\\[2mm]
            \nearrow (\longrightarrow)^j \searrow \qquad \text{if $N(v,T) = 100$}\\[2mm]
            \searrow (\longrightarrow)^j \nearrow \qquad \text{if $N(v,T) = 001$}\\[2mm]
            \longrightarrow \qquad\qquad\quad\;\;\, \text{if $N(v,T) = 010$}
        \end{cases}
  \label{eq.RT.step}
    \end{equation}
           (Here $v$ is neither $v_0$ nor $v_m$, and of course the tree
           is nontrivial.)
\end{itemize}
Notice that:

1) The cases $l = 0$ correspond to the possibilities
when $v$ is the root vertex $v_0$,
while the cases $l > 0$ correspond to the possibilities
when $v$ is a non-root vertex $v_i$ with $i \ge 1$.

2) Similarly, the cases $l+j=n$ correspond to the possibilities
when $v$ is the final vertex $v_m$,
while the cases $l+j<n$ correspond to the possibilities
when $v$ is a non-final vertex $v_i$ with $i \le m-1$.

3) When $N(v,T) = 010$, we must have $j=0$,
by virtue of our condition that vertices with a middle child
are always single-labeled.

4) In \reff{eq.RT.step} with $N(v,T) \neq 010$, we see that
the first step is a rise if $v$ has a left child, and a fall if not;
likewise, the last step is a rise if $v$ has a right child, and a fall if not.

5) The length of $\omega(v)$ is
\be
   |\omega(v)|
   \;=\;
   \begin{cases}
      2j\hphantom{+1}\;\,  \;=\; 2 |L_v| -2   & \textrm{if $l=0$ and $j=n$}   \\
      2j+1 \;=\; 2 |L_v| -1   & \textrm{if $l=0$ and $j<n$}   \\
      2j+1 \;=\; 2 |L_v| -1   & \textrm{if $l>0$ and $l+j =n$}   \\
      2j+2 \;=\; 2 |L_v|      & \textrm{if $l>0$ and $l+j<n$}
   \end{cases}
\ee 
Since $l=0$ must occur exactly once,
and $l+j=n$ must also occur exactly once, we have
\be
   \sum\limits_{i=0}^n |\omega(v_i)|
   \;=\;
   \sum\limits_{i=0}^n 2 |L_{v_i}| \:-\: 2
   \;=\;
   2n  \;.
\ee
So the path $\omega$ is indeed of length $2n$.

\medskip

It is clear that $\omega$ consists of rises, falls and long level steps
and that it starts at $(0,0)$
and ends at $(2n,k)$ for some $k\in \mathbb{Z}$.
To show that $\omega$ is a Schr\"oder path,
we need to show that it always stays on or above the $x$-axis
and that it ends at $(2n,0)$.
We will do this by obtaining a precise interpretation of the heights. 
In particular, we interpret the starting and ending heights of the segments
$\omega(v)$, as follows:

\begin{lemma}[Interpretation of the heights]
   \label{lemma.RT.heights}
Let $v$ be a vertex of $T$ with label set~$L_v$.
If $L_v = [0,n]$ (which corresponds to $v$ being the root of a trivial tree),
then the path consists of $n$ long level steps at height 0.
Otherwise the tree is nontrivial, and:
\begin{enumerate}
    \item The segment $\omega(v_0)$ starts at height $0$ and ends at height $1$.
    \item For $v=v_i$ with $1 \le i \le m$,
    the segment $\omega(v)$ starts at height $2\, \lev(v,T)+1$;
    when~$1 \le i < m$, it ends at height $h$ where $h$ is given by
        \begin{equation}
            h
            \;=\;
            \begin{cases}
                2\, \lev(v,T) +3 \qquad\text{if $N_v(T)=101$} \\
                2\, \lev(v,T) +1 \qquad \text{if $N_v(T) =100$, $001$ or $010$}\\
                2\, \lev(v,T) -1 \qquad \text{if $N_v(T)=000$}
            \end{cases}
  \label{eq.ending_height}
        \end{equation}
    \item The segment $\omega(v_m)$ starts at height $1$ and ends at height $0$.
\end{enumerate}
\label{lem.sbt.height}
In particular,
all the segments $\omega(v_i)$ for $i \neq 0$ start at an odd height,
and all the segments $\omega(v_i)$ for $i \neq m$ end at an odd height.
\end{lemma}

\proof
(a) is clear from the definition.

(b) Let us now look at the segments $\omega(v_i)$ for $1 \le i \le m$.
It suffices to prove the statement about the starting height,
because the statement \reff{eq.ending_height} about the ending height
when $i < m$ then follows immediately from this,
using the definition \reff{eq.RT.step} of the step
(note that for $1 \le i \le m-1$ we are always in the case
 $l > 0$ and $l+j < n$).
We proceed by induction on $i$.

From (a) it follows that the segment $\omega(v_1)$ starts at height $1$.
On the other hand, $\lev(v_1, T) = 0$ because the root has only one child
(this was observed already following \reff{eq.def.level.vertex.sbt}).
This proves the base case $i=1$ of the induction.

Consider now $i \in [2,m]$,
and assume the inductive hypothesis
that the segment $\omega(v_{i-1})$ starts at height $2\lev(v_{i-1},T) + 1$.
We now compare $\lev(v_i,T)$ with $\lev(v_{i-1},T)$,
and claim that the following equality holds:
\be
\lev(v_i,T) - \lev(v_{i-1},T) \;=\; \deg(v_{i-1})-1\;.
\label{eq.diff.lev.deg}
\ee
The proof of \reff{eq.diff.lev.deg} is identical to the proof of
\reff{eq.diff.lev.deg.1.bt} in the preceding subsection:
indeed, the proof holds without alteration
for any increasing tree on a totally ordered vertex set.
On the other hand,
\be
    \deg(v_{i-1})-1
    \;=\;
    \begin{cases}
       +1   & \textrm{when $N_{v_{i-1}}(T) = 101$}   \\
       0    & \textrm{when $N_{v_{i-1}}(T) = 100$, $001$ or $010$}   \\
       -1   & \textrm{when $N_{v_{i-1}}(T) = 000$}
    \end{cases}
\ee
It then follows from \reff{eq.ending_height}
that the ending height of the segment $\omega(v_{i-1})$,
which is also the starting height of the segment $\omega(v_i)$,
is $2\lev(v_i,T) + 1$.

(c) Since $\lev(v_m,T) = 0$, it follows from~(b)
that $\omega(v_m)$ starts at height 1.
From the case $l > 0$ and $l+j=n$ of the definition,
it therefore ends at height 0.
\qed

\begin{corollary}
$\omega$ is a Schr\"oder path.
\end{corollary}

\proof
All that remains to prove is that the path stays always at height $\ge 0$,
even in-between the starting and ending points of steps.
The only potential danger is a step
$\searrow (\longrightarrow)^j \nearrow$ starting at height 0.
But this cannot happen, because Lemma~\ref{lemma.RT.heights} guarantees that
for $1 \le i \le m-1$ the step $\omega(v_i)$ starts at a height $\ge 1$.
\qed

%
%
%
%
%
%
\def\blue{\textcolor{blue}}
\def\red{\textcolor{red}}
\def\green{\textcolor{green}}
\def\yellow{\textcolor{yellow}}

\begin{figure}[t]
\vspace*{1.5cm}
\hspace*{2cm}
\centering
\Scale[0.5]{
\begin{picture}(50,120)(450,-40)
\setlength{\unitlength}{1.2cm}
\linethickness{.3mm}
\put(0,0){\vector(1,0){24.5}}\put(0,0){\vector(0,1){4.5}}
\put(-0.05,-0.5){$0$}
\put(0.7,-0.5){$1$}\put(0.75,0){\circle*{.1}}
\put(1.45,-0.5){$2$}\put(1.5,0){\circle*{.1}}
\put(2.2,-0.5){$3$}\put(2.25,0){\circle*{.1}}
\put(2.95,-0.5){$4$}\put(3,0){\circle*{.1}}
\put(3.7,-0.5){$5$}\put(3.75,0){\circle*{.1}}
\put(4.45,-0.5){$6$}\put(4.5,0){\circle*{.1}}
\put(5.2,-0.5){$7$}\put(5.25,0){\circle*{.1}}
\put(5.95,-0.5){$8$}\put(6,0){\circle*{.1}}
\put(6.7,-0.5){$9$}\put(6.75,0){\circle*{.1}}
\put(7.35,-0.5){$10$}\put(7.5,0){\circle*{.1}}
\put(8.1,-0.5){$11$}\put(8.25,0){\circle*{.1}}
\put(8.85,-0.5){$12$}\put(9,0){\circle*{.1}}
\put(9.6,-0.5){$13$}\put(9.75,0){\circle*{.1}}
\put(10.35,-0.5){$14$}\put(10.5,0){\circle*{.1}}
\put(11.10,-0.5){$15$}\put(11.25,0){\circle*{.1}}
\put(11.85,-0.5){$16$}\put(12,0){\circle*{.1}}
\put(12.6,-0.5){$17$}\put(12.75,0){\circle*{.1}}
\put(13.35,-0.5){$18$}\put(13.5,0){\circle*{.1}}
\put(14.10,-0.5){$19$}\put(14.25,0){\circle*{.1}}
\put(14.85,-0.5){$20$}\put(15,0){\circle*{.1}}
\put(15.60,-0.5){$21$}\put(15.75,0){\circle*{.1}}
\put(16.35,-0.5){$22$}\put(16.5,0){\circle*{.1}}
\put(17.10,-0.5){$23$}\put(17.25,0){\circle*{.1}}
\put(17.85,-0.5){$24$}\put(18,0){\circle*{.1}}
\put(18.60,-0.5){$25$}\put(18.75,0){\circle*{.1}}
\put(19.35,-0.5){$26$}\put(19.5,0){\circle*{.1}}
\put(20.10,-0.5){$27$}\put(20.25,0){\circle*{.1}}
\put(20.85,-0.5){$28$}\put(21,0){\circle*{.1}}
\put(21.60,-0.5){$29$}\put(21.75,0){\circle*{.1}}
\put(22.35,-0.5){$30$}\put(22.5,0){\circle*{.1}}
\put(23.1,-0.5){$31$}\put(23.25,0){\circle*{.1}}
\put(23.85,-0.5){$32$}\put(24,0){\circle*{.1}}
\put(-0.4,-0.1){$0$}
\put(-0.4,0.65){$1$}\put(0,0.75){\circle*{.15}}
\put(-0.4,1.4){$2$}\put(0,1.5){\circle*{.15}}
\put(-0.4,2.15){$3$}\put(0,2.25){\circle*{.15}}
\put(-0.4,2.9){$4$}\put(0,3){\circle*{.15}}
\put(-0.4,3.65){$5$}\put(0,3.75){\circle*{.15}}
\put(0,0){\red{\line(1,0){1.5}}}
\put(1.5,0){\red{\line(1,1){0.75}}}
\put(2.25,0.75){\red{\line(1,1){0.75}}}
\put(3,1.5){\red{\line(1,0){1.5}}}
\put(4.5,1.5){\red{\line(1,1){0.75}}}
\put(5.25,2.25){\red{\line(1,1){0.75}}}
\put(6,3){\red{\line(1,-1){0.75}}}
\put(6.75,2.25){\red{\line(1,1){0.75}}}
\put(7.5,3){\red{\line(1,1){0.75}}}
\put(8.25,3.75){\red{\line(1,0){1.5}}}
\put(9.75,3.75){\red{\line(1,0){1.5}}}
\put(11.25,3.75){\red{\line(1,-1){0.75}}}
\put(12,3){\red{\line(1,1){0.75}}}
\put(12.75,3.75){\red{\line(1,0){1.5}}}
\put(14.25,3.75){\red{\line(1,-1){0.75}}}
\put(15,3){\red{\line(1,-1){0.75}}}
\put(15.75,2.25){\red{\line(1,0){1.5}}}
\put(17.25,2.25){\red{\line(1,-1){0.75}}}
\put(18,1.5){\red{\line(1,0){1.5}}}
\put(19.5,1.5){\red{\line(1,1){0.75}}}
\put(20.25,2.25){\red{\line(1,-1){0.75}}}
\put(21,1.5){\red{\line(1,0){1.5}}}
\put(22.5,1.5){\red{\line(1,-1){0.75}}}
\put(23.25,0.75){\red{\line(1,-1){0.75}}}
\put(0,0){\red{\circle*{.17}}}
\put(2.25,0.75){\red{\circle*{.17}}}
\put(5.25,2.25){\red{\circle*{.17}}}
\put(6.75,2.25){\red{\circle*{.17}}}
\put(8.25,3.75){\red{\circle*{.17}}}
\put(9.75,3.75){\red{\circle*{.17}}}
\put(11.25,3.75){\red{\circle*{.17}}}
\put(12.75,3.75){\red{\circle*{.17}}}
\put(14.25,3.75){\red{\circle*{.17}}}
\put(15.75,2.25){\red{\circle*{.17}}}
\put(17.25,2.25){\red{\circle*{.17}}}
\put(20.25,2.25){\red{\circle*{.17}}}
\put(23.25,0.75){\red{\circle*{.17}}}
\put(24,0){\red{\circle*{.17}}}
\end{picture}
}
\vspace*{-5mm}
\caption{Schr\"oder path corresponding to the IRT
shown in Figure~\ref{fig.IRT.big}.
The red dots indicate the endpoints of the segments $\omega(v)$ 
corresponding to the vertices $v$ of the tree.
}
\vspace*{7mm}
   \label{fig.IRT.Schroder}
\end{figure}
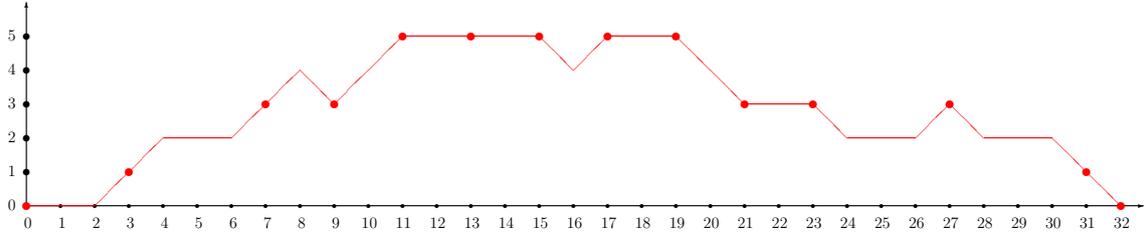

\begin{example}
\rm
The Schr\"oder path $\omega$ corresponding to the IRT
shown in Figure~\ref{fig.IRT.big}
has been drawn in Figure~\ref{fig.IRT.Schroder}.
\myendremark
\end{example}

{\bf Remark.} When $m \ge 1$ (i.e.~the tree consists of more than just a root),
the quantities $\lev(v_1,T), \ldots, \lev(v_m,T)$
describe the heights of a Motzkin path of length $m-1$.
\myendremark

\bigskip

\bigskip

\noindent {\bf Step 2: Definition of the labels $\xi_i$.}
Fix a consistent tree-traversal algorithm $\bfA$.
We will now describe the labels. 
However, notice first that for every vertex $v \neq v_0$,
only the first step in the segment $\omega(v)$ 
starts at an odd height, 
and hence by \reff{def.abc.sbt} it will be
the only step of the segment $\omega(v)$ that gets a choice of labels;
all the other steps have only one choice (namely, $\xi=0$).
Thus, it suffices to assign a label $\xi_v$ to each vertex $v$,
which we define as follows:
\begin{equation}
\xi_v 
\;\eqdef\;
\#\{ w \colon\: p(w) < v < w \text{ and } w<_\bfA v \}\;,
\label{eq.labels.sbt}
\end{equation}
or in other words
\be
\xi_v \;=\; \nid(v,T)
\label{eq.labels.sbt.nid}
\ee
[exactly as in \reff{eq.labels.abt.bt}/\reff{eq.labels.abt.nid.bt}].
To verify that the inequalities~\eqref{def.abc.sbt} are satisfied,
we need to check that 
$0\leq \xi_v \leq \,  \lfloor h/2 \rfloor$,
where $h$ is the starting height of the segment $\omega(v)$.
The lower bound is immediate from the definition of $\xi_v$.
For the upper bound, notice that from Lemma~\ref{lem.sbt.height} it follows
$\lfloor h/2 \rfloor = \lev(v,T)$;
and $\xi_v \leq \lev(v,T)$ is immediate
from the definitions~\eqref{eq.labels.sbt}/\eqref{eq.def.level.vertex.sbt}.

In fact, we also have an interpretation of the difference in terms of the 
statistic croix:
\be
\left\lfloor \dfrac{h}{2}\right\rfloor -\xi_v \;=\; \croix(v,T)
\label{eq.labels.sbt.croix}
\ee
by \eqref{eq.croix.plus.nid}.

\bigskip

\noindent {\bf Step 3: Proof of bijection.}
We now prove that the mapping $T \mapsto (\omega,\xi)$ is a bijection,
by describing the inverse bijection.

To construct the inverse bijection,
we first use the path $\omega$
to determine the label sets $L_{v_0},\ldots,L_{v_m}$
and the node types $N(v_0,T),\ldots,N(v_m,T)$.
This is not entirely trivial (in contrast to the corresponding step
in Section~\ref{subsec.bijection.J.abt});
we present it in Step~3a.
Next we use the labels $\xi$ to glue the vertices together 
and construct our increasing interval-labeled restricted ternary tree
$T\in \IRT_n$;  we present this in Step~3b.
In what follows, we find it convenient to think of the Schr\"oder path
as a word in the alphabet $\{\nearrow, \searrow, \longrightarrow\}$.

If the path $\omega$ consists of $n$ long level steps at height $0$
[i.e., the word $(\longrightarrow)^n$],
then $T$ is the trivial tree ($m=0$) consisting of a root $v_0$
with label set $L_{v_0} = [0,n]$.
In what follows we assume that $\omega$ has at least one rise and one fall.

\medskip

\noindent{\bf Step 3a: Description of vertices.}
We begin by partitioning the label set $[0,n]$ into some intervals,
which will give us our set of vertex labels.
We do this by splitting our path $\omega$ into segments 
$\omega_0,\omega_1,\ldots, \omega_m$ 
such that the word $\omega$ can be factorized as 
\be
\omega \;=\; \omega_0 \omega_1 \cdots \omega_{m-1} \omega_m
   \;;
\ee
then each segment $\omega_i$ will correspond to a vertex $v_i$.
The segments are determined as follows:
\begin{itemize}
    \item $\omega_0$ is the segment of $\omega$ that consists of
all steps before and including the first rise from height $0$.
Thus, it corresponds to a word $(\longrightarrow)^j \nearrow$
for some $j\geq 0$.

    \item The last segment $\omega_m$ is the segment of $\omega$
that consists of all steps of $\omega$ starting at the last fall to height $0$
and including all steps after it.
Thus, it corresponds to a word $\searrow(\longrightarrow)^j$
for some $j\geq 0$.

    \item Now consider the path $\omega^{\flat}$ obtained from $\omega$
by removing the prefix $\omega_0$ and the suffix $\omega_m$:
it starts and ends at height $1$ and never goes below the $x$-axis.
We then split $\omega^{\flat}$ into minimal collections of steps
starting and ending at odd heights. This gives us the factorization 
$\omega^{\flat} = \omega_1\cdots \omega_{m-1}$.
    
\end{itemize}

We will now obtain the label sets for the vertices $v_0,\ldots, v_m$,
as follows:
Consider the unique set-partition
$[0,n] = l_0 \cup l_1 \cup \ldots \cup l_m$
where the set $l_i$ is an interval 
such that $\max l_i<\min l_{i+1}$
and has cardinality given by
\begin{itemize}
   \item If $\omega_0 = (\longrightarrow)^j \nearrow$,
then $|l_0| = j+1$.
   \item If $\omega_m = \searrow (\longrightarrow)^j$,
then $|l_m| = j+1$.
   \item For $1 \le i \le m-1$:
\begin{itemize}
   \item[(a)] If $\omega_i$ contains $j$ long level steps at even height,
      then $|l_i| = j+1$.
   \item[(b)] If $\omega_i$ contains of a single long level step at odd height,
      then $|l_i| = 1$.
\end{itemize}
\end{itemize}
Equivalently,
let $\leven(\omega_i)$ count the number of 
long level steps of $\omega_i$ at an even height.
The cardinality $|l_i|$ is simply given by
$|l_i|\; =\; \leven(\omega_i)+1$.
Then, to the vertex $v_i$ we assign the label set $L_{v_i} = l_i$.

We are now ready to assign node types to our vertices. 
The vertex $v_0$ gets node type $N(v_0,T) = 100$,
and the vertex $v_m$ gets node type $N(v_m,T) = 000$.
For a vertex $v_i$ with $1 \le i \le m-1$,
its node type is determined by the segment $\omega_i$ as follows:
\begin{itemize}
    \item If $\omega_i \,=\: \nearrow (\longrightarrow)^j \nearrow$,
    then $N(v_i,T) = 101$.

    \item If $\omega_i \,=\: \searrow (\longrightarrow)^j \searrow$ ,
    then $N(v_i,T) = 000$.

    \item If $\omega_i \,=\: \nearrow (\longrightarrow)^j \searrow$,
    then $N(v_i,T) = 100$.

    \item If $\omega_i \,=\: \searrow (\longrightarrow)^j \nearrow$,
    then $N(v_i,T) = 001$.

    \item If $\omega_i \,=\: \longrightarrow$,
    then $N(v_i,T) = 010$.
\end{itemize}

\noindent{\bf Step 3b: Constructing the tree using the vertices.}
We now use our labels $\xi$ to glue together the vertices obtained in Step~3a. 
Before doing this, 
notice that each segment $\omega_i$ with $1 \le i \le m-1$
has exactly one step starting at an odd height, namely, its first step. 
From~\eqref{def.abc.sbt}, we know that only this step 
may get a non-unique choice of labels. 
We refer to its label as $\xi_{\omega_i}$.

To describe the construction, we first define a class of intermediate objects
in our bijection:
a {\em slotted interval-labeled restricted ternary tree}\/
is an increasing interval-labeled restricted ternary tree
whose vertices may also have a label set $\{\infty\}$
and we now allow the label set $\{\infty\}$ to be assigned to multiple vertices
(all of which must be leaves).
(We think of the vertices labeled $\{\infty\}$ as slots where new vertices
 may be inserted into the tree in the future.)

Given a tree $T\in \IRT_n$ with vertices $v_0,v_1, \ldots,v_m$,
and $i\in [0,m]$, we define the tree $T|_{i}$ to be the subtree of $T$ 
consisting of vertices $v_0,v_1,\ldots, v_i$ along with their children,
in which the vertices with labels $> \max L_{v_i}$
are relabeled to $\{\infty\}$.
Clearly, $T|_i$ is a slotted interval-labeled restricted ternary tree,
and $T|_m = T$.
By the {\em history}\/ of tree $T$, we will mean the sequence of
slotted interval-labeled restricted ternary trees
$T|_0 \to T|_1 \to \ldots \to T|_{m}$.
Notice that for $i<m$ the number of vertices labeled $\{\infty\}$
in the tree $T|_i$ is $\lev(v_{i+1},T)+1$.

We begin with the tree $T|_0$,
which consists of a root $v_0$ labeled $[0,j]$
with $j = {(|\omega(v_0)| - 1)/2}$
and a left child labeled $\{\infty\}$.
For $i>0$, we will now construct the tree $T|_i$ from $T|_{i-1}$ by
using the label $\xi_{\omega_i}$.
We select the $(\xi_{\omega_i}+1)$-th vertex in the tree-traversal order
chosen in Step~2, and we replace this vertex with $v_i$.
And then we choose to add children (labeled $\{\infty\}$) to $v_i$
according to the node type $N(v_i,T)$.
The resulting tree is $T|_i$.

Finally, at the end of this process,
the tree $T|_{m-1}$ will contain only a single vertex labeled $\{\infty\}$;
we rename this vertex to $v_m$, thereby constructing $T = T|_m$.

\bigskip



\noindent{\bf Step 4: Computation of the weights.}
We can now compute the weights associated to the Schr\"oder path $\omega$
in Theorem~\ref{thm.flajolet_master_labeled_Schroder},
which we recall are
$a_{h,\xi}$ for a rise starting at height $h$ with label $\xi$, 
$b_{h,\xi}$ for a fall starting at height $h$ with label $\xi$, 
and $c_{h,\xi}$ for a long level step starting at height $h$ with label $\xi$.
The weight $\wt(v)$ defined in
\reff{eq.weight.master.sbt.I}--\reff{eq.weight.master.sbt}
is in general a product of factors;
we will distribute this weight among the steps of the segment $\omega(v)$,
as follows: 
\begin{subeqnarray}
   a_{2k-1,\xi}  & = &  \sfahat_{k-1-\xi, \xi}    \\[1mm]
   a_{2k,\xi}    & = &  \mu_k    \\[1mm]
   b_{2k-1,\xi}  & = &  \sfbhat_{k-1-\xi, \xi}    \\[1mm]
   b_{2k,\xi}    & = &  \nu_{k-1}    \\[1mm]
   c_{2k-1,\xi}  & = &  \sff_{k-1-\xi, \xi}    \\[1mm]
   c_{2k,\xi}    & = &  \sfe_k
 \label{eq.defs.abcweights}
\end{subeqnarray}
Let us now verify that these step weights give to each vertex $v$
the correct weights \reff{eq.weight.master.sbt.I}--\reff{eq.weight.master.sbt}
when taking the product over all the steps in the segment $\omega(v)$.

We examine individually each type of step, starting with
the steps starting at odd heights:
\begin{itemize}
    \item[(a)] Rise from height $2k-1$ to height $2k$:
        \begin{itemize}
            \item  By definition of the Schr\"oder path,
we know that this step must correspond to the first step
of a segment $\omega(v)$ for some vertex $v \neq v_0,v_m$ in the tree.
            \item Since the first step of $\omega(v)$ is a rise,
$v$ can have node type $N(v,T) = 100$ or 101. 
In either case we see from \reff{eq.weight.master.sbt}
that it will need one factor $\sfahat$.
            \item From~\eqref{eq.labels.sbt}/\eqref{eq.labels.sbt.nid},
we know that $\xi_v = \nid(v,T)$;
and from~\eqref{eq.labels.sbt.croix}, we get ${k-1-\xi_v = \croix(v,T)}$.
        \end{itemize}
We therefore assign to this step a weight
        \be
        a_{2k-1,\xi} \;=\; 
        \sfahat_{k-1-\xi, \xi}
        \ee
    \item[(b)] Fall from height $2k-1$ to height $2k-2$:
        \begin{itemize}
            \item  By definition of the Schr\"oder path,
we know that this step must correspond to the first step
of a segment $\omega(v)$ for some vertex $v \neq v_0$ in the tree.
            \item Since the first step of $\omega(v)$ is a fall,
$v$ can have node type $N(v,T) = 000$ or $N(v,T) = 001$. 
            In either case it contributes a letter $\sfbhat$.
            \item Once again we have $\xi_v = \nid(v,T)$
and $k-1-\xi_v = \croix(v,T)$.
        \end{itemize}
We therefore assign to this step a weight
        \be
        b_{2k-1,\xi} \;=\; 
        \sfbhat_{k-1-\xi, \xi}
        \ee
    \item[(c)] Long level step at height $2k-1$:
        \begin{itemize}
            \item By definition of the Schr\"oder path,
we know that this step must correspond to $\omega(v)$ for some vertex $v$
with node type $N(v,T) = 010$.  Thus, it contributes the letter $\sff$.
            \item Once again we have $\xi_v = \nid(v,T)$
and $k-1-\xi_v = \croix(v,T)$.
        \end{itemize}
We therefore assign to this step a weight
        \be
        c_{2k-1,\xi} \;=\; 
        \sff_{k-1-\xi, \xi}
        \ee
\end{itemize}
The remaining steps begin at an even height
and hence have $\xi=0$ [cf.~\eqref{def.abc.sbt}].
\begin{itemize}
    \item[(d)] Rise from height $2k$ to height $2k+1$:
        \begin{itemize}
            \item By definition of the Schr\"oder path,
we know that this step must correspond to the last step
of a segment $\omega(v)$ for some vertex $v$ in the tree.
(Here $v=v_0$ is a possibility.)
            \item As the last step of $\omega(v)$ is a rise,
$v$ can have node type $N(v,T) = 101$ or $N(v,T) = 001$. 
In either case, we see from~\reff{eq.weight.master.sbt} that 
we will need one factor $\mu$.
            \item If $N(v,T) = 101$,
the segment $\omega(v)$ started at height $2k-1$,
so from Lemma~\ref{lemma.RT.heights} we have $\lev(v,T)=k-1$;
if $N(v,T)  = 001$,
the segment $\omega(v)$ started at height $2k+1$,
so from Lemma~\ref{lemma.RT.heights} we have $\lev(v,T)=k$.
In either case, \eqref{eq.weight.master.sbt} tells us to assign a weight
        \be
        a_{2k,\xi} \;=\; \mu_{k}
        \ee
\end{itemize}
      \item[(e)] Fall from height $2k$ to height $2k-1$:
        \begin{itemize}
            \item By definition of the Schr\"oder path,
we know that this step must correspond to the last step
of a segment $\omega(v)$ for some vertex $v \neq v_0,v_m$ in the tree.
            \item As the last step of $\omega(v)$ is a fall,
$v$ can have node type $N(v,T) = 000$ or 100. 
In either case, we see from~\eqref{eq.weight.master.sbt} that 
it contributes the letter $\nu$.
            \item If $N(v,T) = 000$,
the segment $\omega(v)$ started at height $2k+1$,
so from Lemma~\ref{lemma.RT.heights} we have $\lev(v,T)=k$;
if $N(v,T)  = 100$,
the segment $\omega(v)$ started at height $2k-1$,
so from Lemma~\ref{lemma.RT.heights} we have $\lev(v,T)=k-1$.
In either case, \eqref{eq.weight.master.sbt} tells us to assign a weight
        \be
        b_{2k,\xi} \;=\; \nu_{k-1}
        \ee
\end{itemize}
      \item[(f)] Long level step at height $2k$:
        \begin{itemize}
            \item By definition of the Schr\"oder path,
we know that this step must correspond to one of the long level steps
in a segment $\omega(v)$ for some vertex $v$ with $|L_v|>1$.
            \item The vertex $v$ can have any node type except $010$.
Thus, we know from \reff{eq.weight.master.sbt.I}--\reff{eq.weight.master.sbt}
that it contributes a letter $\sfe$.
            \item If $v$ is the root, it must have node type 000
(if the tree is trivial) or 100 (if it is nontrivial),
by definition of IRT.  In either case we have $k=0$ and $\lev(v,T) = 0$.
            \item If $v$ is not the root:  If $N(v,T) = 000$ or $001$,
the segment $\omega(v)$ started at height $2k+1$,
so from Lemma~\ref{lemma.RT.heights} we have $\lev(v,T)=k$;
if $N(v,T)  = 100$ or $101$,
the segment $\omega(v)$ started at height $2k-1$,
so from Lemma~\ref{lemma.RT.heights} we have $\lev(v,T)=k-1$.
        \end{itemize}
In all cases, \reff{eq.weight.master.sbt.I}--\reff{eq.weight.master.sbt}
tell us to assign a weight
        \be
        c_{2k,\xi} \;=\; \sfe_{k}
        \ee
The number of such long level steps in the segment $\omega(v)$
is $j= |L_v| - 1$, which agrees with
\reff{eq.weight.master.sbt.I}--\reff{eq.weight.master.sbt}.
\end{itemize}

We can now check, for each of the eight types of segments $\omega(v)$
shown in \reff{eq.RT.step.I}--\reff{eq.RT.step},
that the weight of the segment,
taken as the product of the step weights \reff{eq.defs.abcweights},
indeed coincides in all cases with $\wt(v)$ as defined in
\reff{eq.weight.master.sbt.I}--\reff{eq.weight.master.sbt}.

Putting this all together in Theorem~\ref{thm.flajolet_master_labeled_Schroder},
we obtain a T-fraction with
\begin{eqnarray}
\alpha_{2k-1} &=& 
(\hbox{rise from $2k-2$ to $2k-1$})
   \,\times\,
   (\hbox{fall from $2k-1$ to $2k-2$})
\nonumber\\
&=& \mu_{k-1} \biggl(\sum_{\xi=0}^{k-1} \sfbhat_{k-1-\xi,\xi} \biggr)
\\[2mm]
\alpha_{2k} &=& 
(\hbox{rise from $2k-1$ to $2k$})
   \,\times\,
   (\hbox{fall from $2k$ to $2k-1$})
\nonumber\\
&=&  \biggl(\sum_{\xi=0}^{k-1} \sfahat_{k-1-\xi, \xi} \biggr) \nu_{k-1}
\\[2mm]
\delta_{2k-1} &=& 
\hbox{long level step at height $2k-2$}
\nonumber\\
&=& \sfe_{k-1}
\\[2mm]
\delta_{2k} &=& 
\hbox{long level step at height $2k-1$}
\nonumber\\
&=& \sum_{\xi=0}^{k-1} \sff_{k-1-\xi,\xi}
\label{eq.computation.master.sbt}
\end{eqnarray}
This completes the proof of Theorem~\ref{thm.sbt.master}.
\qed

\medskip

We can now deduce Theorem~\ref{thm.sbt.simple} 
as a corollary:

\proofof{Theorem~\ref{thm.sbt.simple}}
Consider the following substitutions:
\begin{subeqnarray}
\sfahat_{\ell,\ell'} &=& x \\
\sfbhat_{\ell,\ell'} &=& y \\
\sff_{\ell,\ell'} &=& w \\
\mu_{\ell} &=& 1\\
\nu_{\ell} &=& 1\\
\sfe_{\ell} &=& z
\label{eq.substitutions.simple}
\end{subeqnarray}
It is clear that inserting these into~\reff{eq.sbt.master.weights}
yields the continued fraction coefficients~\reff{eq.sbt.simple.weights}.

To finish the proof we will need to argue that substituting
\reff{eq.substitutions.simple} into the polynomials
$Q_n(\bsfahat,\bsfbhat,\bmu,\bnu,\bsfe,\bsff)$
defined in~\reff{eq.def.Qn.irtt}
yields the polynomials $P_n(x,x, y, y, w,z)$
defined in~\reff{eq.Pn.sbt}.
To do this, we show that the summands corresponding to the tree $T\in\IRT_n$
are the same in both polynomials.

Let us first examine the weights of the vertices of $T$
in the polynomials\\
$Q_n(\bsfahat,\bsfbhat,\bmu,\bnu,\bsfe,\bsff)$
--- which are given by
\reff{eq.weight.master.sbt.I}--\reff{eq.weight.master.sbt} ---
under the substitutions \reff{eq.substitutions.simple}.
Consider a vertex $v$ with $L_v = \{l,l+1, \ldots, l+j\}$.
If $l=0$ or $l+j = n$ or both,
by substituting \reff{eq.substitutions.simple}
into \reff{eq.weight.master.sbt.I}--\reff{eq.weight.master.sbt.III} we obtain
\be
\wt(v) \;=\;
\begin{cases}
    z^n \quad\;\, \text{ if $l=0$ and $l+j=n$}  \\
    z^j \;\;\,\,\quad \text{if $l=0$ and $l+j < n$} \\
    y z^j \quad\;\, \text{if $l > 0$ and $l+j = n$}
\end{cases}
  \label{eq.wt.l0jn}
\ee
Note that these three cases correspond, respectively,
to the root of a trivial tree, the root of a nontrivial tree,
and the final vertex of a nontrivial tree
(note that this is a non-root vertex and a leaf).
In the other cases (i.e.\ $l > 0$ and $l+j < n$),
by substituting \reff{eq.substitutions.simple}
into \reff{eq.weight.master.sbt} we obtain
\be
\wt(v) 
\;=\;
\begin{cases}
x z^j \quad \text{ if $N(v,T) = 101$ }\\
y z^j \quad \text{ if $N(v,T) = 000$ }\\
x z^j \quad \text{ if $N(v,T) = 100$ }\\
y z^j \quad \text{ if $N(v,T) = 001$ }\\
w       \quad\quad\!\! \text{ if $N(v,T) = 010$ }
\end{cases}
\ee
These vertices are all non-root and non-final.

In all cases, the power of $z$ contributed by a vertex $v$
is $j = |L_v| - 1$ (the label surplus of $v$),
so by \reff{eq.label.surplus.dist} the total power of $z$
contributed by the tree $T$ is $I_T(\varepsilon)$.
We then get a factor $x$ for each non-root vertex
with node type 101 or 100,
a factor $y$ for each non-root vertex
with node type 000 or 001,
and a factor $w$ for each non-root vertex
with node type 010.
Therefore,
\be
\prod_{v\in V(T)} \wt(v)
\;=\; 
x^{I'_T(100)+I'_T(101)} \, y^{I'_T(000)+I'_T(001)} \, w^{I'_T(010)} \,
z^{I_T(\varepsilon)}\;.
\ee
This matches the contribution of tree $T$
in the polynomial $P_n(x,x, y, y, w,z)$ defined in~\reff{eq.Pn.sbt}.
\qed

\section{Algebraic proofs of Theorems~\ref{thm.Simple.J.abt},
           \ref{thm.sbt.simple} and \ref{thm.sbt.simple.new}}
   \label{sec.algebraic.proofs}

In the preceding section we deduced the ``simple'' continued fractions
(Theorems~\ref{thm.Simple.J.abt} and \ref{thm.sbt.simple})
as special cases of the ``master'' continued fractions
(Theorems~\ref{thm.abt.master.J} and \ref{thm.sbt.master}),
which were proved bijectively.
Here we would like to show how these ``simple'' continued fractions
can be given an alternate, and extremely simple, algebraic proof
using the theory of exponential Riordan arrays and their production matrices.
We begin (Section~\ref{subsec.algebraic.proofs.gen})
by recalling these two concepts and their application
to the enumeration of rooted trees.
Then we apply this theory to prove
Theorem~\ref{thm.Simple.J.abt} (Section~\ref{subsec.algebraic.proofs.RTT})
and Theorems~\ref{thm.sbt.simple} and \ref{thm.sbt.simple.new}
(Section~\ref{subsec.algebraic.proofs.IRT}).
In the latter proof, a crucial role is played by
Proposition~\ref{prop.transformation}.

\subsection{Exponential Riordan arrays, production matrices,
   and the enumeration of rooted trees}
  \label{subsec.algebraic.proofs.gen}

Let $R$ be a commutative ring containing the rationals,
and let $F(t) = \sum_{n=0}^\infty f_n t^n/n!$
and $G(t) = \sum_{n=1}^\infty g_n t^n/n!$ be formal power series
with coefficients in $R$; we set $g_0 = 0$.
Then the \textbfit{exponential Riordan array}
\cite{Deutsch_04,Deutsch_09,Barry_16,Shapiro_22}
associated to the pair $(F,G)$
is the infinite lower-triangular matrix
$\scrr[F,G] = (\scrr[F,G]_{nk})_{n,k \ge 0}$ defined by
\be
   \scrr[F,G]_{nk}
   \;=\;
   {n! \over k!} \:
   [t^n] \, F(t) G(t)^k
   \;.
 \label{def.RFG}
\ee
That is, the $k$th column of $\scrr[F,G]$
has exponential generating function $F(t) G(t)^k/k!$.

Let us now explain briefly about production matrices
\cite{Deutsch_05,Deutsch_09}
\cite[sections~2.2 and 2.3]{forests_totalpos}.
Let $P = (p_{ij})_{i,j \ge 0}$ be an infinite matrix
with entries in a commutative ring $R$,
and assume that $P$ is either row-finite or column-finite
(so that powers of $P$ are well-defined).
Now define a matrix $A = (a_{nk})_{n,k \ge 0}$ by
\be
   a_{nk}  \;=\;  (P^n)_{0k}
 \label{def.iteration}
\ee
(note in particular that $a_{0k} = \delta_{k0}$).
We call $P$ the \textbfit{production matrix}
and $A$ the \textbfit{output matrix},
and we write $A = \scro(P)$.
It is not difficult to see that $AP = \Delta A$,
where $\Delta$ is the matrix with 1 on the superdiagonal and 0 elsewhere.
Note that if $P$ is lower-Hessenberg
(i.e.~vanishes above the first superdiagonal),
then $\scro(P)$ is lower-triangular;
this is the most common case, and will be the case here.
Conversely, when $A$ is lower-triangular with invertible diagonal entries
and $a_{00} = 1$,
then there exists a unique production matrix generating $A$,
and it is the lower-Hessenberg matrix $P = A^{-1} \Delta A$.

The production matrix of an exponential Riordan array $\scrr[F,G]$
is given as follows
(see e.g.~\cite[Theorem~2.19]{forests_totalpos}
\cite[Theorem~6.1]{Shapiro_22} for a proof):

\begin{theorem}[Production matrices of exponential Riordan arrays]
   \label{thm.riordan.exponential.production}
Let $L$ be a lower-triangular matrix
(with entries in a commutative ring $R$ containing the rationals)
with invertible diagonal entries and $L_{00} = 1$,
and let $P = L^{-1} \Delta L$ be its production matrix.
Then $L$ is an exponential Riordan array
if and only~if $P = (p_{nk})_{n,k \ge 0}$ has the form
\be
   p_{nk}
   \;=\;
   {n! \over k!} \: (z_{n-k} \,+\, k \, a_{n-k+1})
 \label{eq.thm.riordan.exponential.production}
\ee
for some sequences $\ba = (a_n)_{n \ge 0}$ and $\bz = (z_n)_{n \ge 0}$
in $R$.  (We set $a_n = z_n = 0$ for $n < 0$.)

More precisely, $L = \scrr[F,G]$ if and only~if $P$
is of the form \reff{eq.thm.riordan.exponential.production}
where the ordinary generating functions
$A(s) = \sum_{n=0}^\infty a_n s^n$ and $Z(s) = \sum_{n=0}^\infty z_n s^n$
are connected to $F(t)$ and $G(t)$ by
\be
   G'(t) \;=\; A(G(t))  \;,\qquad
   {F'(t) \over F(t)} \;=\; Z(G(t))
 \;.
 \label{eq.prop.riordan.exponential.production.1}
\ee
\end{theorem}

\noindent
We refer to $A(s)$ and $Z(s)$ as the
\textbfit{$\bm{A}$-series} and \textbfit{$\bm{Z}$-series}
of the exponential Riordan array $\scrr[F,G]$.

We now apply this theory to the enumeration of rooted trees,
following \cite{latpath_lah}.
Recall first \cite[p.~573]{Stanley_12}
that an \textbfit{ordered tree} (also called {\em plane tree}\/)
is a rooted tree in which the children of each vertex are linearly ordered.
An {\em unordered forest of ordered trees}\/
is an unordered collection of ordered trees.
An \textbfit{increasing ordered tree} is an ordered tree
in which the vertices carry distinct labels from a linearly ordered set
(usually some set of integers) in such a way that
the label of each child is greater than the label of its parent;
otherwise put, the labels increase along every path downwards from the root.
An \textbfit{unordered forest of increasing ordered trees}
is an unordered forest of ordered trees with the same type of labeling.

Now let $\bphi = (\phi_i)_{i \ge 0}$ be indeterminates,
and let $L_{n,k}(\bphi)$ be the generating polynomial for
unordered forests of increasing ordered trees on the vertex set $[n]$,
having $k$ components (i.e.\ $k$ trees),
in which each vertex with $i$ children gets a weight $\phi_i$.
Clearly $L_{n,k}(\bphi)$ is a homogeneous polynomial of degree $n$
with nonnegative integer coefficients;
it is also quasi-homogeneous of degree $n-k$
when $\phi_i$ is assigned weight~$i$.
The polynomials $L_{n,k}(\bphi)$ are called
the \textbfit{generic Lah polynomials};
the lower-triangular matrix $\sfL = (L_{n,k}(\bphi))_{n,k \ge 0}$
is called the \textbfit{generic Lah triangle}.

We now follow \cite[sections~7 and 8]{latpath_lah}.
Define the exponential generating function for trees:
\be
   G(t)  \;=\;  \sum_{n=1}^\infty L_{n,1}(\bphi) \, {t^n \over n!}
   \;.
 \label{def.U}
\ee
It is easy to see that the exponential generating function
for $k$-component unordered forests is then
\be
   {G(t)^k \over k!}  \;=\;  \sum_{n=0}^\infty L_{n,k}(\bphi) \, {t^n \over n!}
   \;.
 \label{eq.Uk}
\ee
Therefore, the generic Lah triangle $\sfL$
is the exponential Riordan array $\scrr[1,G]$.
Furthermore, standard enumerative arguments \cite[Theorem~1]{Bergeron_92}
show that $G(t)$ satisfies the ordinary differential equation
\be
   G'(t) \;=\; \Phi(G(t))
   \;,
 \label{eq.bergeron.ODE}
\ee
where $\Phi(s) \eqdef \sum_{k=0}^\infty \phi_k s^k$
is the ordinary generating function for $\bphi$.
Therefore, from \reff{eq.prop.riordan.exponential.production.1}
we see that $A(s) = \Phi(s)$.

Here we would like to get the trees into column~0 rather than column~1
of the output matrix, and shifted down to start at $n=0$ rather than $n=1$.
This is easy:  it suffices to consider
the exponential Riordan array $\scrr[F,G]$ with
\be
   F(t)  \;\eqdef\; G'(t)
     \;=\;  \sum_{n=0}^\infty L_{n+1,1}(\bphi) \, {t^n \over n!}
   \;.
\ee
Differentiating \reff{eq.bergeron.ODE}, we deduce that
\be
   G''(t)  \;=\;  \Phi'(G(t)) \: G'(t)
\ee
and hence that
\be
   {F'(t) \over F(t)} \;=\; \Phi'(G(t))
   \;.
\ee
Therefore, from \reff{eq.prop.riordan.exponential.production.1}
we see that $Z(s) = \Phi'(s)$.

Inserting these formulae for $A(s)$ and $Z(s)$ into
\reff{eq.thm.riordan.exponential.production},
we obtain the extraordinarily simple formula
\be
   p_{nk}
   \;=\;
   {(n+1)! \over k!} \: \phi_{n-k+1}
 \label{eq.pnk.trees}
\ee
for the production matrix of $\scrr[G',G]$.

\subsection{Application to increasing restricted ternary trees:
       Proof of Theorem~\ref{thm.Simple.J.abt}}
   \label{subsec.algebraic.proofs.RTT}

It is now easy to apply the foregoing theory to restricted ternary trees
with the weights \reff{eq.Pn.abt},
in which the variables $x_1, y_1, x_2, y_2, w$
are associated to the node types
$101$, $000$, $100$, $001$, $010$, respectively.
It suffices to take
\be
   \phi_0 \:=\: y_1 \,,\quad
   \phi_1 \:=\: x_2 + y_2 + w \,,\quad
   \phi_2 \:=\: x_1 \,,\quad
   \phi_i \:=\: 0 \;\hbox{for $i \ge 3$}
   \;.
\ee
The production matrix \reff{eq.pnk.trees} is then tridiagonal with
\begin{subeqnarray}
   p_{n,n+1}  & = &  y_1  \\[2mm]
   p_{n,n}    & = &  (n+1) (x_2 + y_2 + w)  \\[2mm]
   p_{n,n-1}  & = &  n(n+1) x_1
\end{subeqnarray}
This generates Motzkin paths with these weights
for rises, level steps and falls, respectively.
On the other hand, for walks that end at height~0
we can transfer all the weights from rises to falls;
we thus obtain a J-fraction with coefficients
\be
   \gamma_n  \;=\;  (n+1) (x_2 + y_2 + w)  \,,\qquad
   \beta_n  \;=\;  n(n+1) x_1 y_1
 \;.
\ee
These are precisely the coefficients \reff{eq.def.simple.abt.J};
we have therefore proven Theorem~\ref{thm.Simple.J.abt}.

\subsection{Application to increasing interval-labeled restricted ternary trees:
       Proof of Theorems~\ref{thm.sbt.simple} and \ref{thm.sbt.simple.new}}
   \label{subsec.algebraic.proofs.IRT}


An increasing interval-labeled restricted ternary tree $T'$
on the label set $[0,n']$ can be obtained from an
increasing restricted ternary tree $T$ on the vertex set $[n]$, as follows:
\begin{itemize}
   \item If $n=0$, then $T$ is the empty tree;
      it gets a weight 1 in \reff{eq.Pn.abt}.
      The tree $T'$ is a trivial tree consisting of only a root
      with label set $[0,n']$;
      it gets a weight $z^{n'} = z^{I_{T'}(\varepsilon)}$ in \reff{eq.Pn.sbt}.
   \item If $n \ge 1$, the tree $T$ is nonempty.
      The tree $T'$ is obtained from $T$ by creating a root with
      label set $[0,j]$ for some $j \ge 0$, and making the root of $T$
      be the left child of the root of $T'$;
      also, every vertex of $T$ {\em other than one with a middle child}\/
      can become interval-labeled.
      The weight of $T'$ in \reff{eq.Pn.sbt} is the same as the
      weight of $T$ in \reff{eq.Pn.abt},
      multiplied by $z^{I_{T'}(\varepsilon)}$.
\end{itemize} 
It follows that the ordinary generating functions of
the polynomials \reff{eq.Pn.abt} and \reff{eq.Pn.sbt}
are related as follows:
\be
   \sum_{n=0}^\infty  P_n(x_1,x_2,y_1,y_2,z,w) \, t^n
   \;=\;
   {1 \over 1-zt}
   \sum_{n=0}^\infty  P_n\Bigl( {x_1 \over 1-zt}, {x_2 \over 1-zt},
                                {y_1 \over 1-zt}, {y_2 \over 1-zt}, w \Bigr)
          \, t^n
   \;.
\label{eq.ogf.rtt.irtt}
\ee
Theorem~\ref{thm.sbt.simple.new} is then an immediate consequence of
Theorem~\ref{thm.simple.abt.T} together with
Proposition~\ref{prop.transformation}:
the key fact is that $x_1 + y_1 = x_2 + y_2$,
so that $\delta_{2k}$ does not receive any factor $1/(1-zt)$.
In particular, specializing $x_1=x_2=x$ and $y_1=y_2=y$
gives us Theorem~\ref{thm.sbt.simple}.

\section{Other interpretations}
\label{sec.other.interpretations}

In this section we use bijections to reinterpret our results
in terms of other combinatorial objects:
increasing binary trees are in bijection with permutations
(Section~\ref{subsec.other.interpretations.1}),
and increasing restricted ternary trees
are in bijection with binary free multilabeled increasing trees
(Section~\ref{subsubsec.free}).

\subsection{Increasing binary trees $\simeq$ Permutations}
   \label{subsec.other.interpretations.1}

It is well known \cite[pp.~44--45]{Stanley_12}
that the set $\mathcal{B}_n$ of increasing binary trees
on the vertex set $[n]$ is in bijection
with the set $\mathfrak{S}_n$ of permutations of $[n]$.
In this section, we will translate our statistics for increasing binary trees
to permutation statistics via this bijection.
This will allow us to translate our continued fractions
for increasing binary trees
(Section~\ref{subsec.bt})
to continued fractions counting various statistics on permutations:
namely, linear statistics and vincular patterns.

Let $\sigma = (\sigma_1 \cdots \sigma_n) \in \mathfrak{S}_n$
be a permutation of $[n]$, which we shall consider principally as a word.
We declare $\sigma_0 =\sigma_{n+1} = 0$.
An index $i\in [n]$ (or a letter $\sigma_i \in [n]$) is called a
\begin{itemize}
    \item {\em peak} (pk) if 
    $\sigma_{i-1}<\sigma_i>\sigma_{i+1}$,
    
    \item {\em valley} (val) if 
    $\sigma_{i-1}>\sigma_i<\sigma_{i+1}$,
    
    \item {\em double ascent} (dasc) if 
    $\sigma_{i-1}<\sigma_i<\sigma_{i+1}$,
    
    \item {\em double descent} (ddes) if 
    $\sigma_{i-1}>\sigma_i>\sigma_{i+1}$.
\end{itemize}
Clearly every index $i$ belongs to one of these four types;
we refer to this classification as the {\bf linear classification}.

\medskip

{\bf Remark.}
The boundary condition $\sigma_0 =\sigma_{n+1} = 0$
plays a key role in our definition of the linear classification.
Other boundary conditions can also be used:
for instance, Han, Mao and Zeng \cite{Han_21} use
$\sigma_0 = 0$, $\sigma_{n+1} = n+1$.
Different boundary conditions give rise to {\em different}\/
linear classifications.
\myendremark

\smallskip

Next we introduce certain permutation statistics 
given by the occurrence of some vincular patterns \cite{Babson_00}
(see also the survey article \cite{Steingrimsson_10}).
We recall that a vincular pattern is similar in meaning
to an ordinary permutation pattern,
except that the absence of a dash indicates that the two letters
are required to be consecutive in the word.
For instance, the pattern 31--2 means that we have
three letters $\ell_1 < \ell_2 < \ell_3$
such that $\ell_3$ occurs immediately before $\ell_1$ in the word,
while $\ell_2$ occurs after $\ell_1$
(possibly but not necessarily immediately after it).\footnote{
   There is an alternative (and perhaps preferable) notation
   for vincular patterns in which terms that must be adjacent
   are underlined.
   For instance, the pattern 31--2 would be written in this notation as
   $\underline{31}2$.
   This notation has the advantage of reducing to the ordinary notation
   for permutation patterns when there are no underlinings.
}

For a letter $\ell\in [n]$, 
we define
$(\textrm{31--2})(\ell,\sigma)$ and $(\textrm{2--13})(\ell,\sigma)$
as follows:
\begin{eqnarray}
    (\textrm{31--2})(\ell,\sigma) 
    \;\eqdef\;
    \#\{j \colon\: 1 < j < \sigma^{-1}_\ell \hbox{ and } \sigma_j < \ell < \sigma_{j-1}\}
    \\[2mm]
    (\textrm{2--13})(\ell,\sigma) 
    \;\eqdef\;
    \#\{j \colon\: \sigma^{-1}_{\ell} < j < n \hbox{ and } \sigma_j < \ell < \sigma_{j+1}\}
\end{eqnarray}
Thus, $(\textrm{31--2})(\ell,\sigma)$ counts the number of 
occurrences of the vincular pattern $\textrm{31--2}$
in which the letter $\ell$ is the $2$ in the pattern,
while $\sigma_j$ and $\sigma_{j-1}$ are the $1$ and the $3$.
Similarly, $(\textrm{2--13})(\ell,\sigma)$ counts the number of 
occurrences of the vincular pattern $\textrm{2--13}$
in which the letter $\ell$ is the $2$ in the pattern,
while $\sigma_j$ and $\sigma_{j+1}$ are the $1$ and the $3$.
Note that these definitions do not require any boundary condition.

Let $\Phi_n\colon \mathcal{B}_n \to \mathfrak{S}_n$ 
denote the reverse bijection in \cite[pp.~44--45]{Stanley_12}
from increasing binary trees to permutations.
The permutation $\Phi_n(T)$ in word form
is obtained by writing out the vertices of tree $T$
as per the inorder (= symmetric) traversal:
that~is, (left, root, right), implemented recursively.
The letters of the word $\Phi_n(T)$
thus correspond to the vertex labels of the tree $T$,
and the order of those letters in the word $\Phi_n(T)$
corresponds to the inorder traversal $\bfA$ on those vertex labels:
$\ell <_\bfA \ell'$ $\Longleftrightarrow$
$\sigma^{-1}_\ell < \sigma^{-1}_{\ell'}$.
%
%
%
Furthermore, the table in \cite[p.~45]{Stanley_12}
records the correspondence between the node types of a tree $T\in \mathcal{B}_n$
and the linear classification in the permutation $\Phi_n(T)$,
using the boundary condition $\sigma_0 =\sigma_{n+1} = 0$:
\begin{eqnarray}
    \begin{tabular}{c|c}
        \hline
        Node type $N(\sigma_i,T)$ in tree $T$ &
        Linear classification of index $i$ in $\sigma = \Phi_n(T)$ \\
        \hline
        00 & Peak \\
        11 & Valley \\
        10 & Double descent \\
        01 & Double ascent \\
        \hline
    \end{tabular}
\label{eq.tab.linear.translate}
\end{eqnarray}

Next, we translate the crossing and nesting statistics;
we stress that they are defined with respect to the inorder traversal.

\begin{proposition}
   \label{prop.translate.croix.nid}
Let $T\in \mathcal{B}_n$ be an increasing binary tree
and let $\sigma = \Phi_n(T)$.
Then the following identities hold:
\begin{subeqnarray}
    \nid(\ell,T)   &=& (\mbox{\rm 31--2})(\ell,\sigma)\\
    \croix(\ell,T) &=& (\mbox{\rm 2--13})(\ell,\sigma)
\label{eq.translate.croix.nid}
\end{subeqnarray}
where $\nid$ and $\croix$ are defined with respect to the inorder traversal.
\end{proposition}

\proof
We recall that 
\be
\croix(\ell,T) \;=\; \#\{w\colon\: p(w)<\ell<w \text{ and } \ell <_{\bfA} w\}
\ee
where $p(w)$ denotes the parent of $w$ in $T$,
and $\bfA$ is the inorder traversal.
Let $w$ be a vertex in $T$ that contributes to the right-hand side above.
We first claim that either $w$ is a right child of its parent
or it has an ancestor that is the right child of its parent.
For if not, then $w$ lies on the leftmost branch of $T$,
and the only possibility of having $\ell<_{\bfA} w$
is if $\ell$ occurs in the left subtree of $w$
(recall that we are using inorder traversal);
but this contradicts $\ell<w$ since the tree $T$ is increasing.

Let $\widehat{w}$ be the largest ancestor of $w$
(that is, the ancestor closest to $w$)
that contains $w$ in its right subtree.
In other words, $\widehat{w}$ is the parent vertex of the first right edge
that occurs on the path from $w$ to the root.
In particular, $\widehat{w} = p(w)$
if $w$ is the right child of its parent $p(w)$.

Let $\widetilde{w}$ denote the leftmost descendant of $w$ in $T$
(that is, $\widetilde{w}$ is a descendant of $w$ that does not 
have a left child, and the path downward from $w$ to $\widetilde{w}$
consists only of left edges).
In particular, $\widetilde{w}=w$ if $w$ does not have a left child.

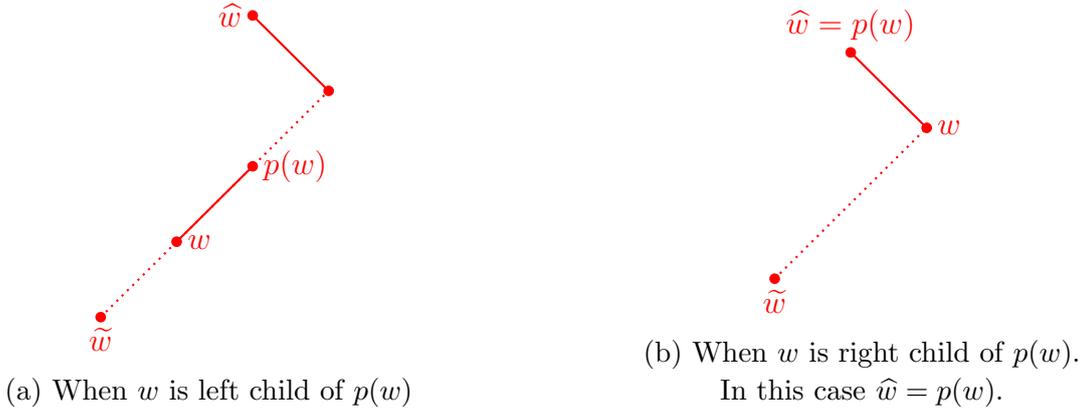
\begin{figure}[t]
\centering
	\begin{tabular}[b]{c}
	\begin{tikzpicture}
	\draw[red,thick,dotted] (0,0) -- (3,3);
	\draw[red,thick] (1,1) -- (2,2);
	\draw[red,thick] (3,3) -- (2,4);
	\fill[red] (0,0) circle (2pt) node[below] {$\widetilde{w}$};
	\fill[red] (1,1) circle (2pt) node[right] {$w$};
	\fill[red] (2,2) circle (2pt) node[right] {$p(w)$};
	\fill[red] (2,4) circle (2pt) node[left] {$\widehat{w}$};
	\fill[red] (3,3) circle (2pt) node {};
	\end{tikzpicture}\\
	\small{(a) When $w$ is left child of $p(w)$}
	\end{tabular} \qquad\qquad\qquad
	\begin{tabular}[b]{c}
	\begin{tikzpicture}
        \draw[red,thick,dotted] (0,0) -- (2,2);
        \draw[red,thick] (2,2) -- (1,3);
	\fill[red] (0,0) circle (2pt) node[below] {$\widetilde{w}$};
        \fill[red] (2,2) circle (2pt) node[right] {$w$};
	\fill[red] (1,3) circle (2pt) node[above] {$\widehat{w}=p(w)$};
        \end{tikzpicture}\\
	\small{(b) When $w$ is right child of $p(w)$.}\\
	\small{In this case $\widehat{w} = p(w)$.}
	\end{tabular} 
\caption{An illustration of $\widehat{w}$ and $\widetilde{w}$ for a given $w$ 
contributing in the definition of $\croix(\ell,T)$.
}
\label{fig.what.wtilde}
\end{figure}

See Figure~\ref{fig.what.wtilde} for an illustration of
$\widehat{w}$ and $\widetilde{w}$ for a given $w$.

It is immediate that $\widehat{w}\leq p(w)<\ell<w\leq \widetilde{w}$.
Let $j = \sigma^{-1}_{\widehat{w}}$.
It is clear, from the definition of inorder traversal,
that the letter $\widehat{w}$ is 
immediately followed by $\widetilde{w}$ in $\sigma$;
therefore $\sigma_{j+1} = \widetilde{w}$.
Thus, the triple
$(\ell_1,\ell_2,\ell_3) = (\widehat{w},\ell,\widetilde{w})$
forms a vincular pattern $\textrm{2--13}$ in $\sigma$.

Also, notice that if we have two distinct vertices $w_1 \neq w_2$
contributing to $\croix(\ell,T)$, 
then necessarily $\widehat{w_1}\neq \widehat{w_2}$
(since for a given vertex $\widehat{w}$
 there is at most vertex $w$ in the relevant chain
 that satisfies $p(w) < \ell < w$).
Therefore, for every $w$ that contributes to $\croix(\ell,T)$,
we have obtained a {\em distinct}\/ $j$
that contributes to $(\textrm{2--13})(\ell,\sigma)$;
this shows that
\be
\croix(\ell,T) \leq  (\textrm{2--13})(\ell,\sigma)\;.
\ee

On the other hand, given $j$ such that 
$\sigma^{-1}_\ell<j<n$ and $\sigma_j<\ell<\sigma_{j+1}$,
we observe that $\sigma_j$ must be an ancestor of $\sigma_{j+1}$ in tree $T$.
For if not, let $v$ be the closest common ancestor
of $\sigma_j$ and $\sigma_{j+1}$;
then the vertices $\sigma_j$ and $\sigma_{j+1}$
are in different subtrees of $v$.
As $\sigma$ is obtained from $T$ by inorder traversal,
the letter $v$ must occur between $\sigma_j$ and $\sigma_{j+1}$ in $\sigma$,
which is a contradiction.

As $\sigma_j$ is an ancestor of $\sigma_{j+1}$,
it is clear that $\sigma_{j+1}$ is in the right subtree of $\sigma_j$
and is the leftmost vertex on the right subtree of $\sigma_j$.
Also notice that the vertex $\ell$ is not in the right subtree of $\sigma_j$,
as that would contradict $\sigma^{-1}_\ell<j$.
In particular,
$\ell$ does not occur in the path from $\sigma_j$ to $\sigma_{j+1}$.
Since $\sigma_j<\ell<\sigma_{j+1}$,
on the path from $\sigma_j$ to $\sigma_{j+1}$
there is a unique vertex $w$ such that $p(w)<\ell<w$.
This vertex $w$ is in the right subtree of $\sigma_j$ and hence
$\sigma_j <_\bfA w$.
As we also have $\ell<_\bfA \sigma_j$
(since this is equivalent to $\sigma^{-1}_\ell<j$),
we get $\ell <_\bfA w$.
Furthermore, since $w$ lies in the left branch
of the right child of $\sigma_j$,
we have $\sigma_j = \widehat{w}$;
therefore $j = \sigma^{-1}_{\widehat{w}}$ is uniquely determined by $w$.
Thus, for every $j$ that contributes to $(\textrm{2--13})(\ell,\sigma)$,
we have obtained a {\em distinct}\/ $w$ that contributes to $\croix(\ell,T)$.
This shows that
\be
 (\textrm{2--13})(\ell,\sigma) \;\le\; \croix(\ell,T)\;.
\ee
This finishes the proof of~(\ref{eq.translate.croix.nid}b).

The proof for~(\ref{eq.translate.croix.nid}a) is 
obtained by simply switching left and right in the above proof;
we leave it as an exercise for the reader.
\qed

{\bf Remark.}
The fact that $\nid(\ell,T)$ translates to $(\textrm{31--2})(\ell,\sigma)$ 
is implicitly mentioned by Viennot \cite[Chapter~4b]{Viennot_abjc1}
using ``$x$-decomposition'' on permutations.\footnote{
   Video link address: \url{https://www.youtube.com/watch?v=Cp8adiOL_6Q&t=865}
}
\myendremark

\medskip

Now define the polynomials $P^\star_n(\bsfa, \bsfb,\bsfc,\bsfd)$ as
\begin{eqnarray}
\hspace*{-6mm}
P^\star_n(\bsfa, \bsfb,\bsfc,\bsfd)
&=&
\sum_{\sigma\in \mathfrak{S}_n}
\prod_{\ell\in \Val(\sigma)}
\sfa_{(\textrm{2--13})(\ell,\sigma),(\textrm{31--2})(\ell,\sigma)}
\prod_{\ell\in \Pk(\sigma) }
\sfb_{(\textrm{2--13})(\ell,\sigma),(\textrm{31--2})(\ell,\sigma)} \times
\nonumber\\
&& \quad\;\;\, \prod_{\ell\in \Ddes(\sigma) }
\!\sfc_{(\textrm{2--13})(\ell,\sigma),(\textrm{31--2})(\ell,\sigma)}
\prod_{ \ell\in \Dasc(\sigma)}
\!\!\sfd_{(\textrm{2--13})(\ell,\sigma),(\textrm{31--2})(\ell,\sigma)}
\label{eq.def.Pstarn.bt}
\end{eqnarray}
where $\Val(\sigma)$ denotes the set of all valley {\em letters}\/ of $\sigma$,
and likewise for the others.
From~\reff{eq.tab.linear.translate}
and Proposition~\ref{prop.translate.croix.nid}
we immediately obtain the following theorem:

\begin{theorem}[Master J- and T-fractions for permutations with linear statistics]
   \label{thm.permutations}
The polynomials~$P^\star_n(\bsfa, \bsfb,\bsfc,\bsfd)$
defined in \reff{eq.def.Pstarn.bt}
and the polynomials~$Q_n(\bsfa, \bsfb,\bsfc,\bsfd)$
defined in~\reff{eq.def.Qn.bt}
are equal:
\be
P^\star_n(\bsfa, \bsfb,\bsfc,\bsfd)
\;=\;
Q_n(\bsfa, \bsfb,\bsfc,\bsfd)\;.
\ee
Therefore, the ordinary generating function of the polynomials
$P^\star_n(\bsfa, \bsfb,\bsfc,\bsfd)$
has the same J-fraction of Theorem~\ref{thm.bt.master}
and the same T-fraction of Theorem~\ref{thm.bt.master.T}.
\end{theorem}

\begin{example}
\rm
For $\sigma = 57316284$, the quantities 
$(\textrm{31--2})(\ell,\sigma)$ and $(\textrm{2--13})(\ell,\sigma)$
are as follows:
%
\begin{equation}
\begin{tabular}{c|c|c}
	$\ell$ & 
	$(\textrm{31--2})(\ell,\sigma)$
	&
	$(\textrm{2--13})(\ell,\sigma)$\\
	\hline
	1 & 0 & 0\\
	2 & 1 & 0\\
	3 & 0 & 2\\
	4 & 2 & 0\\
	5 & 0 & 2\\
	6 & 1 & 1\\
	7 & 0 & 1\\
	8 & 0 & 0	
\end{tabular}
  \label{eq.tab.perm.vincular}
\end{equation}
Notice that when $T$ is increasing binary tree shown in Figure~\ref{fig.bt},
we have $\Phi_8(T) = \sigma = 57316284$.
Comparing \reff{eq.tab.bt.croix.nid} and~\reff{eq.tab.perm.vincular},
we can verify that~\reff{eq.translate.croix.nid} holds for this example.
\myendremark
\end{example}

{\bf Remark.}
Substituting 
$\sfa_{\ell,\ell'}=1$,
$\sfb_{\ell,\ell'} = u$,
$\sfc_{\ell,\ell'} = w$
and $\sfd_{\ell,\ell'} = v$
gives us \cite[Theorem~3A]{Flajolet_80}.
\myendremark

\medskip

Here is an important special case: by setting
\be
 \sfa_{\ell,\ell'} \,=\, 
 \sfb_{\ell,\ell'} \,=\, 
 \sfc_{\ell,\ell'} \,=\, 
 \sfd_{\ell,\ell'} \,=\, p^\ell q^{\ell'}
\ee
in Corollary~\ref{cor.bt.master.T}
(which is a special case of Theorem~\ref{thm.bt.master.T})
and using Theorem~\ref{thm.permutations},
we obtain an S-fraction
for the joint generating function of the statistics 2--13 and 31--2:

\begin{corollary}
   \label{cor.permutations}
We have the S-fraction
\be
   \sum_{n=0}^\infty \, \sum_{\sigma\in\mathfrak{S}_n}
       p^{(\textrm{\rm 2--13})(\sigma)} \, q^{(\textrm{\rm 31--2})(\sigma)} \, t^n
   \;=\;
   \cfrac{1}{1 -  \cfrac{[1]_{p,q} \, t}{1 - \cfrac{[1]_{p,q} \, t}{ 1 -
\cfrac{[2]_{p,q} \, t }{1- \cfrac{[2]_{p,q} \, t}{1-\ldots}
}}}}
\ee
with coefficients
\be
   \alpha_{2k-1} \;=\; \alpha_{2k} \;=\; [k]_{p,q}
   \;.
\ee
\end{corollary}


Furthermore, by combining this continued fraction
with previous work of Claesson and Mansour \cite{Claesson_02},
we can learn more about the joint distributions of various vincular patterns.
For starters, Claesson \cite[Proposition~1]{Claesson_01} showed that
the four vincular patterns
$\textrm{2--13}$, $\textrm{2--31}$, $\textrm{13--2}$, $\textrm{31--2}$
are equidistributed.
Now consider the eight possible ordered pairs formed by taking
one pattern of the form 2--ab and one of the form ab--2:
\begin{samepage}
\begin{quote}
\begin{itemize}
   \item[1.]  $(\textrm{2--13},\, \textrm{31--2})$
   \item[2.]  $(\textrm{31--2},\, \textrm{2--13})$
   \item[3.]  $(\textrm{2--31},\, \textrm{13--2})$
   \item[4.]  $(\textrm{13--2},\, \textrm{2--31})$  \\[-2mm]
   \item[5.]  $(\textrm{2--13},\, \textrm{13--2})$
   \item[6.]  $(\textrm{31--2},\, \textrm{2--31})$
   \item[7.] $(\textrm{2--31},\, \textrm{31--2})$
   \item[8.] $(\textrm{13--2},\, \textrm{2--13})$
\end{itemize}
\end{quote}
\end{samepage}
It is easy to see that the first four of these are equidistributed,
and also the last four:
\begin{itemize}
   \item The equivalences $1 \leftrightarrow 3$,\, $2 \leftrightarrow 4$,\,
                          $5 \leftrightarrow 7$ and $6 \leftrightarrow 8$
         are obtained by using complementation $\sigma \mapsto \sigma^c$
         (that is, mapping {\em letters}\/ $\ell \mapsto n+1-\ell$).
   \item The equivalences $1 \leftrightarrow 2$,\, $3 \leftrightarrow 4$,\,
                          $5 \leftrightarrow 6$ and $7 \leftrightarrow 8$
         are obtained by using reversal $\sigma \mapsto \sigma^r$
         (that is, mapping {\em indices}\/ $i \mapsto n+1-i$).
\end{itemize}

But here is the surprise:  it turns out that {\em all eight}\/
ordered pairs are equidistributed!
This follows from the fact that
Claesson and Mansour's continued fraction
for the generating function
$\sum_{n=0}^\infty \sum_{\sigma\in\mathfrak{S}_n}
p^{(\textrm{2--31})(\sigma)} \, q^{(\textrm{31--2})(\sigma)} \, t^n$
\cite[Theorem~22 with $x=y=1$]{Claesson_02}
coincides with ours in Corollary~\ref{cor.permutations}
for $\sum_{n=0}^\infty \sum_{\sigma\in\mathfrak{S}_n}
p^{(\textrm{2--13})(\sigma)} \, q^{(\textrm{31--2})(\sigma)} \, t^n$.
The combination of these two continued fractions therefore proves:\footnote{
   Let us remark that Vajnovszki \cite{Vajnovszki_18}
   proved a special case of the equivalence $3 \leftrightarrow 5$,
   and thus of Proposition~\ref{prop.permutations.equidistribution}:
   namely, he proved that $\textrm{\rm 2--13}$ and $\textrm{\rm 2--31}$
   are equidistributed
   {\em among permutations avoiding 132}\/
   (or equivalently, avoiding $\textrm{\rm 13--2}$)
   --- that is, for the subclass for which
   $(\textrm{\rm 13--2})(\sigma) = 0$.
   We thank Anders Claesson for drawing our attention to
   Vajnovszki's work.
}

\begin{proposition}
   \label{prop.permutations.equidistribution}
In permutations of $[n]$, the pairs of statistics
$(\textrm{\rm 2--13},\, \textrm{\rm 31--2})$
and
$(\textrm{\rm 2--31},\, \textrm{\rm 31--2})$
are equidistributed.
\end{proposition}

\noindent
But we do not know any direct bijective proof of this equidistribution;
we therefore pose it as an open problem:

\begin{openproblem}
\rm
Find a bijective proof of Proposition~\ref{prop.permutations.equidistribution}.
\end{openproblem}

Indeed, we can go farther, by considering the joint distribution
of all four vincular patterns.
Define the polynomials in four variables
\be
   P_n(p,q,r,s)
   \;=\;
   \sum_{\sigma\in\mathfrak{S}_n}
       p^{(\textrm{\rm 13--2})(\sigma)} \, q^{(\textrm{\rm 31--2})(\sigma)} \,
       r^{(\textrm{\rm 2--13})(\sigma)} \, s^{(\textrm{\rm 2--31})(\sigma)}
   \;.
\ee
Then the identity $P_n(p,q,r,s) = P_n(q,p,s,r)$
is a consequence of complementation symmetry,
while the identity $P_n(p,q,r,s) =  P_n(s,r,q,p)$
is a consequence of reversal symmetry.
By combining these we obtain the $\Z_2 \times \Z_2$ symmetry
\be
   P_n (p, q, r, s) \;=\; P_n (q, p, s, r) \;=\; P_n (s, r, q, p) \;=\;
   P_n (r, s, p, q)
   \;.
 \label{eq.permutations.symmetry}
\ee
Furthermore, it can be checked that for $n \ge 5$
these are the only permutations of the four variables
that leave the polynomial $P_n$ invariant.

But by setting one of the variables equal to 1,
we obtain empirically some interesting identities,
which we state as a conjecture:

\begin{conjecture}[Trivariate symmetries on vincular patterns]
\hfill\break
   \label{conj.trivariate}
$\!\!\!$We have the relations
\begin{subeqnarray}
   P_n (1, q, r, s)  & = &  P_n (1, q, s, r)   \\[1mm]
   P_n (p, 1, r, s)  & = &  P_n (p, 1, s, r)   \\[1mm]
   P_n (p, q, 1, s)  & = &  P_n (q, p, 1, s)   \\[1mm]
   P_n (p, q, r, 1)  & = &  P_n (q, p, r, 1) 
 \label{eq.conj.trivariate}
\end{subeqnarray}
\end{conjecture}

\noindent
We have verified Conjecture~\ref{conj.trivariate} for $n \le 11$.
Please note that the four conjectured relations are equivalent
by virtue of the symmetries \reff{eq.permutations.symmetry},
so it suffices to prove one of them.
Please note also that Proposition~\ref{prop.permutations.equidistribution}
is the identity $P_n(1,q,r,1) = P_n(1,q,1,r)$,
which is the $s=1$ special case of (\ref{eq.conj.trivariate}a).
So Conjecture~\ref{conj.trivariate} generalizes
Proposition~\ref{prop.permutations.equidistribution}.

\bigskip

Finally, it is worth comparing our master polynomials~\reff{eq.def.Pstarn.bt}
with the polynomials \cite[eqs.~(3.1) and (3.2)]{Han_21}
of Han, Mao and Zeng, for which they demonstrate a master J-fraction
\cite[eq.~(1.19)]{Han_21}.
There are three differences between their polynomials and ours:
\begin{itemize}
   \item Their polynomials involve the joint distribution
of the patterns $(\textrm{2--31},\,\textrm{31--2})$
--- as in Claesson and Mansour \cite[Theorem~22]{Claesson_02} ---
whereas ours involve the joint distribution of
$(\textrm{2--13},\,\textrm{31--2})$.
   \item They took the boundary conditions
$\sigma_0 = 0$, $\sigma_{n+1} = n+1$,
while we took the boundary conditions $\sigma_0 = \sigma_{n+1} = 0$.
This difference in the boundary conditions
does not affect the meanings of the vincular patterns,
but it does affect the meanings of the linear classification
and hence of the master polynomials.
   \item Their linear classification was more refined than ours:
they refined double ascents into
{\em foremaxima}\/ (double ascents that are also records) and the rest.
\end{itemize}
Finally, their master J-fraction is different from ours,
because it concerns the ogf of $Q_n$
and coincides with that of \cite[Theorem~2.9]{Sokal-Zeng_masterpoly},
which generalizes the J-fraction for the sequence $(n!)_{n\geq 0}$,
while our master J-fraction in Theorem~\ref{thm.bt.master}
concerns the ogf of $Q_{n+1}$
and hence generalizes the J-fraction for the sequence $((n+1)!)_{n\geq 0}$.

If we call \cite[Theorems~3.1 and~3.2 with Theorem~1.10]{Han_21}
the \textbfit{first and second master J-fractions for permutations
with linear statistics},
then we can call the result obtained by combining
Theorem~\ref{thm.bt.master} with Theorem~\ref{thm.permutations}
the \textbfit{third master J-fraction for permutations with linear statistics}.
Also, the result obtained by combining
Theorem~\ref{thm.bt.master.T} with Theorem~\ref{thm.permutations}
can be called the
\textbfit{master T-fraction for permutations with linear statistics}.

\subsection{Increasing restricted ternary trees $\simeq$
            Binary free multilabeled increasing trees}
   \label{subsubsec.free}

In \cite[Section~5.2, Example~5]{Kuba_16}, 
Kuba and Panholzer consider binary free multilabeled increasing trees.
A \textbfit{binary free multilabeled increasing tree} with label set $[n]$
is a binary tree $T$ in which each vertex $v$ is assigned
a nonempty set of labels $L_v$ such that
(a) $\{L_v\}_{v \in V(T)}$ is a set-partition of $[n]$,
and
(b) every label of a child is larger than every label of its parent.

There is a simple bijection between binary free multilabeled increasing trees
on the label set $[n]$
and increasing restricted ternary trees on the vertex set $[n]$.
In fact, this bijection is a slight modification of the bijection
in \cite[Theorem~10]{Kuba_16} that is illustrated in \cite[Fig.~7]{Kuba_16}.
Namely, consider a binary free multilabeled increasing tree $T$
on the label set $[n]$.
We define an increasing restricted ternary tree $T'$ on the vertex set $[n]$
as follows:
Replace each vertex $u$ in $T$, having label set $L_u = \{u_1< \ldots< u_i\}$,
by a chain of vertices $u_1- u_2 - \cdots - u_i$ in $T'$,
where the node $u_j$ has $u_{j+1}$ as its middle child for $1\leq j\leq i-1$. 
If $u$ has a left (resp.~right) child $v$ in $T$,
then in $T'$ the final vertex $u_i$
has the initial vertex $v_1$ as its left (resp.~right) child.
The reverse bijection can be obtained by simply contracting the middle edges
in a restricted ternary tree to a single multilabeled vertex.

\begin{example}
\rm
Figure~\ref{fig.bmt} is an example of a
binary free multilabeled increasing tree on the label set $[6]$.
This tree is in bijective correspondence with the
increasing restricted ternary tree shown in Figure~\ref{fig.rt}.

\begin{figure}[t!]
\centering
\newcommand{\nodea}{\node[draw,circle] (a) {$1$}
;}\newcommand{\nodeb}{\node[draw,circle] (b) {$2,4$}
;}\newcommand{\nodec}{\node[draw,circle] (c) {$3$}
;}\newcommand{\noded}{\node[draw,circle] (d) {$5$}
;}\newcommand{\nodee}{\node[draw,circle] (e) {$6$}
;}
\begin{tikzpicture}
\matrix[column sep=.1cm, row sep=0.5cm,ampersand replacement=\&]{
\& \& \nodea \& \\
\& \nodeb \& \& \nodec \\
\&  \& \noded \& \\
\& \& \& \nodee\\
};
        \draw[ultra thick] (a) -- (b) node {};
        \draw[ultra thick] (a) -- (c) node {};
        \draw[ultra thick] (c) -- (d) node {};
        \draw[ultra thick] (d) -- (e) node {};
\end{tikzpicture}
  \vspace*{-2mm}
\caption{An example of a binary free multilabeled increasing tree on the label set $[6]$.}
  \vspace*{1cm}
\label{fig.bmt}
\end{figure}
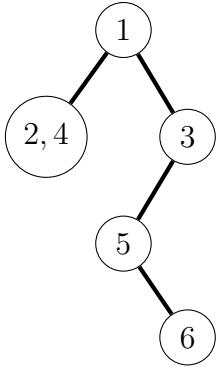
\myendremark
\end{example}

\section*{Acknowledgments}

A crucial role was played in this work by
the On-Line Encyclopedia of Integer Sequences \cite{OEIS}.
We warmly thank Neil Sloane for founding this indispensable resource,
and the hundreds of volunteers for helping to maintain and expand it.

The first author (V.B.)\ was an undergraduate at University College London
when this work was performed;
she thanks the Department of Mathematics
for two undergraduate research bursaries.

The second author (B.D.)\ was supported in 2023--24
by DIMERS project ANR-18-CE40-0033 funded by Agence Nationale de la Recherche
(ANR, France).

We also thank Anders Claesson for drawing our attention to
the work of Vajnovszki \cite{Vajnovszki_18} on vincular patterns.

\appendix

\section{Table of OEIS matches for some T-fractions with quasi-affine coefficients}
\label{sec.OEIS}

In order to come up with our conjectures,
we used a reverse-engineering approach.
That is, rather than starting from a combinatorial family
and attempting to find a T-fraction enumerating it,
we started instead from a ``nice'' T-fraction
and attempted to find a combinatorial interpretation for it.
More precisely, we considered a generic quasi-affine T-fraction of the form
\reff{eq.alphadelta.quasiaffine}
in which we let each of the variables $x,y,u,v,a,b,c,d$
be either 0 or 1.
For each of the $2^8 = 256$ cases we generated the first 10~terms $a_n$
using a {\sc Mathematica} code,
and we then searched the OEIS \cite{OEIS} for the given sequence,
deleting the initial term~$a_0 = 1$.
Of course, the lack of an OEIS entry does not necessarily mean
that the sequence is combinatorially uninteresting;
indeed, the sequence \reff{eq.c=0} is a counterexample.
And conversely, an OEIS entry sometimes lacks a combinatorial interpretation.
But we figured that the OEIS would be a good place to start.
Our idea was that if the OEIS entry for an {\em integer}\/ sequence
arising from $x,y,u,v,a,b,c,d \in \{0,1\}$
mentions some combinatorial interpretation,
then by refining that combinatorial interpretation
we might conjecture (and then prove) a T-fraction
in which one or more of the coefficients are variables
conjugate to some statistic in the combinatorial model.
This strategy turned out to pay off.

In reality we imposed some restrictions.
To begin with, we imposed $x=y=1$,
because taking either $x=0$ or $y=0$
would make $\alpha_1 = 0$ or $\alpha_2 = 0$,
leading to a finite (and quite trivial) continued fraction;
this restriction reduces the number of cases to $2^6 = 64$.
Likewise, we disallowed $a=b=c=d=0$, since that would give an S-fraction.
We excluded those because the main goal of this project
was to discover new T-fractions that {\em cannot}\/ be represented 
as an S-fraction or J-fraction.
This further brought down the number of cases to 60.


We performed an automated search for these 60 sequences on the OEIS.
Among the resulting searches, 
we further omitted the sequences with $c=d=0$,
since they correspond to simple linear transforms of S-fractions\footnote{
   See \cite[Propositions~3 and 15]{Barry_09} for some special cases,
   and \cite{Sokal_totalpos} for the general case.
};
this reduced the number of sequences to 48.
Our {\sc Mathematica} code returned the following OEIS matches:
\begin{itemize}
    \item A187251 --- On setting 
        $x=y=v=c=1$ and $u=a=b=d=0$
        in~\reff{eq.alphadelta.quasiaffine}
        we obtain a sequence that starts as
        \be
        1,1,2,6,22,94,460,2532,15420,102620,
        739512,\ldots .
        \ee
        This sequence matches the OEIS entry A187251,
        which has the description 
        ``Number of permutations of $[n]$ having no cycle with 3 or more alternating runs (it is assumed that the smallest element of a cycle is in the first position).''
    
    \item A105072 --- On setting 
        $x=y=v=a=c=1$ and $u=b=d=0$
        in~\reff{eq.alphadelta.quasiaffine}
        we obtain a sequence that starts as 
        \be
        1,2,5,16,63,290,1511,8756,55761,386394,
        2889181,\ldots .
        \ee
        This sequence matches the OEIS entry A105072,
        which has the description 
        ``Number of permutations on $[n]$ whose local maxima are in ascending order.''
    
    \item A230008 --- On setting 
        $x=y=u=v=b=d=1$ and $a=c=0$
        in~\reff{eq.alphadelta.quasiaffine}
        we obtain a sequence that starts as 
        \be
        1,1,3,11,51,295,2055,16715,155355,1624255,
        18868575,\ldots
        \ee
        This sequence [which is \reff{eq.a=c=0}]
        matches the OEIS entry A230008,
        which has as a comment (by Markus Kuba)
        ``Counts binary free multilabeled increasing trees with $m$ labels.''
        This comment was the starting point for the present research.
        In this paper we have shown in Theorem~\ref{thm.Simple.J.abt},
        and independently in Section~\ref{subsubsec.free},
        that this sequence also counts increasing restricted ternary trees:
        $a_n = |\RT_n|$.
\end{itemize}




After finishing the bulk of this research,
we performed yet another automated search,
this time including also the value $2$ in addition to $0$ and $1$.
This gives $3^8 = 6561$ cases.
We imposed $x,y\in \{1,2\}$  to avoid a finite continued fraction;
this reduced the number of cases to $4\times 3^6 = 2946$.
Then we disallowed $a=c=0$ and $b=d=0$,
as either would give us a J-fraction by using the
odd and even contraction formulae, respectively.
The resulting number of cases is further reduced to 2304.
Our {\sc Mathematica} code returned 13 OEIS matches;
we list these in Table~\ref{tab.OEIS.search}.

\begin{table}[t]
    \centering
    \Scale[0.77]{
    \begin{tabular}{c|l|c|l}
        OEIS A number & First few terms & 
        $(x,y,u,v,a,b,c,d)$ & Description\\
        \hline
          A258173 & $1, 1, 3, 12, 58, 321, 1975, 13265$ & (1,1,0,0,0,1,2,2) &
               Sum over all Dyck paths of semilength $n$ \\
          &&&  of products over all peaks $p$ of $y_p$, where \\
          &&& $y_p$ is the $y$-coordinate of peak $p$. 
          \\[2mm]
          A006318 & 1, 2, 6, 22, 90, 394, 1806, 8558&(1,1,0,0,1,1,0,0)  & Large Schr\"oder numbers\\[2mm]
          A302285 & 1, 2, 7, 33, 185, 1170, 8121 &(1,1,0,0,1,2,2,2) & No combinatorial description  \\[2mm]
          A047891 & 1, 3, 12, 57, 300, 1686, 9912 &(1,1,0,0,2,2,0,0) & Number of planar rooted trees \\
          &&&
          with $n$ nodes and tricolored end nodes. \\[2mm]
          A155866 & 1, 2, 6, 22, 91, 413, 2032 &(1,1,0,1,1,1,0,0) & 
          A `Morgan-Voyce' transform\\
          &&&
          of the Bell numbers
          \\[2mm]
          A155857 & 1, 2, 6, 23, 107, 590, 3786 &(1,1,1,1,1,1,0,0) &  
          Row sums of triangle \\
          &&&
          A155856 
          ($\binom {2 n - k} {k} (n - k)!$)
          \\[2mm]
          A000311 & 1, 1, 4, 26, 236, 2752, 39208 &(1,2,2,2,0,1,2,2) &  
          Schr\"oder's fourth problem;\\
          &&&
          also series-reduced rooted trees\\
          &&&
          with $n$ labeled leaves;\\
          &&&
          also number of total partitions 
            of $n$.\\[2mm]
          A001515 & 1, 2, 7, 37, 266, 2431, 27007 &(1,2,2,2,1,1,0,0) &  Bessel polynomial $y_n(x)$ \\
          &&&
          evaluated at $x=1$. \\[2mm]
          A006351 & 1, 2, 8, 52, 472, 5504, 78416 &(1,2,2,2,1,2,2,2) & 
          Number of series-parallel networks\\
          &&&
          with $n$ labeled edges.\\
          &&&
          Also called yoke-chains\\
          &&&
          by Cayley and MacMahon.
          \\[2mm]
          A043301 & 1, 3, 13, 77, 591, 5627, 64261 &(1,2,2,2,2,2,0,0) & 
          $a(n) = 2^n\sum_{k = 0}^{n} (n + k)!/((n - k)!k!4^k).$ \\[2mm]
          A155867 & 1, 3, 13, 65, 355, 2061, 12501 &(2,1,0,0,1,1,0,0) &  
          A `Morgan-Voyce' transform \\
          &&&
          of the large Schr\"oder numbers A006318.
          \\[2mm]
          A103210 & 1, 3, 15, 93, 645, 4791, 37275 &(2,2,0,0,1,1,0,0) &  $a(0)=1$\\
          &&&
          $a(n) = (1/n)\sum_{i = 0}^{n - 1} 
          \binom{n}{i} \binom{n}{i + 1} 2^i 3^{n - i}$\\[2mm]
          A156017 &  1, 4, 24, 176, 1440, 12608 &(2,2,0,0,2,2,0,0) &  Schr\"oder paths with two rise colors\\
          &&&
          and two level colors.
    \end{tabular}
    }
    \caption{OEIS entries having various T-fractions.}
    \vspace*{1cm}
    \label{tab.OEIS.search}
\end{table}




\bigskip

As an aid to researchers who may wish to implement similar
reverse-engineering searches
--- possibly with a very large number of test sequences ---
we share our {\sc Mathematica} code for automated searches of the OEIS:

\begin{verbatim}
llToString[ll_] := TextString[ll, ListFormat -> {"", ",", ""}]

OEISquerystring[ll_] := 
   "http://oeis.org/search?q=" <> llToString[ll] <> "&fmt=json"

getOEISJSON[ll_] := Import[OEISquerystring[ll], "JSON"]

ToAForm[num_] :=
   Module[{str1 = ToString[num]},
          If[StringLength[str1] >= 6,
             Return["A" <> str1],
             Return["A" <> StringRepeat["0", 6 - StringLength[str1]]
                        <> str1]
            ]
         ]

SequenceInOEIS[ll_, delete_:1, verbose_: False] :=
   Module[{res = getOEISJSON[Drop[ll, delete]]},
          If[verbose, Print[Length[res], " results"]];
          If[res === Null,
             Return[{}],
             Return[{Map[ToAForm,  ("number" /. res)], ll}]
            ]
         ]   
\end{verbatim}

\section{Context-free grammars (= derivative operators) for our generating polynomials}
\label{sec.cfg.operator}

In this appendix, 
we show how the sequences enumerating our families of trees can be
obtained in a very simple way using context-free (Chen) grammars \cite{Chen_93} (see also \cite{Dumont_96}),
also known as derivative operators.
This approach was very helpful to us in guessing our families of trees,
and also in helping us to check the correctness of our T-fractions.
Indeed, whenever we conjectured a family of trees
to match a given sequence arising from a T-fraction, 
we checked our conjecture by first writing out a 
derivative operator to generate that family of trees;
this was a quick way of obtaining a weighted count of the trees in our family
without having to construct them explicitly.
We then checked whether the weighted count of the trees
matched the sequence generated by our T-fraction.

We applied this approach to increasing binary trees
and increasing restricted ternary trees.
However, for increasing interval-labeled restricted ternary trees,
our approach was different:
instead, we generated the first few terms
of the ordinary generating function for increasing restricted ternary trees
using a derivative operator,
and then checked our construction by using the identity~\reff{eq.ogf.rtt.irtt}.

In the rest of this section, we write out our grammar rules
and derivative operators for binary trees and restricted ternary trees.

\subsection{Increasing binary trees}

We shall construct a differential operator that generates the polynomials
$P_n(x_1,x_2,y_1,y_2)$ defined in \reff{eq.Pn.bt},
which enumerate increasing binary trees with weights $x_1, y_1, x_2, y_2$
for the node types $11$, $00$, $10$, $01$, respectively.

To do this, we start with an increasing binary tree $T$
on the vertex set $[n-1]$,
and we consider the various ways in which a new vertex with label $n$
can be attached to this tree as a child of some vertex $j \in [n-1]$:
\begin{itemize}
    \item If $j$ is a leaf in $T$, with weight $y_1$,
    then $n$ can be attached as either a left child or a right child of $j$.
    In either case $n$ has node type $00$,
    and $j$ gets node type $10$ or $01$, respectively.
    This gives the grammar rule 
    $y_1\mapsto y_1 x_2 + y_1 y_2 $.

    \item If $j$ has node type 10 (resp.~01) in $T$,
    with weight $x_2$ (resp.~$y_2$),
    then $n$ can be attached as a right (resp.~left) child of $j$.
    After $n$ is attached, $j$ has node type $11$
    and $n$ has node type $00$.
    This gives us the grammar rules
    $\{x_2 \mapsto x_1 y_1 
    \,,\,
    y_2 \mapsto x_1 y_1\}$.

    \item If $j$ has node type 11 in $T$, then $n$ cannot be attached there.
\end{itemize}
Combining this information, we define the differential operator
\be
\scrd
\;\eqdef\;
y_1(x_2+y_2)\,\dfrac{\partial}{\partial y_1}
\:+\:
x_1 y_1\left(\dfrac{\partial}{\partial x_2} +
             \dfrac{\partial}{\partial y_2} \right)\;.
\label{eq.delta.bt}
\ee
Since $P_1 = y_1$, we have proven the following proposition:

\begin{proposition}
The polynomials $P_n(x_1, x_2, y_1, y_2)$ defined in~\reff{eq.Pn.bt}
satisfy
\be
   P_n(x_1, x_2, y_1, y_2)
   \;=\;
   \scrd^{n-1} y_1
   \quad\hbox{for $n \ge 1$}
\ee
where $\scrd$ is defined by \reff{eq.delta.bt}.
\end{proposition}

{\bf Remark.}
Specializing our grammar rule 
$
\{ {y_1\mapsto y_1 x_2 + y_1 y_2}
\,,\,
{x_2 \mapsto x_1 y_1}
\,,\,
{y_2 \mapsto x_1 y_1}
\}
$
to $x_1 = 1, y_1 = x$, $x_2 = y_2 =y$
gives the grammar rule $\{ {x \mapsto 2xy}\,,\, {y \mapsto x} \}$
studied in \cite[section~2.4]{Dumont_96}.
\myendremark

\subsection{Increasing restricted ternary trees}

\begin{sloppypar}
We shall construct a differential operator that generates the polynomials
$P_n(x_1,x_2,y_1,y_2,w)$ defined in \reff{eq.Pn.abt},
which enumerate increasing restricted ternary trees
with weights $x_1, y_1, x_2, y_2, w$
for the node types $101$, $000$, $100$, $001$, $010$, respectively.
\end{sloppypar}

To do this, we start with an increasing restricted ternary tree $T$
on the vertex set $[n-1]$,
and we consider the various ways in which a new vertex with label $n$
can be attached to this tree.
The reasoning is the same as for increasing binary trees, except that:
\begin{itemize}
    \item If $j$ is a leaf in $T$, then $n$ can also be attached
       as a middle child of $j$.
       This gives the grammar rule $y_1\mapsto y_1 (x_2 + y_2 +w)$.
    \item If $j$ has node type 010 in $T$, then $n$ cannot be attached there.
\end{itemize}
We therefore define the differential operator
\be
\scrd
\;\eqdef\;
y_1(x_2+y_2+w)\,\dfrac{\partial}{\partial y_1}
\:+\:
x_1 y_1\left(\dfrac{\partial}{\partial x_2} +
             \dfrac{\partial}{\partial y_2} \right)\;.
\label{eq.delta.rt}
\ee
Since $P_1 = y_1$, we have proven the following proposition:

\begin{proposition}
The polynomials $P_n(x_1, x_2, y_1, y_2, w)$ defined in~\reff{eq.Pn.abt}
satisfy
\be
   P_n(x_1, x_2, y_1, y_2, w)
   \;=\;
   \scrd^{n-1} y_1
   \quad\hbox{for $n \ge 1$}
\ee
where $\scrd$ is defined by \reff{eq.delta.rt}.
\end{proposition}

\end{document}